\documentclass[reqno, 10pt]{amsart}
\usepackage{tikz-cd}
\usepackage{enumitem}
\usepackage{amscd}
\usepackage{mathtools}
\usepackage{amssymb}
\usepackage{MnSymbol}
\usepackage[all]{xy}
\usepackage{fge}
\usepackage{hyperref}
\usepackage{graphics}

\usepackage[textsize=tiny]{todonotes}

\RequirePackage[T1]{fontenc}
\RequirePackage{amsfonts,latexsym,amssymb}

\let\cal\mathcal   

\let\goth\mathfrak

\usepackage{bbm}
\usepackage{hyperref}

\newtheorem{theorem}[equation]{Theorem}
 \newtheorem{lemma}[equation]{Lemma}
 \newtheorem{proposition}[equation]{Proposition}
 \newtheorem{corollary}[equation]{Corollary}

\theoremstyle{definition}

\newtheorem{remark}[equation]{Remark}

\newtheorem{hypothesis}[equation]{Hypothesis}
\newtheorem{example}[equation]{Example}
\theoremstyle{remark}
\newtheorem*{acknowledgments}{Acknowledgments}

\def\jcdot{\scriptscriptstyle\bullet}

\def\invlim{\mathop{\vtop{\ialign{##\crcr$\hfill{\lim}\hfil$\crcr
\noalign{\kern1pt\nointerlineskip}\leftarrowfill\crcr\noalign
{\kern -3pt}}}}\limits}
\def\dirlim{\mathop{\vtop{\ialign{##\crcr$\hfill{\lim}\hfil$\crcr
\noalign{\kern1pt\nointerlineskip}\rightarrowfill\crcr\noalign
{\kern -3pt}}}}\limits} 
\def\lomapr#1{\smash{\mathop{\relbar\joinrel\longrightarrow}\limits^{#1}}}

\def\phi{\varphi} 
\def\epsilon{\varepsilon}

\let\goth\mathfrak
\newcommand{\ovk}{\overline{K} }

\newcommand{\dr}{\operatorname{dR} } 
\newcommand{\pdr}{\operatorname{pdR} }

      \newcommand{\nr}{\operatorname{nr} }   
 
 \newcommand{\colim}{\operatorname{colim} }

 \newcommand{\proeet}{\operatorname{pro\acute{e}t} } 
 \newcommand{\eet}{\operatorname{\acute{e}t} }
  \newcommand{\Spa}{\operatorname{Spa} }

 \newcommand{\Hom}{{\rm{Hom}} }
  \newcommand{\uHom}{\underline{{\rm{Hom}} }}
  \newcommand{\Hhom}{{\cal Hom}} 
 \newcommand{\Ext}{\operatorname{Ext} }

\newcommand{\Gal}{\operatorname{Gal} }
\newcommand{\pFF}{\operatorname{pFF} }

\newcommand{\st}{\operatorname{st} }

 \newcommand{\crr}{\operatorname{cr} }
  \newcommand{\hk}{\operatorname{HK} }
   \newcommand{\LL}{\operatorname{L} }
     \newcommand{\FF}{\operatorname{FF} }

 \newcommand{\sff}{{\mathcal{F}}}

 \newcommand{\sh}{{\mathcal{H}}}
 \newcommand{\sg}{{\mathcal{G}}}
 \newcommand{\stt}{{\mathcal{T}}}

 \newcommand{\scc}{{\mathcal{C}}}

 \newcommand{\so}{{\mathcal O}}
 
 \newcommand{\se}{{\mathcal{E}}}

\newcommand{\sd}{{\mathcal{D}}} 
 
\newcommand{\spp}{{\mathcal{P}}}

\newcommand{\bA}{{\mathbb{A}}}

\newcommand{\cO}{{\mathcal{O}}}
\newcommand{\an}{{\mathrm{an}}}
\newcommand{\alg}{{\mathrm{alg}}}
\newcommand{\tr}{{\mathrm{tr}}}

 \newcommand{\wt}{\widetilde}
 \newcommand{\wh}{\widehat}
   \numberwithin{equation}{section}

\def\R{{\mathrm R}}
\def\O{{\cal O}}
  \def\B{{\bf B}}
\def\Q{{\bf Q}} \def\Z{{\bf Z}}

\def\N{{\bf N}}
\def\O{{\cal O}}
 
\def\Qbar{\overline{\bf Q}}
\def\sur{\twoheadrightarrow}

\def\rg{{\rm R}\Gamma}

\def\epsilon{\varepsilon}
\def\dual{{\boldsymbol *}}

\numberwithin{equation}{section}
\begin{document}
\title[$C_{\st}$-conjecture]
 {$C_{\st}$-conjecture for compactly supported pro-\'etale cohomology}
 \author{Pierre Colmez} 
\address{CNRS, IMJ-PRG, Sorbonne Universit\'e, 4 place Jussieu, 75005 Paris, France}
\email{pierre.colmez@imj-prg.fr} 
\author{Sally Gilles} 
\address{LMB, Universit\'e Marie et Louis Pasteur, 16 route de Gray, 25030 Besan\c{c}on CEDEX, France}
\email{sally.gilles@umlp.fr}
 \author{Zhenghui Li} 
\address{Institute for Advanced Study, 1 Einstein Drive, Princeton, N.J. 08540, U.S.A.}
\email{zhenghuili@ias.edu}
\author{Wies{\l}awa Nizio{\l}}
\address{CNRS, IMJ-PRG, Sorbonne Universit\'e, 4 place Jussieu, 75005 Paris, France}
\email{wieslawa.niziol@imj-prg.fr}
 \date{\today}
\thanks{The research of P.C. and W.N.    was supported in part by the  Simons Collaboration on Perfection in Algebra, Geometry, and Topology. This material is based upon work supported by the National Science Foundation under Grant No. DMS-2424139, while W.N.  was in residence at the Simons Laufer Mathematical Sciences Institute in Berkeley, California, during the Fall 2026 semester. }
\maketitle
 \begin{abstract}
We prove the \'etale-de Rham direction of the $C_{\st}$-conjecture for compactly supported 
cohomology of smooth partially proper rigid analytic varieties over local fields of mixed characteristic.  More precisely, we recover filtered de Rham cohomology and
(under a finite-rank assumption) Hyodo--Kato cohomology
 from compactly supported $p$-adic pro-\'etale cohomology 
as derived Hom (in a category of Arithmetic Topological Vector Spaces)
with values in $\B_{\pdr}$ and $\B_{\pFF}$ respectively
(period rings related to $\B_{\dr}$ and $\B_{\st}$ whose higher Galois cohomology vanishes).

The proof uses, via the geometrized basic comparison theorem, Poincar\'e  duality for Hyodo--Kato and (filtered) de Rham cohomology,   explicit extension calculations for period sheaves, and fully-faithfulness of the passage from Galois equivariant perfect complexes on the Fargues--Fontaine curve  to 
Arithmetic Topological Vector Spaces.
 \end{abstract}
\tableofcontents
 \section{Introduction}

Let \(p\) be a prime, let \(K\) be a finite extension of $\Q_p$,  and let \(C\) be the
\(p\)-adic completion of an algebraic closure of \(K\). We set
\(G_K=\operatorname{Gal}(\overline K/K)\).

\subsection{Motivation}
  We will review briefly what motivated our paper. 

\subsubsection{Algebraic  setting: Conjecture $C_{\st}$} 
Comparison  theorems  of \(p\)-adic Hodge theory of algebraic varieties state  that the
de Rham plus Hyodo--Kato cohomologies and the $p$-adic \'etale cohomology   of an algebraic  variety can be recovered from each other 
 after extension of scalars to suitable
period rings. Known as the 
 \(C_{\mathrm{dR}}\)- and \(C_{\mathrm{st}}\)-comparison
theorems\footnote{We use the formulation of Beilinson from \cite{Beil}.} they construct, for an algebraic variety $X$ over $K$,  functorial, compatible under the base change $(-)\otimes^{\LL}_{\B_{\st}}\B_{\dr}$, quasi-isomorphisms
\begin{align}\label{utah1}
& (C_{\st}):\quad \R\Gamma_{\proeet}(X_C,\Q_p)\otimes^{\LL}_{\Q_p}\B_{\st}  \simeq \R\Gamma_{\hk}(X_C)\otimes^{\LL}_{C^{\nr}}\B_{\st},\\
& (C_{\dr}):\quad \R\Gamma_{\proeet}(X_C,\Q_p)\otimes^{\LL}_{\Q_p}\B_{\dr}  \simeq \R\Gamma_{\dr}(X)\otimes^{\LL}_{K}\B_{\dr}.\notag
\end{align}
These quasi-isomorphisms   are compatible with all the underlying structures: Frobenius, monodromy, and Galois action for the first one and filtration and Galois action for the second one.

\begin{remark}  
  (1) ({\em \'Etale-to-de Rham.}) We can  recover de Rham and Hyodo--Kato cohomologies from \'etale cohomology   by taking ``fixed" points of  
the suitable operators in \eqref{utah1}. We get the isomorphisms\footnote{``sm'' stands for
``smooth'', i.e. fixed by an open subgroup of $G_K$.}: 
\begin{align}\label{utah2}
(C_{\dr}):\quad H^i_{\dr}(X) & \simeq (H^i_{\eet}(X_C,\Q_p)\otimes_{\Q_p}\B_{\dr})^{G_K} \quad \quad + F^{\jcdot},\\
(C_{\st}):\quad H^i_{\hk}(X_C) & \simeq (H^i_{\eet}(X_C,\Q_p)\otimes_{\Q_p}\B_{\st})^{G_K-\mathrm{sm}} \quad + (\phi,N,G_K).\notag
\end{align}

(2) ({\em De Rham-to-\'etale version.}) A slightly more refined version of comparison theorems allows  to recover \'etale cohomology from de Rham and Hyodo--Kato cohomologies using just the "+" 
period rings\footnote{If $B_?$ is a period ring, the ``non +'' version is obtained
from the ``+'' version by inverting $t$: $B_?=B_?^+[\frac{1}{t}]$, where $t$ is
Fontaine's $2\pi i$.} :
 For  $i\leq r$, we have the exact sequence
of Galois modules
$$
0\to H^i_{\eet}(X_C,\Q_p(r))\to (H^i_{\hk}(X_C)\otimes_{C^{\nr}}\B^+_{\st})^{N=0,\phi=p^r}\to (H^i_{\dr}(X)\otimes_K\B^+_{\dr})/F^r\to 0
$$
\end{remark}

   There is also a  derived version of isomorphisms \ref{utah2} but 
for that we need  almost period rings: analogs of $\B_{\dr}$ and $\B_{\st}$ that are Galois acyclic. 
 For de Rham cohomology, recall first that we have 
 $$
 H^i(G_K,C(j))\simeq \begin{cases} K & \mbox{ if j=0, i=0,1},\\
 0 & \mbox{ otherwise.} 
 \end{cases}
$$
This is inherited by Galois cohomology of $\B_{\dr}$ yielding a nontrivial class in  $H^1$
(given by $\sigma\mapsto \log \chi_{\mathrm{cycl}}(\sigma)$). 
To kill this extra class, 
one introduces, after Fontaine~\cite{FonA}, the ring 
$\B_{\pdr}:=\B_{\dr}[\log t]$ with $\log t$ 
transcendental over $\B_{\dr}$ and satisfying (as in~\cite{Csen})
$\sigma(\log t)= \log t +\log \chi_{\mathrm{cycl}}(\sigma)$, for $\sigma\in G_K$,  yielding 
$\B_{\pdr}^{\R G_K}:=\R\Gamma(G_K,\B_{\pdr})\simeq K$.

For Hyodo--Kato cohomology (see \cite{CGN3} for details), 
we set $\B:=\so(Y_{\FF})$, where $Y_{\FF}$ is the Fargues-Fontaine punctured open disc,  
$\B^+_{\FF}:=\B[\log[p^{\flat}]],$  and $\B^+_{\pFF}$ -- the unipotent envelope of $\B$: the first period ring is a Fargues-Fontaine analog of $\B_{\crr}^+$, the second one -- of $\B^+_{\st}$, and the third one -- of the envisioned $\B^+_{\mathrm{pst}}$. We have
$\B_{\pFF}^{\R G_K-\mathrm{sm}}\simeq C^{\nr}$, where $C^{\nr}=\colim_{[L:K]<\infty}\mathrm{Frac(W(k_L))}$.

Having that, we easily derive from \eqref{utah1}:
\begin{align}\label{gateway2}
(C_{\dr}):\quad \R\Gamma_{\dr}(X) & \simeq (\R\Gamma_{\eet}(X_C,\Q_p)\otimes_{\Q_p}\B_{\pdr})^{\R G_K} \quad\quad + F^{\jcdot},\\
(C_{\st}):\quad \R\Gamma_{\hk}(X_C) & \simeq (\R\Gamma_{\eet}(X_C,\Q_p)\otimes_{\Q_p}\B_{\pFF})^{\R G_K-\mathrm{sm}} \quad + (\phi,N,G_K).\notag
\end{align}

\subsubsection{Analytic setting: Conjecture $C_{\st}$ for usual cohomology}  
For rigid analytic varieties, the corresponding comparison theorems exist mostly for 
pro-\'etale
cohomology (which is usually infinite-dimensional). See \cite{CN5} for a detailed exposition and \cite{Tate} for a survey. 

    Let $X$ be a  smooth,  partially proper rigid analytic variety over $K$. For $i\leq r$, the basic comparision theorem yields an exact sequence  of solid $\Q_p$-modules equipped with a Galois action
\begin{equation}\label{utah3}
\mathrm{HK}^{r,i-1}\to \mathrm{DR}^{r,i-1}\to H^i_{\proeet}(X_C,\Q_p(r))\to {\rm HK}^{r, i}\to\mathrm{DR}^{r,i},
\end{equation}where we set
\begin{align*}\mathrm{HK}^{r,j} & =(H^j_{\hk}(X_C)\otimes^{\Box}_{C^{\nr}}\B^+_{\st})^{N=0,\phi=p^r},\\
\mathrm{DR}^{r,j} & =H^j((\R\Gamma_{\dr}(X)\otimes^{{\Box}}_K\B^+_{\dr})/F^r).
\end{align*}
We like to think of \eqref{utah3} as analog of the fundamental exact sequence ($r\geq 0$)
$$
0\to \Q_p(r)\to (\B_{\st}^{+})^{N=0,\phi=p^r}\to\B^+_{\dr}/F^r\to 0
$$
\begin{example}
(1) If $X$ is proper then $H^i_{\proeet}(X_C,\Q_p(r))$ is of finite rank over $\Q_p$
which forces the exact the sequence \eqref{utah3} to split. 
As a result we obtain a short exact sequence
$$
0\to H^i_{\proeet}(X_C,\Q_p(r))\to  \mathrm{HK}^{r, i}\to\mathrm{DR}^{r,i}\to 0,
$$
which is exactly the same as in the algebraic setting. 

(2)  If $X$ is Stein the sequence and $i=r$ \eqref{utah3} yields the top exact sequence in the following fundamental diagram
$$
\xymatrix@C=5mm@R=4mm{
0\ar[r] & \Omega^{r-1}(X_C)/\ker d\ar[r] \ar@{=}[d]& H^r_{\proeet}(X_C,\Q_p(r))\ar[r] \ar[d] & (H^r_{\hk}(X_C)\otimes^{\Box}_{C^{\nr}}\B^+_{\st})^{N=0,\phi=p^r}\ar[r] \ar[d] & 0 \\
0\ar[r] & \Omega^{r-1}(X_C)/\ker d\ar[r] & \Omega^r(X_C)^{d=0}\ar[r] & H^r_{\dr}(X_C)\ar[r] & 0
}
$$
This diagram represents a map of exact sequences; as a result the right square is bicartesian allowing to recover pro-\'etale cohomology from Hyodo--Kato cohomology and pieces of de Rham complex. 

  (3) In cases (1) and (2) pro-\'etale cohomology can be recovered from de Rham data. This has a uniform formulation for smooth partially proper rigid analytic varieties  and holds  for proper and Stein varieties,  for analytifications of algebraic varieties, for almost proper varieties, and for products of a smooth proper variety with a smooth Stein variety \cite{CN5}. However it fails in general (see \cite{CN6} for a counterexample).
\end{example}

   For usual pro-\'etale cohomology, if $X$ is Stein of dimension $d$, the fundamental diagram allows us
to  recover de Rham and Hyodo--Kato
cohomology by taking continuous morphisms into the de Rham and
semistable period rings. We get in $\sd(K_{\Box})$ and $\sd(C^{\nr}_{\Box})$, respectively
\footnote{We underline a set to indicate that the set comes with a condensed structure
in case this is not obvious.}:
\begin{align}\label{utah4}
\uHom_{G_K}(H^n_{\proeet}(X_C,\Q_p(d)), \B_{\dr}) & \simeq H^{2d-n}_{\dr,c}(X) \quad \quad + F^{\jcdot},\\
\uHom_{G_K-\mathrm{sm}}(H^n_{\proeet}(X_C,\Q_p(d)), \B_{\st}) & \simeq H^{2d-n}_{\hk,c}(X_C) \quad + (\phi,N,G_K).\notag
\end{align}
We do not yet  have a derived version of this result\footnote{The reader may wonder why the same argument as presented in this paper for compactly supported cohomology  would not work and the main reason for this is that we crucially use the fact that the derived dual of a space of compact type is equal to the naive dual and the same does not hold for nuclear Fr\'echet spaces.}.  

\subsubsection{Analytic setting: Conjecture $C_{\st}$ for compactly supported  cohomology} 
For compactly supported pro-\'etale cohomology of smooth partially proper varieties  we still have a  basic comparison theorem yielding an analog of the exact sequence \eqref{utah3}. This sequence does not however split in the Stein case so we do not have an analog of the  fundamental diagram. See \cite{AGN} for details.
\begin{example} Let  $X$ be  a smooth Stein curve over $K$. Early heuristic computations in yet to be defined category of {\em Arithmetic Topological Vector Spaces} (ATVS's for short) gave us  (see \cite{CGN3}) that in $\sd(K_{\Box})$
 \begin{align*}
\underline{\Hom}_{\mathrm{ATVS}}(\mathbb{H}_{\proeet,c}^1(X_C,\Q_p(1)),\mathbb{B}_{\pdr}) &  \simeq H^1_{\rm dR}(X)^{(1)}
\subset H^1_{\rm dR}(X),\\
\underline{\Ext}^1_{\mathrm{ATVS}}(\mathbb{H}_{\proeet,c}^2(X_C,\Q_p(1)),\mathbb{B}_{\pdr}) & \simeq {\rm Ker}\big(H^1_{\rm dR}(X)\to H^1_{\rm dR}(X)^{(1)}\big),
\end{align*}
where $H^1_{\rm dR}(X)^{(1)}$ denotes the  subspace of de Rham cohomology with Frobenius slopes $1$. So, unlike the case of usual cohomology, $\underline{\Hom}_{\mathrm{ATVS}}(\mathbb{H}_{\proeet, c}^1(X_C,\Q_p(1)),\mathbb{B}_{\pdr})$ recovers only a part of de Rham cohomology; hence an analog of \eqref{utah4}  is not possible. 

On the other hand, the above computation  suggests that,   since we should   have the exact sequence 
\begin{align*}
0\to \underline{\rm Ext}^1_{\mathrm{ATVS}}(\mathbb{H}^{2}_{\proeet,c}(X_C,\Q_p(1)),\mathbb{B}_{\pdr}) & \to H^1(\R\underline{\Hom}_{\mathrm{ATVS}}(\mathbb{R}_{\proeet,c}(X_C,\Q_p(1))[2],\mathbb{B}_{\pdr})\\
& \to \underline{\rm Hom}_{\mathrm{ATVS}}(\mathbb{H}^{1}_{\proeet,c}(X_C,\Q_p(1)),\mathbb{B}_{\pdr})\to 0,
\end{align*}
the right object to look at is
$\R\underline{\Hom}_{\mathrm{ATVS}}(\mathbb{R}_{\proeet,c}(X_C,\Q_p(1))[2],\mathbb{B}_{\pdr})$:  we expect to  have $$H^1(\R\underline{\Hom}_{\mathrm{ATVS}}(\mathbb{R}_{\proeet,c}(X_C,\Q_p(1))[2],\mathbb{B}_{\pdr})\simeq H^1_{\dr}(X).$$

We note that if we had used $\mathbb{B}_{\dr}$ instead of $\mathbb{B}_{\pdr}$ we would have added a parasitic term to $\underline{\Ext}^1_{\mathrm{ATVS}}(\mathbb{H}_{\proeet,c}^2(X_C,\Q_p(1)),\mathbb{B}_{\pdr}) $. 
\end{example}

\subsection{Main theorem}  
The   problem of what to choose for the category of ATVS's turned out to be a delicate one. It was tempting to  just work in the category of solid $G_K$-modules (in analogy with \eqref{utah4}). However, as we show in the appendix, $\Ext^1(C,C)$ in this category is not isomorphic to $K$ as one would expect\footnote{Following the computations of Fontaine in the category of almost $C_p$-representations (see \cite{FoC}).}  but is actually infinite dimensional over $K$. 
We settled in the end on a
\(G_K\)-equivariant version of Topological Vector Spaces\footnote{Topologically enriched presheaves on strictly totally disconnected  perfectoid spaces over \(C\).} defined in \cite{TVS} and used in geometric duality statements of \cite{CGN2}.
The compactly supported pro-\'etale cohomology of \(X_C\) defines an
object $
\mathbb R_{\mathrm{pro\acute et},c}
   (X_C,\mathbf Q_p) $
of this category.

   Our main result proves the \'etale-to-de Rham direction of the compact-support \(C_{\mathrm{st}}\)-conjecture
in the setting of ATVS's.

\begin{theorem}[Conjecture $C_{\st}$: \'etale-to-de Rham] 
\label{main-intro}
Let \(X\) be a smooth partially proper rigid analytic variety over
\(K\), of dimension \(d\).

\begin{enumerate}
\item
There is a natural quasi-isomorphism\footnote{We defined the filtered derived category:  
$
\sff^{\mathrm{aug}}\sd(K_\square)
:=
\mathrm{Fun}\bigl(
(\mathbf Z\cup\{-\infty\})^{\rm op},
\sd(K_\square)
\bigr). $}  in
\(\sff^{\mathrm{aug}}\sd(K_\square)\)
\[
F^\bullet \R\Gamma_{\mathrm{dR}}(X)
\simeq
\R{\uHom}_{\mathrm{ATVS}}
\bigl(
\mathbb{R}_{\mathrm{pro\acute et},c}
  (X_C,\mathbf Q_p(d))[2d],
F^\bullet \mathbb{B}_{\mathrm{pdR}}
\bigr).
\]
Here the right-hand side denotes the filtered object whose
\(a\)-th term is
\[
\R\uHom_{\mathrm{ATVS}}
\bigl(
\mathbb{R}_{\mathrm{pro\acute et},c}
  (X_C,\mathbf Q_p(d))[2d],
F^a\mathbb{B}_{\mathrm{pdR}}
\bigr),
\qquad
a\in\mathbf Z\cup\{-\infty\},
\]
where we set $ F^{-\infty}\R\Gamma_{\mathrm{dR}}(X):= \R\Gamma_{\mathrm{dR}}(X) $ and \(F^{-\infty}\mathbb B_{\mathrm{pdR}}
:=\mathbb B_{\mathrm{pdR}}\).

\item
Assume that the de Rham cohomology of \(X\) has finite rank. Then
there is a natural quasi-isomorphism
\[
\R\Gamma_{\mathrm{HK}}(X_C)
\simeq
\R\underline{\operatorname{Hom}}^{\mathrm{sm}}_{\mathrm{ATVS}}
\left(
 \mathbb R_{\mathrm{pro\acute et},c}
     (X_C,\mathbf Q_p(d))[2d],
 \mathbb B_{\mathrm{pFF}}
\right)
\]
in the derived category of solid
\((\varphi,N,G_K)\)-modules over
\(C^{\,\mathrm{nr}}\).
\end{enumerate}
\end{theorem}

\begin{remark}[Conjecture $C_{\st}$: \'etale-to-Hodge] (1) Let $a\in\Z$. 
Passing to associated
graded objects in claim (1) of Theorem \ref{main-intro} 
yields a natural quasi-isomorphism in $\sd(K_\Box)$
$$
\mathrm{gr}_F^a\R\Gamma_{\mathrm{dR}}(X)
\simeq
\R\underline{{\Hom}}_{\mathrm{ATVS}}
\left(
 \mathbb{R}_{\proeet,c}
     (X_C,\mathbf Q_p(d))[2d],
\mathbb{G}_a(a)[\log t]
\right).
$$
If $X$ is a  Stein space, this becomes a natural quasi-isomorphism in $\sd(K_\Box)$
\begin{equation}\label{gateway3}
\Omega^a(X)[-a]
\simeq
\R\underline{{\Hom}}_{\mathrm{ATVS}}
\left(  \mathbb R_{\mathrm{\proeet,c}}
     (X_C,\mathbf Q_p(d))[2d],
 \mathbb{G}_a(a)[\log t] \right).
\end{equation}
Similarly (Corollary~\ref{gateway1}) one can also recover 
the differential $d: \Omega^a\to \Omega^{a+1}$. 
Moreover (see section~\ref{soso12}), 
conditionally on compatibility of the trace in the six-functor formalism with the geometric trace, a multiplicative refinement recovers the de Rham differential graded algebra.

  (2)  Let $X$ be a smooth algebraic variety over $K$ compactified  to a smooth proper variety $Y$ via a normal crossing divisor $D$. Assume that $X$ is affine. Then the classical conjecture $C_{\st}$ (see \eqref{gateway2}) yields an isomorphism 
 $$
\R\Gamma(Y,  \Omega^a(\log D))[-a]\simeq  \R\Gamma(G_K,\R\Gamma_{\eet}(X_C)\otimes_{\Q_p}^{\LL}C(a)[\log t])
 $$
 recovering the logarithmic differentials from \'etale cohomology:
 $$
 H^0(Y,  \Omega^a(\log D))\simeq  H^a\R\Gamma(G_K,\R\Gamma_{\eet}(X_C)\otimes_{\Q_p}^{\LL}C(a)[\log t]).
 $$
On the other hand, starting from the pro-\'etale cohomology of $X^{\mathrm{an}}$,
quasi-isomorphism \eqref{gateway3} allows to recover the whole of
$\Omega^a(X^{\rm an})$ (bigger than  $\Omega^a(X)$ as 
differential forms on $X^{\rm an}$ can have an essential singularity at the boundary). 
\end{remark}

\subsection{Outline of the proof}
The proof has two main steps: Conjecture $C_{\st}$ on the Fargues-Fontaine curve and a descent to Arithmetic Topological Vector Spaces. 

\subsubsection{Conjecture $C_{\st}$ on the Fargues-Fontaine curve} 
The proof takes place first on the Fargues--Fontaine curve. In previous
work (see \cite{CGN2}), we  associated  to \(X\) three \(G_K\)-equivariant solid
quasi-coherent complexes on the Fargues-Fontaine curve $X_{\FF}$:
\begin{equation*}
\mathcal E_{\mathrm{pro\acute et},c}
(X_C,\mathbf Q_p(r)),
\qquad
\mathcal E_{\mathrm{HK},c}(X_C,r),
\qquad
\mathcal E_{\mathrm{dR},c}(X_C,r),
\end{equation*}
whose coherent cohomology are: 
\begin{align*}
\R\Gamma(X_{\FF},\mathcal E_{\mathrm{pro\acute et},c}
   (X_C,\mathbf Q_p(r))) & \simeq  \R\Gamma_{\proeet,c}(X_C,\Q_p(r)),\\
\R\Gamma(X_{\FF}, \mathcal E_{\mathrm{HK},c}(X_C,r) )& \simeq [\R\Gamma_{\hk,c}\otimes_{C^{\nr}}^{\LL_{\Box}}\B^+_{\st}]^{N=0,\phi=p^r},\\
\R\Gamma(X_{\FF},\mathcal E_{\mathrm{dR},c}(X_C,r) )& \simeq (\R\Gamma_{\dr,c}(X)\otimes^{\LL_\Box}_{K}\B_{\dr}^+)/F^r.
\end{align*}
The complex $\mathcal E_{\mathrm{dR},c}(X_C,r)$ is supported at the point at infinity $x_\infty$
(defined by the equation $t=0$) of $X_{\FF}$, and is $t$-torsion. 
For $r\geq 2d$,   these complexes are related  by the  basic comparison triangle
\begin{equation}\label{utah5}
\mathcal E_{\mathrm{pro\acute et},c}
   (X_C,\mathbf Q_p(r))
\longrightarrow
\mathcal E_{\mathrm{HK},c}(X_C,r)
\longrightarrow
\mathcal E_{\mathrm{dR},c}(X_C,r).
\end{equation}

   For the following geometric calculations, assume first that \(X\) is Stein with finite-rank de Rham cohomology. We dualize this triangle against two different coefficient objects.
For the de Rham cohomology, the coefficient object is
$ i_{\infty,*}F^a\B_{\mathrm{dR}}$,
 for $i_{\infty}: x_{\infty}\hookrightarrow X_{\FF}$. 
Hyodo--Kato Poincar\'e duality, filtered de Rham Poincar\'e duality, and
explicit calculations with skyscraper sheaves give a natural
quasi-isomorphism
\begin{equation}\label{utah8}
 \R\mathcal Hom_{\mathrm{FF}}
\left(
 \mathcal E_{\mathrm{pro\acute et},c}
     (X_C,\mathbf Q_p(d))[2d],
 i_{\infty,*}F^a\B_{\mathrm{dR}}
\right)\simeq F^a\mathcal E_{\mathrm{dR}}(X),
\end{equation}
where $F^a\mathcal E_{\mathrm{dR}}(X):=i_{\infty,*}F^a(\R\Gamma_{\dr}(X)\otimes_K^{\LL_{\Box}}\B_{\dr})$. 
 The computation  is particularly simple in the case of $a=-\infty$ because then the dual of the de Rham term (i.e., $\mathcal E_{\mathrm{dR},c}(X_C,r)$) 
from \eqref{utah5} is trivial (since $t$ is invertible in $\B_{\dr}$) 
and we just need to dualize the Hyodo--Kato term. 

For the Hyodo--Kato cohomology, morally speaking\footnote{The precise definition has to take into account the fact that the Fargues-Fontaine punctured disc $Y$ is Stein.}, we use the coefficient object
$
\mathcal{B}_{\pFF}:=\R\pi_{{\FF},*}\mathcal O_{\mathrm{pFF}},
$
where
\(\pi_{\FF}:Y_{\mathrm{FF}}\to X_{\mathrm{FF}}\)
is the usual projection from the \(Y\)-curve. Since \(t\) is
invertible in \(\mathcal O_{\mathrm{pFF}}\), the de Rham term of the
basic comparison triangle has zero dual. Hyodo--Kato Poincar\'e duality
then yields
\begin{equation}\label{utah7}
\R\mathcal Hom_{\mathrm{FF}}
\left(
 \mathcal E_{\mathrm{pro\acute et},c}
     (X_C,\mathbf Q_p(d))[2d],
\mathcal{B}_{\pFF}
\right)
\simeq
\R\pi_{{\FF},*}
\left(
 \R\Gamma_{\mathrm{HK}}(X_C)
 \otimes_{ C^{\,\mathrm{nr}}}^{\LL}
 \mathcal O_{\mathrm{pFF}}
\right).
\end{equation}
Taking Galois invariants of \eqref{utah8} (extended to the period rings $\B_{\pdr}$) and \eqref{utah7}  we derive: 
\begin{theorem}[Conjecture $C_{\st}$ on the Fargues-Fontaine curve] \label{utah6} Let $X$ be a smooth partially proper rigid analytic variety over $K$. 
\begin{enumerate}
\item  Let $a\in \Z\cup \{-\infty\}$. We have in $\sd(K_{\Box})$
$$
\R\Hom_{\FF,G_K}(\se_{\proeet,c}(X_C,\Q_p(d))[2d],i_{\infty,*}F^{a}\B_{\pdr})\simeq F^{a}\R\Gamma_{\dr}(X).
$$
\item Assume that the de Rham cohomology of \(X\) has finite rank. We have in $\sd_{\phi,N,G_K}(C^{\nr}_{\Box})$
$$
\R\Hom_{\FF,G_K-\mathrm{sm}}(\se_{\proeet,c}(X_C,\Q_p(d))[2d],\mathcal{B}_{\pFF})\simeq \R\Gamma_{\hk}(X_C).
$$
\end{enumerate}
\end{theorem}

\subsubsection{Descent to ATVS} 
A substantial part of the paper is devoted to transporting these
curve-level identities to ATVS's via the equivariant global sections functor 
$$
\R\tau_*: \quad \mathrm{\mathrm{QCoh}}(X_{\FF})^{G_K}\to \mathrm{ATVS}.
$$
We do it actually on the level of the global section  functor 
$$
\R\tau_*: \quad \mathrm{\mathrm{QCoh}}(X_{\FF})\to \mathrm{TVS}.
$$
and then impose  Galois equivariance. 
 The existing  fully-faithfulness results apply directly only  to perfect
complexes (see \cite{TVS}), whereas the coefficient objects
\(i_{\infty,*}F^a\B_{\mathrm{pdR}}\) and \(\mathcal{B}_{\mathrm{pFF}}\) are not perfect. Hence we 
 need to  establish several complementary fully-faithfulness
statements.

     A key ingredient is the calculation of derived Homs in VS's from
\(\mathbb G_a\) into period sheaves. Since the relevant period
objects involve colimits, these calculations cannot be obtained
formally from the Banach--Colmez case. We use instead the functorial
Mac Lane resolution of \(\mathbb G_a\). Its terms are represented by
affine spaces over perfectoid bases, and the required Ext-vanishing is
reduced to the affine-space comparison theorem for period sheaves.
This gives, among other formulas, that
\[
\R{\Hhom}_{\mathrm{VS}} (\mathbb G_a,\mathbb{B}_{\mathrm{dR}}) \simeq 0,
\qquad
\R{\Hhom}_{\mathrm{VS}} (\mathbb G_a,\mathbb{B}[1/t]) \simeq 0,
\]
together with their \(\mathbb{B}_{\mathrm{pdR}}\)- and
\(\mathbb{B}_{\mathrm{pFF}}\)-analogues and the filtered computation
\[
\R{\Hhom}_{\mathrm{VS}} (\mathbb G_a,F^a\mathbb{B}_{\mathrm{pdR}})
\simeq \mathrm{gr}^{a-1}\mathbb{B}_{\mathrm{pdR}}[-1], \quad a\in\Z.
\]

   These computations extend the relevant fully-faithfulness statements to the compactly supported pro-\'etale complexes and the period coefficients used here. Taking equivariant global sections then transports the curve-level reconstruction formulas of Theorem \ref{utah6} to the ATVS reconstruction formulas of Theorem \ref{main-intro}.
\subsubsection*{Organization of the paper}

In Section~2, we compute derived Homs between the period sheaves that
occur in the proof. We first perform the corresponding calculations
on the Fargues--Fontaine curve, then pass to Vector Spaces.  We  use  Mac Lane resolutions to treat noncompact
period sheaves.

Section~3 contains the geometric heart of the proof. We dualize the
basic comparison triangle on the Fargues--Fontaine curve and prove the
de Rham and Hyodo--Kato reconstruction formulas in this setting.

In Section~4, we introduce Arithmetic Topological Vector Spaces and discuss
the equivariant categories of solid quasi-coherent sheaves on the
Fargues--Fontaine curve. We study equivariant internal Homs, global
sections, and the fully-faithfulness properties of the global section
functors.

In Section~5, we descend the curve-level formulas to ATVS's and prove
Theorem~\ref{main-intro}. We also spell out the resulting Hodge
cohomology formula and illustrate the calculation in the case of
smooth Stein curves.

   In Section~6, we
apply the preceding computations to revisit  
the geometric and arithmetic dualities for $p$-adic pro-\'etale cohomology of  
rigid analytic varieties from \cite{CGN2} and \cite{Li}.

    Finally, in  Appendix \ref{firstA},
 we show that there is a discrepancy between $\Ext^1$-groups in almost $C_p$-representations and Galois-equivariant solid $\Q_p$-modules. This   discrepancy is caused by non-almost-split extensions: Fontaine's category retains only extensions which
can be split be quotienting by a finite-dimensional $ \Q_p$-representation, 
whereas the ambient Banach category contains extension classes 
of unbounded complexity along the cyclotomic tower. 
In Appendix \ref{second1} we  construct the Fr\'echet avatar \(\widetilde \B_{\mathrm{pdR}}\) 
of $\B_{\dr}$ and compute its Galois cohomology.
\begin{acknowledgments}  
 We would  like to thank Johannes Ansch\"utz, Joseph Ayoub, Shizhang Li, Emanuel Reinecke, Bogdan Zavyalov   for  helpful discussions concerning the content of this paper.
  Z.L.  would like to thank Institute for Advanced Study and Simons Foundation for their support.  
 
   The authors used GPT-5.6 and GPT-6 (OpenAI)  during the preparation of this manuscript for mathematical proofreading, checks of  internal consistency and cross-references, and exploratory discussion of several arguments.  
\end{acknowledgments}

 \subsubsection*{Notation and conventions.}\label{Notation} 
Let $p$ be a fixed prime, $\Q_p$ -- the $p$-adic rational numbers. 
 Let $K$ be a finite extension of $\Q_p$ with residue field $k$. Let $F=\mathrm{Frac}(W(k))$. Let $\ovk$ be an algebraic closure of $K$ and let $C=\wh{\ovk}$ be the $p$-adic completion of $\ovk$.  Let $G_K$ be the Galois group of $K$. 
We will denote by $ \B_{\crr}, \B_{\st},\B_{\dr}$ the crystalline, semistable, and  de Rham period rings of Fontaine; we will denote by $ \mathbb{B}_{\crr}, \mathbb{B}_{\st}, \mathbb{B}_{\dr} $ the corresponding period sheaves on the $v$-site of $\mathrm{Spa}(C)$. 

 The category of affinoid perfectoid spaces over an affinoid perfectoid space $S$ over $C$ will be denoted by  ${\rm Perf}_S$. The pro-\'etale site of $S$ will be  the category  ${\rm Perf}_S$ endowed with the quasi-pro-\'etale topology and will be denoted by $S_{\proeet}$.   
All rigid analytic spaces considered will be over K or C. We assume that they are separated,
taut, and countable at infinity.

   For $S\in {\rm Perf}_C$, we will denote by $X_{\FF,S^{\flat}}$  
(the relative) Fargues-Fontaine curve over $S^{\flat}$; by 
    $\pi_{\FF}: Y_{\FF,S^{\flat}}\to X_{\FF,S^{\flat}}$ -- the canonical 
projection (we have $X_{\FF,S^{\flat}}=Y_{\FF,S^{\flat}}/\varphi^\Z$). 
If $S={\rm Spa}(C)$ we will drop the subscript $S^{\flat}$ from the notation. For an analytic  adic space $X$, we write     ${\rm QCoh}(X)$, ${\rm Nuc}(X)$, ${\rm Perf}(X)$  
for  the categories of solid quasi-coherent, nuclear, and perfect  sheaves on $X$, respectively (see  \cite{And21}).  
    
    We denote by  ${\rm Cond}$, ${\rm CondAb}$, ${\rm Solid}$ the condensed sets, condensed abelian groups, and solid abelian groups, respectively. To simplify arguments, we fix a cut-off cardinal $\kappa$ and assume that  perfectoid spaces and condensed sets are $\kappa$-small. 

 \section{Vector Spaces}
The purpose of this section is to compute $\R\Hhom$'s in Vector Spaces between Banach-Colmez spaces and certain period sheaves.

 \subsection{Extensions of quasi-coherent sheaves on the Fargues-Fontaine curve} \label{better2} 
We start with computing  $\R\Hhom$'s  between certain quasi-coherent sheaves on the Fargues-Fontaine curve that are induced from period rings. 

    \subsubsection{Absolute period rings} \label{period1}
We will briefly recall  the definitions of  period rings $\B_{\pdr}$ and $\B_{\pFF}$  from \cite{FonA} and  \cite{CGN3}. 

We start with $\B_{\pdr}$.
Let $\B^+_{\pdr}$ be the algebra $\B^+_{\dr}[\log t]$ defined in \cite[Sec. 4.3]{FonA} and let $\B_{\pdr}:=\B^+_{\pdr}[1/t]$.
The Galois group $G_K$ acts on $\B^+_{\pdr}$ via ring homomorphisms extending its usual action on $\B^+_{\dr}$ and such that $ g (\log t)=\log t+\log(\chi_{\rm cycl}( g ))$, $ g \in G_K$.
This action naturally extends to $\B_{\pdr}$.
The filtration on $\B_{\pdr}$ is given by $F^a\B_{\pdr}=t^a\B_{\pdr}^+$.
We have 
$$\R\Gamma(G_K, \B^+_{\pdr})\simeq \R\Gamma(G_K, \B_{\pdr})\simeq K[0].  $$

 For two integers $a\leq b$, let $\B^{[a,b[}_{\dr}:=F^a\B_{\dr}/F^b\B_{\dr}$ and similarly for $\B_{\pdr}$. For $i\in\Z$, we have a $G_K$-equivariant isomorphism $\B^{[a-i,b-i[}_{\dr}\simeq\B^{[a,b[}_{\dr}\otimes\Q_p(-i)$. 
 
   We pass now to  the period rings $\B_{\FF}$ and $\B_{\pFF}$.  Define  $U=\log[p^{\flat}]$ as a transcendental element over  $\B$, with 
$\varphi(U)=p\,U$ and $g(U)=U+ \log [g(p^{\flat})/p^{\flat}].$
We set:
$$\B_{\rm FF}^+:=\B[U],\quad\B_{\rm FF}:=\B_{\rm FF}^+[\tfrac{1}{t}]. $$
The ring
${\B}^+_{\rm FF}$ is endowed with a surjective monodromy operator $N=-\frac{d}{dU}$, commuting with the action of $G_K$, such that 
$N\varphi=p\varphi N$ and ${\B}= {\B}_{\rm FF}^{+,N=0}$. We extend the canonical  injection
$\B\hookrightarrow\B_{\dr}^+$ to a ring morphism 
$\B_{\rm FF}^+\otimes_FK\to\B_{\dr}^+$ by sending $U$ to $\log\frac{[p^{\flat}]}{p}$.
This morphism is  $G_K$-equivariant and  injective. 

Let us recall the main steps of the construction of $\B_{\pFF}^+$ 
(see~\cite[\S\,3.3]{CGN3} for details).
Take any open interval $I\subset]0,\infty[$ with endpoints in $\Q^\dual$
(or choose an increasing sequence of such, covering $]0,\infty[$),
 let $Y_{\FF,I}$ be the corresponding open of $Y_{\FF}$, and $\B_I:=\O(Y_{\FF,I})$.
For each finite extension \(L/K\) in \(\overline K\), let
\(\B_{I,\mathrm{uni},L}\) denote the universal unipotent
\(G_L\)-algebra over \(\B_{I}\) constructed in \cite{CGN3},
by successively adjoining primitives of Galois 1-cocycles; 
denote by $\B^{I,\leq m}_{\mathrm{uni},L}$ its sub-module of elements of weight~$\leq m$.
These are $(\phi,N,G_L)$-modules over $C^{\nr}$ (but $N\varphi\neq p\varphi N$),
and the natural action of $g\in G_K$ sends $\B^{I,\leq m}_{\mathrm{uni},L}$ to
$\B^{I,\leq m}_{\mathrm{uni},gL}$. These structures can be made compatible with 
restrictions $Y_{\FF}^J\to Y_{\FF}^I$ (see~\cite[\S\,3.2]{CGN3}),
and finite field extensions $K\subset L\subset M$
(see~\cite[Lemme 3.17]{CGN3}).
Set
\begin{align*}
 & \B_{\mathrm{pFF},L}^{+,\leq m}
:=
\lim_I \B_{\mathrm{uni},L}^{I,\leq m},\quad 
\B^+_{\mathrm{pFF},L}
:=
\colim_m \B^{+,\leq m}_{\pFF,L},\\
 & \B_{\mathrm{pFF}}^{+}
:=
\operatorname*{colim}_{K\subset L\subset\overline K,\,[L:K]<\infty}
\B_{\mathrm{pFF},L}^{+}.
\end{align*}
We set 
$ \B_{\mathrm{pFF},L}^{\leq m}:=  \B_{\mathrm{pFF},L}^{+,\leq m}[1/t]$, 
$\B_{\mathrm{pFF},L}:=\B_{\mathrm{pFF},L}^{+}[1/t]$, 
$\B_{\mathrm{pFF}}:=\B_{\mathrm{pFF}}^{+}[1/t].$
The $(\phi, N,G_K)$-structure extends from $\B^+_{\rm pFF}$ to $\B_{\rm pFF}$. 
   We choose a $G_K$-equivariant and injective morphism 
$\B_{\rm pFF}\otimes_FK\to\B_{\rm pdR}$ .
Moreover, we have~\cite[Th.\,3.19]{CGN3}:
 $$\R\Gamma(G_K, \B^+_{\rm pFF})\simeq \R\Gamma(G_K, \B_{\rm pFF})\simeq F[0].
   $$

\subsubsection{Relative period rings}\label{relative1}
   We start with the de Rham period rings. Let $S\in\operatorname{Perf}_C$.  Let
$\mathbf{B}_{\mathrm{dR}}^+(S)$ and
$\mathbf{B}_{\mathrm{dR}}(S)=\mathbf{B}_{\mathrm{dR}}^+(S)[1/t]$
be the usual relative de Rham period rings, with their canonical
solid structures. Write $\ell=\log t$ for a formal indeterminate.
For $m\geq 0$, set
\[
\mathbf{B}_{\mathrm{pdR}}^{+,\leq m}(S)
 :=\bigoplus_{j=0}^{m}\mathbf{B}_{\mathrm{dR}}^+(S)\,\ell^j
\]
and define
\begin{align*}
\mathbf{B}_{\mathrm{pdR}}^+(S)
 &:=\colim_m\mathbf{B}_{\mathrm{pdR}}^{+,\leq m}(S)
   =\mathbf{B}_{\mathrm{dR}}^+(S)[\ell],\\
\mathbf{B}_{\mathrm{pdR}}(S)
 &:=\mathbf{B}_{\mathrm{pdR}}^+(S)[1/t]
   =\colim_{r,m}t^{-r}
       \mathbf{B}_{\mathrm{pdR}}^{+,\leq m}(S)
   =\mathbf{B}_{\mathrm{dR}}(S)[\ell].
\end{align*}
For a morphism $S'\to S$ over $C$, pullback acts on the coefficients
and fixes $\ell$. In particular, the construction is functorial in $S$
and recovers the absolute period rings at $S=\operatorname{Spa}(C)$.

 For the filtration, for $a\in\mathbf{Z}$, set
$$
F^a\mathbf{B}_{\mathrm{pdR}}(S)
 :=t^a\mathbf{B}_{\mathrm{pdR}}^+(S)
 =\colim_m\bigoplus_{j=0}^{m}
      F^a\mathbf{B}_{\mathrm{dR}}(S)\,\ell^j.
$$
For $a\leq b$, this gives
\[
\mathbf{B}_{\mathrm{pdR}}^{[a,b[}(S)
 :=F^a\mathbf{B}_{\mathrm{pdR}}(S)
      /F^b\mathbf{B}_{\mathrm{pdR}}(S)
 =\mathbf{B}_{\mathrm{dR}}^{[a,b[}(S)[\ell].
\]
We put
$F^{-\infty}\mathbf{B}_{\mathrm{pdR}}(S)
 :=\mathbf{B}_{\mathrm{pdR}}(S)$.

Let $\mathbb{B}_{\mathrm{dR}}^+$ denote the usual de Rham period
sheaf on $\operatorname{Perf}_{C,v}$. We define
$$
\mathbb{B}_{\mathrm{pdR}}^+
 :=\colim_m
       \left(\bigoplus_{j=0}^{m}\mathbb{B}_{\mathrm{dR}}^+\,\ell^j\right),\quad 
\mathbb{B}_{\mathrm{pdR}}:=\mathbb{B}_{\mathrm{pdR}}^+[1/t],\quad
F^a\mathbb{B}_{\mathrm{pdR}}
 :=t^a\mathbb{B}_{\mathrm{pdR}}^+.
$$
 Equivalently, these are the sheaves
associated to the preceding functorial constructions. Their values
on an affinoid perfectoid $S$ are the rings and filtered modules
defined above.

     We pass now to the $\B_{\pFF}$ period rings.  
Let $S\in\operatorname{Perf}_C$. If $I\subset]0,\infty[$ is an open interval with
endpoints in $\Q^\dual$, write
\[
\mathbf{B}_{I}(S)
 :=\Gamma(Y_{\mathrm{FF},S^\flat,I},\mathcal{O}).
\]
The structural morphism $S\to\operatorname{Spa}(C)$ induces a map
$\mathbf{B}_{I}\to\mathbf{B}_{I}(S)$. For $m\geq 0$ and a finite
extension $L/K$ in $\overline K$, define
\begin{equation}\label{eve1}
\mathbf{B}^{I,\leq m}_{\mathrm{uni},L}(S)
 :=\mathbf{B}_{I}(S)
   \otimes^{\LL_{\Box}}_{\mathbf{B}_{I}}
   \mathbf{B}^{I,\leq m}_{\mathrm{uni},L}.
\end{equation}
Since $\mathbf{B}^{I,\leq m}_{\mathrm{uni},L}$ is finite free over
$\mathbf{B}_{I}$, this tensor product is concentrated in degree $0$.

The restriction maps in $I$, the inclusions in $m$, and the
finite-extension maps in $L$ are induced by the corresponding maps
of the fixed absolute system. These constructions are functorial
in $S$. Multiplication induces maps
\[
\mathbf{B}^{I,\leq m}_{\mathrm{uni},L}(S)
 \otimes^{\LL_{\Box}}_{\mathbf{B}_{I}(S)}
\mathbf{B}^{I,\leq m'}_{\mathrm{uni},L}(S)
 \longrightarrow
\mathbf{B}^{I,\leq m+m'}_{\mathrm{uni},L}(S).
\]

  We set
\begin{align*}
\mathbf{B}^{+,\leq m}_{\mathrm{pFF},L}(S) &:=\lim_I\mathbf{B}^{I,\leq m}_{\mathrm{uni},L}(S),
\quad 
\mathbf{B}^{+}_{\mathrm{pFF},L}(S) :=\colim_m\mathbf{B}^{+,\leq m}_{\mathrm{pFF},L}(S),\\
\mathbf{B}^{+}_{\mathrm{pFF}}(S) &:=\colim_{\substack{K\subseteq L\subset\overline K,[L:K]<\infty}}\mathbf{B}^{+}_{\pFF,L}(S).
\end{align*}
The products above give the last two objects their structures of
commutative $\mathbf{B}(S)$-algebras (recall that 
$\mathbf{B}(S):=\lim_I\mathbf{B}_{I}(S)$). Finally, define
\[
\begin{aligned}
\mathbf{B}^{\leq m}_{\mathrm{pFF},L}(S)
 &:=\mathbf{B}^{+,\leq m}_{\mathrm{pFF},L}(S)[1/t],\quad
\mathbf{B}_{\mathrm{pFF},L}(S)
 :=\mathbf{B}^{+}_{\mathrm{pFF},L}(S)[1/t],\qquad
\mathbf{B}_{\mathrm{pFF}}(S)
 :=\mathbf{B}^{+}_{\mathrm{pFF}}(S)[1/t].
\end{aligned}
\]
For $S=\operatorname{Spa}(C)$, these formulas recover the preceding
absolute definitions.

For the corresponding period sheaves, start with the finite-level
sheaves
\[
\mathbb{B}^{I,\leq m}_{\mathrm{uni},L}:= \{S\mapsto {\B}^{I,\leq m}_{\mathrm{uni},L}(S)\}
\]
 and apply
the same formulas in the category of $v$-sheaves.
The chosen operators $\phi, N$ are extended by scalar extension
at finite level and then through the above limits.

\subsubsection{Extensions of de Rham sheaves on the Fargues-Fontaine curve} 
Let $S\in \mathrm{Perf}_C$. 
  Let $(-)(i)$  be the twisting by the line bundle $\so(1)$  on the category ${\rm QCoh}(X_{\FF, S^{\flat}})$. We will often  use the notation $\so(n,k):=\so(n-k)\otimes\Q_p(k)$. For $n,k\geq 0$, we have the Galois equivariant\footnote{For arbitrary \(S/C\), Galois equivariance is understood functorially with respect to Galois twists of \(S\), as in Section \ref{relative1}.} exact sequence
$$
0\to \so(n,k)\to \so(n)\to i_{\infty,*}(\so/t^k)\to 0 
$$
Thus, we can identify
$$
i_{\infty,*}(\mathbb{B}^+_{\dr}(S)/t^k)\simeq[\so(0,k)\to \so][1]
$$
and, more generally, for $b\geq a$, 
we have
$$
i_{\infty,*}\mathbb{B}^{[a,b[}_{\dr}(S)\simeq [\so(a,b)\to\so(a,a)][1];
$$
as Galois equivariant sheaves. For $i\in\Z$, we have a Galois equivariant isomorphism 
$$i_{\infty,*}\mathbb{B}^{[a-i,b-i[}_{\dr}(S)\simeq i_{\infty,*}\mathbb{B}^{[a,b[}_{\dr}(S)\otimes^{\LL}_{\so}\so(i)\otimes\Q_p(-i)$$
We used here that tensoring a skyscraper sheaf with a line bundle is a neutral operation.

We set $$\R\Hhom_{\FF,S^{\flat}}(-,-):=\R\Hhom_{{\rm QCoh}(X_{\FF,S^{\flat}})}(-,-), \quad \R\Hom_{\FF,S^{\flat}}(-,-):=\R\Hom_{{\rm QCoh}(X_{\FF,S^{\flat}})}(-,-).
$$

\begin{lemma}\label{paris3a}
Let $i\geq 0$ and $a\leq b\in\Z$.  
\begin{enumerate}
\item We have in ${\rm QCoh}(X_{\FF,S^{\flat}})$ 
\begin{align*}
\R\Hhom_{\FF,S^{\flat}}(i_{\infty,*}\mathbb{B}^{[0,i[}_{\dr}(S), i_{\infty,*}\mathbb{B}^{[a,b[}_{\dr}(S))  &  \simeq   [i_{\infty,*}\mathbb{B}^{[a,b[}_{\dr}(S)\to i_{\infty,*}\mathbb{B}^{[a-i,b-i[}_{\dr}(S)].
\end{align*}
In particular
\begin{align*}
\Hhom_{\FF,S^{\flat}}(i_{\infty,*}\mathbb{B}^{[0,i[}_{\dr}(S), i_{\infty,*}\mathbb{B}^{[a,b[}_{\dr}(S))  & \simeq i_{\infty,*}\mathbb{B}^{[\max(a,b-i),b[}_{\dr}(S),\\
 {\mathcal Ext}^1_{\FF,S^{\flat}}(i_{\infty,*}\mathbb{B}^{[0,i[}_{\dr}(S), i_{\infty,*}\mathbb{B}^{[a,b[}_{\dr}(S))  & \simeq i_{\infty,*}\mathbb{B}^{[a-i,\min(a,b-i)[}_{\dr}(S).
\end{align*}
\item We have in ${\rm QCoh}(X_{\FF,S^{\flat}})$ 
$$
\R\Hhom_{\FF,S^{\flat}}(i_{\infty,*}\mathbb{B}^{[0,i[}_{\dr}(S), i_{\infty,*}F^a\mathbb{B}_{\dr}(S)) \simeq  i_{\infty,*}\mathbb{B}^{[a-i,a[}_{\dr}(S)[-1].
$$
\end{enumerate}
We have analogous results with $\mathbb{B}_{\pdr}$ in place of $\mathbb{B}_{\dr}$ in the target. 
\end{lemma}
\begin{proof}  We compute in ${\rm QCoh}(X_{\FF,S^{\flat}})$ 
(with $\R\Hhom:=\R\Hhom_{\FF,S^{\flat}}$, $\iota:=i_{\infty,*}$, 
${\rm B}_{\dr}^?:=\mathbb{B}^?_{\dr}(S)$):
\begin{align*}
\R\Hhom(\iota \mathrm{B}^{[0,i[}_{\dr}, \iota \mathrm{B}^{[a,b[}_{\dr}) & \simeq   \R\Hhom([\so(0,i)\to\so][1],\iota \mathrm{B}^{[a,b[}_{\dr})\\
& \simeq [\R\Hhom(\so,\iota \mathrm{B}^{[a,b[}_{\dr})\to \R\Hhom(\so(0,i)[1],\iota \mathrm{B}^{[a,b[}_{\dr})]\\
& \simeq [\iota \mathrm{B}^{[a,b[}_{\dr}\to \R\Hhom(\so(0,i),[\so(a,b)\to \so(a,a)])]\\
& \simeq [\iota \mathrm{B}^{[a,b[}_{\dr}\to [\so(a,b-i)\to \so(a,a-i)])]\\
& \simeq [\iota \mathrm{B}^{[a,b[}_{\dr}\to \iota \mathrm{B}^{[a,b[}_{\dr}\otimes^{\LL}_{\so}\so(i)\otimes\Q_p(-i))]
\simeq [\iota \mathrm{B}^{[a,b[}_{\dr}\to \iota \mathrm{B}^{[a-i,b-i[}_{\dr}].
\end{align*}
This proves the first quasi-isomorphism of the lemma.

For the second one:  it follows  from the first one that, for $a\in\Z$,  we have 
\begin{align}\label{paris6}
\R\Hhom(\iota \mathrm{B}^{[0,i[}_{\dr}, \iota F^a\mathrm{B}_{\dr})  & \simeq  \R\lim_b \R\Hhom(\iota \mathrm{B}^{[0,i[}_{\dr}, \iota \mathrm{B}^{[a,b[}_{\dr})\\
& \simeq  [ \R\lim_b \iota \mathrm{B}^{[a,b[}_{\dr}\to \R \lim_b\iota \mathrm{B}^{[a-i,b-i[}_{\dr}]\notag\\
    & \simeq   [\iota F^a\mathrm{B}_{\dr}\to \iota F^{a-i}\mathrm{B}_{\dr}]
    \simeq \iota \mathrm{B}^{[a-i,a[}_{\dr}[-1],\notag
\end{align}
as wanted. 

The above implies the following computation: 
\begin{align*}
\R\Hhom(\iota \mathrm{B}^{[0,i[}_{\dr}, \iota F^a\mathrm{B}_{\pdr})  & \simeq    \R\Hhom(\iota \mathrm{B}^{[0,i[}_{\dr}, \colim_m\bigoplus_{j=0}^{m}
\iota F^a\mathrm B_{\mathrm{dR}}\,\ell^j) \\
& \simeq    \colim_m\bigoplus_{j=0}^{m}\R\Hhom(\iota \mathrm{B}^{[0,i[}_{\dr}, \iota 
F^a\mathrm B_{\mathrm{dR}})\ell^j \\
& \simeq   \iota \colim_m\bigoplus_{j=0}^{m}\mathrm{B}^{[a-i,a[}_{\dr}\ell^j[-1]
\simeq \iota \mathrm{B}^{[a-i,a[}_{\pdr}[-1],
\end{align*}
as wanted.
   \end{proof}
\begin{remark} \label{special} We note that, for $a\leq b-i$, the above lemma yields in  ${\rm QCoh}(X_{\FF,S^{\flat}})$ 
\begin{align*}
\Hhom_{\FF,S^{\flat}}(i_{\infty,*}\mathbb{B}^{[0,i[}_{\dr}(S), i_{\infty,*}\mathbb{B}^{[a,b[}_{\dr}(S))  &   \simeq i_{\infty,*}\mathbb{B}^{[b-i,b[}_{\dr}(S),\\
 {\mathcal Ext}^1_{\FF,S^{\flat}}(i_{\infty,*}\mathbb{B}^{[0,i[}_{\dr}(S), i_{\infty,*}\mathbb{B}^{[a,b[}_{\dr}(S))   &  \simeq i_{\infty,*}\mathbb{B}^{[a-i,a[}_{\dr}(S).
 \end{align*}
\end{remark}

\subsubsection{Extensions of Hyodo--Kato sheaves on the Fargues-Fontaine curve} 
If $I\subset]0,+\infty[$ is an open interval with endpoints in $\Q^\dual$, 
set $Y_I:=Y_{\FF,I}$, and let $j_I:Y_I\to Y_{\FF}$ be the corresponding open immersion.
Let \(\spp_{I}\in\mathrm{QCoh}(Y_I)\) denote the solid 
quasi-coherent algebra on \(Y_I\) associated to 
 \(\B^{I}_{\mathrm{uni},K}\); 
we denote by $\spp_{I}^{\leq m}\in \mathrm{QCoh}(Y_I)$ 
the weight~$\leq m$ part of $\spp_{I}$.  We define
\[
\mathcal O_{\mathrm{pFF},K}^{+,\leq m}
:=
\R\lim_I
\R j_{I,*}\spp_{I}^{\leq m}.
\]
It is a solid quasi-coherent sheaf  on the curve $Y_{\FF}$, 
which is also  a $(\phi,N,G_K)$-module. We define similarly 
$\mathcal O_{\mathrm{pFF},L}^{+,\leq m}$ for a finite extension $L$ of $K$.  
If the context permits we will drop the subscript~$K$.

   More generally, for $S\in\operatorname{Perf}_C$, write
$ Y_{S,I}:=Y_{\mathrm{FF},S^\flat,I}$,
$j_{I,S}:Y_{S,I}\hookrightarrow Y_{\mathrm{FF},S^\flat}$,
and let $q_{I,S}:Y_{S,I}\to Y_I$ be the natural morphism. We set
\[
\spp_{I,S}^{\leq m}:=\LL q_{I,S}^{*}\spp_{I}^{\leq m},\qquad
\mathcal O^{+,\leq m}_{\mathrm{pFF},K,S}
 :=\R\lim_I \R j_{I,S,*}P_{I,S}^{\leq m}.
\]
 The sheaf $\spp_{I,S}^{\leq m}$ is finite free and
is associated to the relative module
$\mathbf B^{I,\leq m}_{\mathrm{uni},K}(S)$ of
\eqref{eve1}. We suppress $K$ and $S$ in the
notation when they are clear.
\begin{lemma}\label{tired}
For $m\geq 0$, we have\footnote{We denote by $\sd^{\prime}_{\phi,N,G_K}(C^{\nr}_{\Box})$ the derived category of solid $(\phi,N)$-modules over $C^{\nr}$, where we do not impose any condition on commutation of $\phi$ and $N$.}
$\R\Gamma(Y_{\FF,S^{\flat}},\so^{+,\leq m}_{\pFF})
\simeq \B^{+,\leq m}_{\pFF}(S)$ in $\sd^{\prime}_{\phi,N}(C^{\nr}_{\Box})$.
\end{lemma}
\begin{proof} We compute
\begin{align*}
\R\Gamma(Y_{\FF,S^{\flat}},\so^{+,\leq m}_{\pFF}) 
& \simeq \R\Gamma(Y_{\FF,S^{\flat}},\R\lim_I\R j_{I,*}\spp_{I}^{\leq m})
\simeq  \R\lim_I\R\Gamma(Y_{\FF,S^{\flat}},\R j_{I,*}\spp_{I}^{\leq m})\\
& \simeq  \R\lim_I\R\Gamma(Y^I_S,\spp_{I}^{\leq m})
\simeq  \R\lim_I\B^{I,\leq m}_{\rm uni}(S).
\end{align*}
Hence it suffices to show that $\B^{+,\leq m}_{\pFF}(S)
\simeq \R\lim_I\B^{I,\leq m}_{\rm uni}(S)$, that is, that 
$$\B^{+,\leq m}_{\pFF}(S)\simeq \lim_I\B^{I,\leq m}_{\rm uni}(S), 
\quad \R^1\lim_I\B^{I,\leq m}_{\rm uni}(S)=0. $$
The first isomorphism follows from the definition. 
The vanishing of $\R^1\lim_I$ is proved verbatim as~\cite[Lemme~3.14]{CGN3},
replacing $\B_I$ and $\B^{I,\leq m}_{\rm uni}$ in loc. cit. by $\B_I(S)$ and
$\B^{I,\leq m}_{\rm uni}(S)$ (but keeping $W_J:=H^1(G_K,\B_J)$ as in loc. cit.).
\end{proof}

   We set $\so^+_{n,m}:=t^{-n}\so^{+,\leq m}_{\pFF}$, for $n\in\Z, m\in\N$, and  
\[
\mathcal B_{\mathrm{pFF}}:=
\operatorname*{colim}_{[L:K]<\infty}
\operatorname*{colim}_{n,m}
\R\pi_{{\FF},*}\bigl(\mathcal O_{n,m,L}^{+}\bigr)
\in\mathrm{Nuc}(X_{\mathrm{FF}}).
\]
By Lemma \ref{tired}, we have
\begin{align*}
\R\Gamma(X_{\FF,S^{\flat}},\mathcal{B}_{\pFF}) 
& \simeq \operatorname*{colim}_{[L:K]<\infty}\colim_{n,m}\R\Gamma(X_{\FF,S^{\flat}},\R\pi_*\so^+_{n,m,L})
\simeq \operatorname*{colim}_{[L:K]<\infty}\colim_{n,m} t^{-n} \B^{+,\leq m}_{\pFF,L}(S)\\
& \simeq  \operatorname*{colim}_{[L:K]<\infty}\B_{\pFF,L}(S)\simeq \B_{\pFF}(S).
\end{align*}

\begin{lemma}\label{paris13a} Let $i\geq 0$. We have 
$$
\R\Hhom_{\FF,S^{\flat}}(i_{\infty,*}\mathbb{B}^{[0,i[}_{\dr}(S),\mathcal{B}_{\pFF})  \simeq 0.
$$
\end{lemma}
\begin{proof} Using the exact sequence
$$
0 \to \so(0,i)\to \so\to i_{\infty,*}\mathbb{B}^{[0,i[}_{\dr}(S)\to 0,
$$
we compute (with $\R\Hhom_{X}:=\R\Hhom_{\FF,S^{\flat}}$, 
$\R\Hhom_{Y}:=\R\Hhom_{Y_{\FF},S^{\flat}}$,
$\pi:=\pi_{\FF}$):
\begin{align*}
\R\Hhom_{X}(i_{\infty,*}\mathbb{B}^{[0,i[}_{\dr}(S),\R\pi_*\so^+_{n,m}) 
& \simeq  \R\Hhom_{X}([\so(0,i)\to \so][1],\R\pi_*\so^+_{n,m})\\
   & \simeq  \R\pi_*([\R\Hhom_{Y}(\LL\pi^*\so(0,i),\so^+_{n,m})
\to \R\Hhom_{Y}(\LL\pi^*\so,\so^+_{n,m})])[1]\\
   & \simeq  \R\pi_*([\R\Hhom_{Y}(\so\{-i\}\otimes\Q_p(i),\so^+_{n,m})\to \R\Hhom_{Y}(\so,\so^+_{n,m})])[1]\\
    & \simeq \R\pi_*[\so^+_{n,m}\{i\}\otimes\Q_p(-i)\stackrel{t^i}{\to} \so^+_{n,m}][1].
   \end{align*}
For $n\geq i$, the multiplication map $t^i$ factors as $\so^+_{n,m}\xrightarrow[\sim]{t^i}\so^+_{n-i,m}\to \so^+_{n,m}$. Hence
$$
\colim_{n\geq i} \R\Hhom_{X}(i_{\infty,*}\mathbb{B}^{[0,i[}_{\dr}(S),\R\pi_*\so^+_{n,m})
\simeq \colim_{n\geq i}\R\pi_*[\so^+_{n-i,m}\to  \so^+_{n,m}][1]\simeq 0,
$$
as wanted.
\end{proof}

\subsection{Extensions of Vector Spaces}
We pass now   to computing $\R\Hhom$'s in Vector Spaces between Banach-Colmez spaces and certain period sheaves.

\subsubsection{Extensions of Banach-Colmez spaces} 
Let $S\in {\rm Perf}_C$. We call Vector Spaces over $S$  (${\rm VS}(S)$ for short) the objects in the $\infty$-derived category 
$\sd(S_{\proeet},\Q_p)$ of $\underline{\Q}_p$-modules in the category of pro-\'etale sheaves on ${\rm Perf}_S$. If $S={\rm Spa}(C)$ we will just say Vector Spaces.
Recall the following result from \cite[Cor. 3.10]{AnLB}: 
\begin{proposition}[\rm Algebraic Fargues-Fontaine fully-faithfulness] \label{tired1} Let $S\in {\rm Perf}_C$. The global sections  functor 
$$
\R\tau^{\prime}_*:\quad  {\rm Perf}(X_{\FF,S^{\flat}})\to {\sd}(S_{\proeet},\Q_{p})
$$
is fully faithful. That is, for $\sff,\sg\in {\rm Perf}(X_{\FF,S^{\flat}})$, the natural morphism in  ${\sd}(S_{\proeet},\Q_p)$
$$\R\tau^{\prime}_*\R\Hhom_{\FF,S^{\flat}}(\sff,\sg)\to \R\Hhom_{\rm VS(S)}(\R\tau^{\prime}_*\sff,\R\tau^{\prime}_*\sg)
$$
is a quasi-isomorphism. 
\end{proposition}
The above proposition   means that we can compute extensions of Banach-Colmez spaces or on the Fargues-Fontaine curve or in {\rm VS}'s. Hence, we have in ${\sd}(S_{\proeet},\Q_{p})$:
\begin{align}\label{bc10}
& \R\Hhom_{\rm VS(S)}(\Q_p,\Q_p)  \stackrel{\sim}{\leftarrow}\Q_p, &&\R\Hhom_{\rm VS(S)}(\Q_p,{\mathbb G}_a)\stackrel{\sim}{\leftarrow}{\mathbb G}_a,\\
&  \R\Hhom_{\rm VS(S)}({\mathbb G}_a,{\mathbb G}_a)  \simeq{\mathbb G}_a\oplus{\mathbb G}_a(-1)[-1], 
&& \R\Hhom_{\rm VS(S)}({\mathbb G}_a,\Q_p)  \simeq{\mathbb G}_a(-1)[-1].\notag
\end{align}
Somewhat more generally: 
\begin{corollary}\label{zimno10}
Let $i\geq 0$. 
We have in ${\sd}(S_{\proeet},\Q_{p})$ 
\begin{align*}
\R\Hhom_{\rm VS(S)}({\mathbb B}^{[0,i[}_{\dr}, {\mathbb B}^{[a,b[}_{\dr})  &  \simeq   [{\mathbb B}^{[a,b[}_{\dr}\to {\mathbb B}^{[a-i,b-i[}_{\dr}],\quad a\leq b\in\Z,\\
\R\Hhom_{\rm VS(S)}({\mathbb B}^{[0,i[}_{\dr}, F^a{\mathbb B}_{\dr})  &  \simeq{\mathbb B}_{\dr}^{[a-i,a[}[-1], \quad \quad a\in\Z.
\end{align*}
\end{corollary}
\begin{proof} The first quasi-isomorphisms follow  immediately from Lemma \ref{paris3a} and Proposition \ref{tired1}. The last one from  the computation:
\begin{align*}
\R\Hhom_{\rm VS(S)}({\mathbb B}^{[0,i[}_{\dr}, F^a{\mathbb B}_{\dr}) & \simeq \R\lim_b \R\Hhom_{\rm VS(S)}({\mathbb B}^{[0,i[}_{\dr}, {\mathbb B}^{[a,b[}_{\dr})\\
&\simeq \R\lim_b[{\mathbb B}^{[a,b[}_{\dr}\to {\mathbb B}_{\dr}^{[a-i,b-i[} ]
\simeq  {\mathbb B}^{[a-i,a[}_{\dr}[-1].
\end{align*}
\end{proof}

\subsubsection{Comparison theorem for period sheaves} 
We quote the following comparison theorem from \cite[Sec. 4.2.2]{CGN2}, \cite[Cor. 2.19]{Tate}.

\begin{theorem} \label{hot1K} 
Let $X$ be a smooth dagger  variety over $K$.  
\begin{enumerate}
\item  Let $r\geq 2d$. We have  a natural, functorial in $S\in{\rm Perf}_C$, 
and compatible with Frobenius  quasi-isomorphism 
in  $\sd(\mathbb{B}(S)_{\Box})$: 
\begin{equation}\label{mak1}
\R\Gamma_{\proeet}(X_S,{\mathbb B})(r)  \simeq [\R\Gamma^{\mathbb B}_{\hk}(X_S,r)\lomapr{\iota^{\mathbb B}_{\hk}} \R\Gamma_{\dr}^{\mathbb B}(X_S,r)],
\end{equation}
where we  set
\begin{align*}
\R\Gamma^{{\mathbb B}}_{\hk}(X_S,r) &:=[\R\Gamma_{\hk}(X_C)\{r\}\otimes^{\LL_{\Box}}_{{C}^{\nr}}{\mathbb B}^+_{\rm FF}(S)]^{N=0},\\
\R\Gamma_{\dr}^{{\mathbb B}}(X_S,r) & :=\prod_{\Z}((\R\Gamma_{\dr}(X)\otimes^{\LL_{\Box}}_K{\mathbb B}^+_{\dr}(S))/F^r).
\end{align*}
and $\iota^{\mathbb B}_{\hk}$ is the Hyodo--Kato morphism.
\item If $X$ is quasi-compact, we have  a natural, functorial in $S$, and compatible with Frobenius  quasi-isomorphism 
 in  $\sd(\mathbb{B}(S)_{\Box})$: 
\begin{equation}\label{mak2}
\R\Gamma_{\proeet}(X_S,\mathbb{B}[1/t]) \simeq [\R\Gamma_{\hk}(X_C)\otimes^{\LL_{\Box}}_{{C}^{\nr}}\mathbb{B}_{\rm FF}(S)]^{N=0}.
\end{equation}
\item For $a\in\Z$,  we have  a natural, functorial in $S$  quasi-isomorphism 
 in  $\sd(\mathbb{B}^+_{\dr}(S)_{\Box})$: 
\begin{equation}\label{mak4}
\R\Gamma_{\proeet}(X_S,F^a{\mathbb B}_{\dr}) \simeq F^a(\R\Gamma_{\dr}(X)\otimes^{\LL_{\Box}}_K\mathbb{B}_{\dr}(S)).
\end{equation}
\item If $X$ is quasi-compact, we have  a natural, functorial in $S$  quasi-isomorphism 
 in  $\sd(\mathbb{B}_{\dr}(S)_{\Box})$: 
\begin{equation}\label{mak3}
\R\Gamma_{\proeet}(X_S,{\mathbb B}_{\dr}) \simeq \R\Gamma_{\dr}(X)\otimes^{\LL_{\Box}}_K\mathbb{B}_{\dr}(S).
\end{equation}
\end{enumerate}
\end{theorem}
\begin{remark}
This theorem in loc. cit. is stated for smooth partially proper varieties because that is what was needed there. However the proof works for  general smooth dagger varieties. 
\end{remark}

\subsubsection{Canonical resolutions} 
We quote the following result:
\begin{theorem}{\rm (Mac Lane resolution, \cite[Sec. 3]{Breen})}\label{app2b}
Let $\stt$ be a site, let $({\rm Sh}(\stt),\Lambda)$ be a ringed topos and let $P$ be a $\Lambda$-module on $\stt$. There exists a complex $M_{\stt,\Lambda}(P)$ of $\Lambda$-modules on $\stt$ with augmentation to $P$, functorial in $P$ and called the {\em Mac Lane complex of $P$}, such that:
\begin{itemize}
\item The augmentation $\epsilon: M_{\stt,\Lambda}(P)\to P$ is a quasi-isomorphism, i.e., $M_{\stt,\Lambda}(P)$ is a resolution of $P$.
\item The term $M_{\stt,\Lambda}(P)_i$ of the complex is a finite direct sum of the terms of the form $\Lambda^+[\Lambda^{s_i}\times P^{t_i}]$, where $s_i,t_i$ are integers depending only on $i$ and not on $P$. Here we set $\Lambda^+[\Lambda^s\times P^t]:=\Lambda[\Lambda^s\times P^t]/\Lambda[0]$. 
\end{itemize}
\end{theorem}
\begin{example} (1) Functoriality in $P $ gives, from the zero morphism $P\to 0$, a morphism of complexes $M_{\stt,\Lambda}(P)\to M_{\stt,\Lambda}(0).$
In degree $i$, it is induced by the canonical projection  $\Lambda^{s_i}\times P^{t_i}\to \Lambda^{s_i}$.
Since the augmentation $M_{\stt,\Lambda}(P)\to P$  is a quasi-isomorphism, we have $M_{\stt,\Lambda}(0)\simeq 0$. 
That is,  the  complex $M_{\stt,\Lambda}(0)$ is acyclic.

 (2)  By functoriality, the zero morphism $0\to P$  induces a morphism $M_{\Lambda}(0)\to M_{\Lambda}(P)$.  Since $M_{\Lambda}(0)$
 is acyclic, the quotient $\overline{M}_{\Lambda}(P):=[M_{\Lambda}(0)\to M_{\Lambda}(P)][1]$
is again a resolution of $P$, with terms of the form $\Lambda[\Lambda^s\times P^t]/\Lambda[\Lambda^s]$.
\end{example}

\subsubsection{Extensions of period sheaves} 
The $\R\Hhom$'s  in Corollary \ref{zimno10} involved only limits of Banach-Colmez spaces hence it was possible to derive them from the case of Banach-Colmez spaces (see~\eqref{bc10}).
 It is more difficult to computes such $\R\Hhom$'s when they involve colimits.
 This is different from analogous computations on the Fargues-Fontaine curve that can simply invoke
 compactness of the curve. 
\begin{proposition}\label{app1}  There exist natural quasi-isomorphisms in $\sd(\mathrm{Spa}(C)_{\proeet},\Q_p)$: 
\begin{align*}
& \R\Hhom_{\rm VS}(\Q_p,{\mathbb B}_{\dr})  \stackrel{\sim}{\leftarrow}{\mathbb B}_{\dr},\quad \R\Hhom_{\rm VS}(\Q_p,{\mathbb B}[1/t])\stackrel{\sim}{\leftarrow}{\mathbb B}[1/t],\\
&  \R\Hhom_{\rm VS}({\mathbb G}_a,{\mathbb B}_{\dr})  \simeq 0, \quad 
  \R\Hhom_{\rm VS}({\mathbb G}_a,{\mathbb B}[1/t])  \simeq 0,\\
  &  \R\Hhom_{\rm VS}({\mathbb G}_a,F^a{\mathbb B}_{\dr}) \simeq \mathrm{gr}^{a-1}_F\mathbb{B}_{\dr}[-1]\quad \mbox{for } a\in\Z.
\end{align*}
We have analogous computations with ${\mathbb B}_{\pdr}$ in place of ${\mathbb B}_{\dr}$ and with ${\mathbb B}_{\rm pFF}, {\mathbb B}_{\rm FF}$ in place of ${\mathbb B}[1/t]$. 
\end{proposition}
\begin{proof}
The first two quasi-isomorphisms are obvious.   We will now treat the remaining  ones. 

\vskip2mm
$\bullet$ {\em The case of $({\mathbb G}_a,{\mathbb B}_{\dr})$.}   We want to show that
$$\R\Hhom_{\rm VS}({\mathbb G}_a,{\mathbb B}_{\dr})  \simeq 0.
$$
 Let $S\in {\rm Perf}_{C}$ be an affinoid perfectoid space. It suffices to  show that there exist natural (in $S$) quasi-isomorphisms in   $\sd(\Q_p(S))$
$$
 \R\Hom_{S_{\proeet},\Q_p}({\mathbb G}_a,{\mathbb B}_{\dr})  \simeq 0.
$$
Or, equivalently, that there exist natural (in $S$) quasi-isomorphisms in   $\sd(\Z_p(S))$
\begin{equation}\label{zimno11}
 \R\Hom_{S_{\proeet},\Z_p}({\mathbb G}_a,{\mathbb B}_{\dr})  \simeq 0.
\end{equation}

We will skip the subscript $\proeet$ in the rest of the proof. 

  Set $\Lambda=\underline{\Z}_p$. To compute the left-hand side of \eqref{zimno11} we use the Mac Lane complex $\overline{M}_{S, \Lambda}({\mathbb G}_a)$ from Theorem \ref{app2b}:
$$ \R\Hom_{S,\Lambda}({\mathbb G}_a,{\mathbb B}_{\dr})\simeq  \R\Hom_{S,\Lambda}(\overline{M}_{S,\Lambda}({\mathbb G}_a),{\mathbb B}_{\dr}).
$$
 The terms of the complex $\overline{M}_{S,\Lambda}({\mathbb G}_a)$ are finite direct sums  of the sheaves  $\Lambda[\Lambda^s\times{\mathbb G}_a^t]/\Lambda[\Lambda^s]$.  
Since the sheaf ${\mathbb G}_a$ is represented by 
${\mathbb A}_S$ and the sheaf $\Lambda^s$
is represented by the perfectoid affinoid 
$\wh{\Lambda}^s={\rm Spa}(\scc(\Z_p,C),\scc(\Z_p,\so_C))^s$,   the sheaf $\Lambda^s\times{\mathbb G}_a^t$ is represented by ${\mathbb A}^t_{S_{s}}$, where we set $S_{s}:=S\times\wh{\Lambda}^s$.   By Yoneda Lemma, we have a natural quasi-isomorphism  in   $\sd(\mathbb{B}_{\dr}(S))$ 
$$
\R\Hom_{S,\Lambda}(\Lambda[\Lambda^s\times\mathbb{G}_a^t]/\Lambda[\Lambda^s],{\mathbb B}_{\dr})\simeq \wt{{\R\Gamma}}_{\proeet}({\mathbb A}^t_{S_s},{\mathbb B}_{\dr}),
$$
 where $\wt{\R\Gamma}_{\proeet}({\mathbb A}^t_{S_s},{\mathbb B}_{\dr})$ is the reduced cohomology:
 $$
 \wt{\R\Gamma}_{\proeet}({\mathbb A}^t_{S_s},{\mathbb B}_{\dr}):=[{\R\Gamma}_{\proeet}({\mathbb A}^t_{S_s},{\mathbb B}_{\dr})\stackrel{e^*}{\to} {\R\Gamma}_{\proeet}({S_s},{\mathbb B}_{\dr})],
 $$
 for the zero section $e: S_s\hookrightarrow \mathbb{A}^t_{S_s}$.
 
  To compute this cohomology, choose an exhaustive covering $\{U_n\}$ of ${\mathbb A}^t_K$ by dagger closed discs. We have 
 $$
 \R\Gamma_{\proeet}(U_{n,S_s},{\mathbb B}_{\dr})\simeq \R\Gamma_{\dr}(U_n)\otimes^{\LL_{\Box}}_K\mathbb{B}_{\dr}(S_s)\simeq \mathbb{B}_{\dr}(S_s),
 $$
 where the first quasi-isomorphism follows from  \eqref{mak3} and the second one from the fact that $\R\Gamma_{\dr}(U_n)\simeq H^0_{\dr}(U_n)\simeq K$.
Hence we have
 $$
 \R\Gamma_{\proeet}({\mathbb A}^t_{S_s},{\mathbb B}_{\dr})\simeq \mathbb{B}_{\dr}(S_s),\quad  \wt{\R\Gamma}_{\proeet}({\mathbb A}^t_{S_s},{\mathbb B}_{\dr})\simeq 0.
 $$
Yielding  $\R\Hom_{S,\Lambda}({\mathbb G}_a,{\mathbb B}_{\dr})\simeq 0$, as wanted. 
 
   The argument for $\mathbb{B}_{\rm pdR}$ is identical. 
 
\vskip2mm
 $\bullet$ {\em The case of $({\mathbb G}_a,{\mathbb B}[1/t])$.} 
We argue  as   in  the case of $({\mathbb G}_a, {\mathbb B}_{\dr})$ treated above. 
 We want to show that
$$\R\Hhom_{\rm VS}({\mathbb G}_a,{\mathbb B}[1/t])  \simeq 0.
$$
Let $S\in {\rm Perf}_{C}$ be an affinoid perfectoid space. It suffices to  show that there exist natural (in $S$) quasi-isomorphisms in $\sd(\Q_p(S))$ 
\begin{equation}\label{virus1}
 \R\Hom_{S,\Lambda}({\mathbb G}_a,{\mathbb B}[1/t])  \simeq 0.
\end{equation}

   To compute the left-hand side of \eqref{virus1} we use again the Mac Lane complex $\overline{M}_{S,\Lambda}({\mathbb G}_a)$:
$$ \R\Hom_{S,\Lambda}({\mathbb G}_a,{\mathbb B}[1/t])\simeq  \R\Hom_{S,\Lambda}(\overline{M}_{S,\Lambda}({\mathbb G}_a),{\mathbb B}[1/t]).
$$
  By Yoneda Lemma, we have a natural quasi-isomorphism in $\sd(\mathbb{B}(S)[1/t])$ 
$$
\R\Hom_{S,\Lambda}(\Lambda[\Lambda^s\times\mathbb{G}_a^t]/\Lambda[\Lambda^s],{\mathbb B}[1/t]))\simeq \wt{\R\Gamma}_{\proeet}({\mathbb A}^t_{S_s},{\mathbb B}[1/t]).
$$
To compute this cohomology, 
 choose an exhaustive covering $\{U_n\}$ of ${\mathbb A}^t_K$ by dagger closed discs. We have (we use here that ${\mathbb B}\simeq \R\lim_I\mathbb{B}_I$)
 $$
 \R\Gamma_{\proeet}(U_{n,S_s},{\mathbb B}[1/t])\simeq [\R\Gamma_{\hk}(U_{n,C})\otimes^{\LL_{\Box}}_{{C^{\nr}}}\mathbb{B}_{\rm FF}(S_s)]^{N=0}\simeq [\mathbb{B}_{\rm FF}(S_s)]^{N=0}\simeq \mathbb{B}(S_s)[1/t],
 $$
 where the first quasi-isomorphism follows from  \eqref{mak2} and the second one 
 from the fact that   $\R\Gamma_{\hk}(U_{n,C})\simeq H^0_{\hk}(U_{n,C})\simeq {C^{\nr}}$ with trivial monodromy.  
Hence we have $ \wt{\R\Gamma}_{\proeet}({\mathbb A}^t_{S_s},{\mathbb B}[1/t])\simeq 0. $
 Yielding 
 $\R\Hom_{S,\Lambda}({\mathbb G}_a,{\mathbb B}[1/t])\simeq 0$, as wanted. 
 The claim for $\mathbb{B}_{\rm FF}$ is an easy consequence.

By passing to a colimit, the claim for 
 $\mathbb{B}_{\rm pFF}$ reduces to the one for $\mathbb{B}^{\leq m}_{\rm pFF,L}$. We compute
 \begin{align*}
 \R\Gamma_{\proeet}(U_{n,S_s},{\mathbb B}^{\leq m}_{\rm pFF,L}) 
& \simeq  \R\Gamma_{\proeet}(U_{n,S_s},\lim_I{\mathbb B}^{I,\leq m}_{{\rm uni},L})[1/t]
\simeq \R\Gamma_{\proeet}(U_{n,S_s},\R\lim_I{\mathbb B}^{I,\leq m}_{{\rm uni},L})[1/t]\\
 &  \simeq  (\R\lim_I\R\Gamma_{\proeet}(U_{n,S_s},{\mathbb B}^{I,\leq m}_{{\rm uni},L}))[1/t].
 \end{align*}
  But ${\mathbb B}^{I,\leq m}_{{\rm uni},L}$ is a free module of finite rank  
$\mathbb{B}_{I}[\underline{X}_n]$ over $\mathbb{B}_{I}$,
 and we have  for $r\geq 2d$ (see \eqref{mak1})
  \begin{align*}
   \R\Gamma_{\proeet}(U_{n,S_s},{\mathbb B}_{I}(r)) 
& \simeq [\R\Gamma_{\hk}(U_{n,C})\{r\}\otimes^{\LL_{\Box}}_{{C^{\nr}}}\mathbb{B}_{I,\log}(S_s)]^{N=0}\to M_n(S_s)\{r\}]\\
    &\simeq [\mathbb{B}_{I}(S_s)\{r\}[0]\to  M_n(S_s)\{r\}],
  \end{align*}
where $M_n$ is a (uniformly in $n$) bounded complex 
of $\mathbb{B}_{\dr}^+$-modules killed by $t^r$.
This implies: 
  \begin{align*}
  (\R\lim_I\R\Gamma_{\proeet}(U_{n,S_s},{\mathbb B}^{I,\leq m}_{{\rm uni},L}))[1/t] 
& \simeq (\R\lim_I[\mathbb{B}^{I,\leq m}_{{\rm uni},L}(S_s)[0]
\to  M_n(S_s)[\underline{X}_n]])[1/t]\\
  & \simeq [\R\lim_I\mathbb{B}^{I,\leq m}_{{\rm uni},L}(S_s)[0]
\to  \R\lim_I M_n(S_s)[\underline{X}_n]][1/t]\\
& \simeq \mathbb{B}^{\leq m}_{\pFF,L}(S_s),
\end{align*}
as wanted. 

\vskip2mm
$\bullet$ {\em The case of $({\mathbb G}_a,F^a{\mathbb B}_{\pdr})$.} 
Finally, for the last claim of the proposition: the case of $F^a{\mathbb B}_{\dr}$ 
follows from Corollary \ref{zimno10} (with $i=1$) 
so we just need to focus on $F^a{\mathbb B}_{\pdr}$. We need to show
that, for \(a\in\mathbf Z\), there is a natural quasi-isomorphism
\[
{\R\Hhom}_{\mathrm{VS}}
(\mathbb G_a,F^a\mathbb B_{\mathrm{pdR}})
\simeq
\operatorname{gr}^{a-1}_F\mathbb B_{\mathrm{pdR}}[-1].
\]
Put \(\ell:=\log t\) and
\[
M_m:=F^a\mathbb B_{\mathrm{dR}}[\ell]_{\le m}
=\bigoplus_{j=0}^{m}F^a\mathbb B_{\mathrm{dR}}\,\ell^j.
\]
Thus $ F^a\mathbb B_{\mathrm{pdR}}=
\operatorname*{colim}_{m}M_m. $
By Corollary~\ref{zimno10}, with \(i=1\), it suffices to prove
that the canonical morphism
\[
\operatorname*{colim}_{m}
{\R\Hhom}_{\mathrm{VS}}
(\mathbb G_a,M_m)
\longrightarrow
{\R\Hhom}_{\mathrm{VS}}
\left(\mathbb G_a,\operatorname*{colim}_{m}M_m\right)
\]
is a quasi-isomorphism. Or that, for a fixed  affinoid perfectoid \(S\in\mathrm{Perf}_C\), 
\[\operatorname*{colim}_{m}
\operatorname{RHom}_{S,\Q_p}(\mathbb G_a,M_m)
\xrightarrow{\ \sim\ }
\operatorname{RHom}_{S,\Q_p}
  \left(\mathbb G_a,\operatorname*{colim}_{m}M_m\right).
\]

   Let $\mathbb G_a^+$ denote the sheaf
$T\longmapsto\mathcal O^+(T)$
on \(S_{\mathrm{pro\acute et}}\).
We have $
\mathbb G_a^+\otimes_{\mathbf Z_p}^{L}\mathbf Q_p
\simeq
\mathbb G_a. $
Hence, for every complex \(M\) of \(\mathbf Q_p\)-module sheaves, we have
\[{\R\Hhom}_{S,\mathbf Q_p}(\mathbb G_a,M)
\simeq
\R\Hhom_{S,\mathbf Z_p}(\mathbb G_a^+,M).
\]
Now, apply the reduced Mac Lane resolution of
Theorem~\ref{app2b} and Example~2.17 to
\(\mathbb G_a^+\), with \(\Lambda=\mathbf Z_p\).
This gives a resolution
\(L^\bullet\to\mathbb G_a^+\),  whose terms are finite direct sums of sheaves
of the form
\[
E_{s,t}:=
\Lambda[\Lambda^s\times(\mathbb G_a^+)^t]
\big/
\Lambda[\Lambda^s].
\]
With \(S_s:=S\times\widehat{\Lambda}^{\,s}\) as above,
the sheaf
\(\Lambda^{\,s}\times(\mathbb G_a^+)^t\)
is represented by the relative closed unit polydisc
\(\mathbb D^t_{S_s}\).
The inclusion of \(\Lambda^{\,s}\) is its zero section.
Hence, by derived Yoneda,
\[
{\R\Hom}_{S,\Lambda}(E_{s,t},M)
\simeq [
\R\Gamma_{\mathrm{pro\acute et}}(\mathbb D^t_{S_s},M)
\xrightarrow{\,e^*\,}
\R\Gamma_{\mathrm{pro\acute et}}(S_s,M)].
\]

     Both \(\mathbb D^t_{S_s}\) and \(S_s\) are quasi-compact
and quasi-separated hence  we have
\[
\operatorname*{colim}_{m}
\operatorname{RHom}_{S,\Lambda}(E_{s,t},M_m)
\xrightarrow{\ \sim\ }
\operatorname{RHom}_{S,\Lambda}
  \left(E_{s,t},\operatorname*{colim}_{m}M_m\right).
\]
This is inherited  by
\(\mathbb G_a\): we have
\[
\operatorname*{colim}_{m}
\operatorname{RHom}_{S,\Lambda}(\mathbb G_a,M_m)
\xrightarrow{\ \sim\ }
\operatorname{RHom}_{S,\Lambda}
  \left(\mathbb G_a,\operatorname*{colim}_{m}M_m\right),
\]
as wanted.
 \end{proof}

\subsubsection{Algebraic fully-faithfulness: complements}  
We will need the following  fully-faithfulness result. 
\begin{lemma}\label{deszcz1}
 Let $S\in {\rm Perf}_C$ and $\sff\in {\rm Perf}(X_{\FF,S^{\flat}})$. The natural morphisms in $\sd(S_{\proeet},\Q_p)$
\begin{align*}
\R\tau^{\prime}_*\R\Hhom_{\FF,{S^{\flat}}}(\sff,i_{\infty,*}F^a\mathbb{B}_{\pdr}(S)) & \to \R\Hhom_{{\rm VS}(S)}(\R\tau^{\prime}_*\sff,F^a{\mathbb B}_{\pdr}),\quad a\in \Z\cup \{-\infty\},\\
\R\tau^{\prime}_*\R\Hhom_{\FF,S^{\flat}}(\sff,\mathcal{B}_{\rm pFF}) & \to \R\Hhom_{{\rm VS}(S)}(\R\tau^{\prime}_*\sff,{\mathbb B}_{\rm pFF})
\end{align*}
are quasi-isomorphisms (we set $F^{-\infty}{\mathbb B}_{\pdr}:={\mathbb B}_{\pdr}$).
\end{lemma}
\begin{proof}  We note that $$\R\tau^{\prime}_*i_{\infty,*}F^a\mathbb{B}_{\pdr}(S)\simeq F^a\mathbb{B}_{\pdr},\quad \R\tau^{\prime}_*\mathcal{B}_{\pFF}\simeq \mathbb{B}_{\pFF}. 
$$
The second fact needs a brief argument: 
By definition $
\mathcal B_{\mathrm{pFF}}
=
\colim_{n,m,L}
\R\pi_{\FF,\ast}(\mathcal O^+_{n,m,L}). $
The functor $\R\tau'_\ast$ commutes with this filtered colimit. Hence it is
enough to show that 
\[
\R\tau'_\ast \R\pi_{\FF,\ast}(\mathcal O^+_{n,m,L})\simeq t^{-n}\mathbb{B}^{+,\leq m}_{\mathrm{pFF},L} .
\]
But we have
$ \R\tau'_\ast \R\pi_{\FF,\ast}(\mathcal O^+_{n,m,L})(S) \simeq
\R\Gamma(Y_{\mathrm{FF},S^{\flat}},\mathcal O^+_{n,m,L})$.  
Now $\mathcal O^+_{n,m,L}=t^{-n}\mathcal O^{+,\leq m}_{\mathrm{pFF},L}$, hence,
by Lemma~\ref{tired}, 
\[
\R\Gamma(Y_{\mathrm{FF}},\mathcal O^+_{n,m,L})(S)
\simeq
t^{-n}\mathbb{B}^{+,\leq m}_{\mathrm{pFF},L}(S),
\]
as wanted.

The claim of the lemma  is local for the pro-\'etale topology on $S$, so we may
assume that $S$ is strictly totally disconnected. Let $
  \mathcal P\in
  \left\{
    i_{\infty,*}F^a\mathbb B_{\mathrm{pdR}}(S),
    \mathcal B_{\mathrm{pFF}}
  \right\}. $
The full subcategory of $\mathrm{Perf}(X_{\mathrm{FF},S^\flat})$ consisting
of the objects $\mathcal F$ for which the comparison morphism
\[
  \R\tau'_*
  \R\mathcal Hom_{\mathrm{FF},S^\flat}(\mathcal F,\mathcal P)
  \longrightarrow
  \R\mathcal Hom_{\mathrm{VS}(S)}
  (\R\tau'_*\mathcal F,R\tau'_*\mathcal P)
\]
is a quasi-isomorphism is a thick subcategory.

By the usual d\'evissage of perfect complexes, we may first reduce to
vector bundles. After \'etale localization on $S$, the argument of
\cite[Prop.~II.3.1,  Cor.~II.2.20]{FS}, together with
\cite[Prop.~8.4.7]{KL15}, reduces further to the line bundles
$\mathcal O(s)$, $s\in\Z$. Twists of the Euler sequence
\[
  0\longrightarrow\mathcal O
  \longrightarrow\mathcal O(1)^{\oplus2}
  \longrightarrow\mathcal O(2)
  \longrightarrow0
\]
reduce the claim to $\mathcal O$ and $\mathcal O(1)$, and the
exact sequence
\[
  0\longrightarrow\mathcal O
  \longrightarrow\mathcal O(1)
  \longrightarrow i_{\infty,*}R_S
  \longrightarrow0
\]
reduces it finally to
$\mathcal F\in\{\mathcal O,i_{\infty,*}R_S\}$.
For $\sff=\so$, since $\R\tau^{\prime}_*\so\simeq\Q_p,$ this is a tautology. 
For  $\sff=i_{\infty, *}R_S$: in the case of $i_{\infty,*}\mathbb{B}_{\pdr}(S)$ and $\mathcal{B}_{\rm pFF}$
 both sides are zero by Lemma \ref{paris3a}  and Lemma \ref{paris13a} together with Proposition \ref{app1}.  
 
    For $i_{\infty,*}F^a\mathbb{B}_{\pdr}(S)$ and  $ a\in\Z$, we know, by Lemma  \ref{paris3a} and Proposition \ref{app1} that both sides are $\mathrm{gr}^{a-1}\mathbb{B}_{\pdr}(S)[-1]$. So it remains to check that the canonical map is the identity. But  instead of doing this we will argue in the following way: by  Proposition \ref{app1}  and the fact that solid tensor products commute with  colimits  we may replace $i_{\infty,*}F^a\mathbb{B}_{\pdr}(S)$ with $i_{\infty,*}F^a\mathbb{B}_{\dr}(S)$; then, by passing to limits, 
    we may replace $i_{\infty,*}F^a\mathbb{B}_{\dr}(S)$ with $i_{\infty,*}\mathbb{B}^{[a,b[}_{\dr}(S)$, for $a\leq b$; and invoke Lemma \ref{paris3a} and Corollary \ref{zimno10}. 
\end{proof}

\section{Conjecture $C_{\st}$ on  the Fargues-Fontaine curve} \label{better1}  
We now prove the two curve-level (internal) duality statements from which Conjecture  $C_{\st}$ will be deduced: The filtered de Rham cohomology sheaves are   obtained by dualizing the pro-\'etale sheaves against $i_{\infty,*}F^a\B_{\dr}$,
whereas the Hyodo--Kato cohomology sheaves  are obtained by dualizing against the sheaf $\mathcal{B}_{\pFF}$.  
The main inputs are the basic comparison theorem for pro-\'etale cohomology \cite{CN4}, Poincar\'e duality for Hyodo--Kato and filtered de Rham cohomology \cite{AGN},  and the extension calculations of  Section \ref{better2}.

\subsection{Cohomology sheaves on the Fargues-Fontaine curve}   
Let us look briefly at the cohomology sheaves defined in \cite{CGN2}. 

\subsubsection{The cohomology sheaves}\label{paris11} 
Let $X$ be a smooth partially proper  rigid analytic variety over~$K$. 
Recall (see \cite[Sec.\,4, Sec.\,3.2, Sec.\,3.1]{CGN2}) that, for $r\in\Z$,  we have $G_K$-equivariant solid quasi-coherent sheaves
$$  \se_{\proeet,c}(X_C, \Q_p(r)), \quad  \se_{\hk,c}(X_C, r),  $$
on the Fargues-Fontaine curve, which we call {\it the pro-\'etale sheaf}
 and {\it the twisted Hyodo--Kato sheaf}, respectively. 
Taking their cohomology we get  
\begin{align*}
& \R\Gamma(X_{\FF}, \se_{\proeet,c}(X_C, \Q_p(r)))\simeq \R\Gamma_{\proeet,c}(X_C,\Q_p(r)), \\
&   \R\Gamma(X_{\FF}, \se_{\hk,c}(X_C, r))\simeq[\R\Gamma_{\hk,c}(X_C)\{r\}\otimes^{\LL_{\Box}}_{{C}^{\nr}}\B^+_{\FF}]^{N=0,\phi=1}, 
\end{align*}
where $\{r\}$ means Frobenius divided by $p^r$.  For $r\geq 0$, we also have the twisted de Rham sheaf
$
\se_{\dr,c}(X_C, r)
$
such that 
$$
\R\Gamma(X_{\FF}, \se_{\dr,c}(X_C, r))\simeq \R\Gamma_{\dr,c}(X_C/\B^+_{\dr})/F^r
=:\R\Gamma_{\dr, c}(X_C,r).
$$
Similarly for the usual cohomology (with no compact support condition). 

The twisted de Rham sheaf $\se_{\dr,*}(X_C, r)$ 
is supported  at the point at infinity $x_{\infty}\in X_{\FF}$.   
The twisted Hyodo--Kato sheaf $\se_{\hk,c}(X_C, r)$ 
has  the fiber at $x_{\infty}$ ($m\geq 1$) 
\begin{equation}\label{paris7}
\R\Gamma_{\hk,c}(X_C)\otimes^{\LL_{\Box}}_{{C}^{\nr}}i_{\infty,*}(\B^+_{\dr}/t^m)\simeq i_{\infty,*}(\R\Gamma_{\dr,c}(X)\otimes^{\LL_{\Box}}_K\B^+_{\dr}/t^m),
\end{equation}
This follows from the Galois equivariant computation\footnote{More precisely, 
its $\B_I$-analog for a suitable closed interval $I\subset (0,\infty)$ with rational endpoints.}: 
\begin{align*}
[\R\Gamma_{\hk,*}(X_C)\otimes^{\LL_{\Box}}_{{C}^{\nr}}\B_{\FF}^+]^{N=0}\otimes^{\LL_{\Box}}_{\B}(\B^+_{\dr}/t^m)  & \stackrel{\sim}{\leftarrow} (\R\Gamma_{\hk,*}(X_C)\otimes^{\LL_{\Box}}_{{C}^{\nr}}\B)\otimes^{\LL_{\Box}}_{\B}(\B^+_{\dr}/t^m)\\
& \stackrel{\sim}{\leftarrow} \R\Gamma_{\hk,*}(X_C)\otimes^{\LL_{\Box}}_{{C}^{\nr}}(\B^+_{\dr}/t^m) \\
& \xrightarrow[\sim]{\iota_{\hk}}\R\Gamma_{\dr,*}(X)\otimes^{\LL_{\Box}}_K(\B^+_{\dr}/t^m),
\end{align*}
where $\iota_{\hk}$ is the Hyodo--Kato quasi-isomorphism.  We stress that the first quasi-isomorphism above is indeed Galois equivariant.

\subsubsection{Hyodo--Kato and de Rham dualities}
Let $X$ be as in Section \ref{paris11}. Let $d$ be the dimension of $X$. 
Let $r, r^{\prime}\in \Z$. Set $s=r+r^{\prime}-d$. We have a pairing in $\sd_{\phi,N,G_K}(C^{\nr}_{\Box})$ (see \cite[Sec. 5.1.3]{CGN3}) 
$$
\R\Gamma_{\hk}(X_C,r)\otimes_{C^{\nr}}^{\LL_{\Box}}\R\Gamma_{\hk,c}(X_C,r^{\prime})\to C^{\nr}\{s\}[-2d],
$$
which induces a duality quasi-isomorphism in $\sd_{\phi,N,G_K}(C^{\nr}_{\Box})$
\begin{equation}\label{paris31}
\R\Gamma_{\hk}(X_C,r)\stackrel{\sim}{\to}\R\Hhom_{\sd_{\phi,N,G_K}(C^{\nr}_{\Box})}(\R\Gamma_{\hk,c}(X_C,r^{\prime})[2d], C^{\nr}\{s\}),
\end{equation}
This yields  a pairing on $X_{\FF}$
$$
\se_{\hk}(X_C,r)\otimes^{\LL}_{\so}\se_{\hk,c}(X_C,r^{\prime})\to \so(s)[-2d],
$$
inducing a duality quasi-isomorphism in ${\rm QCoh}(X_{\FF})^{G_K}$ (see \cite[Lemma 5.12]{CGN2}): 
\begin{equation}\label{hk-dual}
\gamma_{X_C}:  \se_{\hk}(X_C,r)\stackrel{\sim}{\to}\R\Hhom_{\FF}( \se_{\hk,c}(X_C,r^{\prime})[2d],\so(s)).
\end{equation}
If $X_C$ has finite rank Hyodo--Kato cohomology then the sheaves $ \se_{\hk,*}(X_C,r^{\prime})$ are perfect and we can write
$$
\R\Hhom_{\FF}( \se_{\hk,c}(X_C,r^{\prime}),\so(s))\stackrel{\sim}{\leftarrow}\R\Hhom_{\FF}( \se_{\hk,c}(X_C,r^{\prime}),\so)\otimes^{\LL}_{\so}\so(s).
$$

\subsection{Conjecture $C_{\st}$: de Rham case} \label{paris10}  
We will prove in this section the de Rham version of Conjecture $C_{\st}$ on the Fargues-Fontaine curve. 
    
\subsubsection{Preliminary computations} 
\begin{lemma}\label{paris3}
Let $i\geq 0$ and $a\in\Z$.
\begin{enumerate}
\item  Let $b\in\Z$ be such that $a<0< b-i$. We have in $\sd(\Q_{p,\Box})$
\begin{align*}
\R\Gamma(G_K,\R\Hom_{\FF}(i_{\infty,*}\B^{[0,i[}_{\dr},i_{\infty,*}\B^{[a,b[}_{\pdr}))\simeq 0.
\end{align*}
\item We have in $\sd(\Q_{p,\Box})$
\begin{align*}
  &  \R\Gamma(G_K,\R\Hom_{\FF}(i_{\infty,*}\B^{[0,i[}_{\dr},i_{\infty,*}F^a\B_{\pdr}))\simeq 0,   \quad a\leq 0;\\
  &     \R\Gamma(G_K,\R\Hom_{\FF}(i_{\infty,*}\B^{[0,i[}_{\dr},i_{\infty,*}\B_{\pdr}))\simeq 0.
\end{align*}
\end{enumerate}
\end{lemma}
\begin{proof} 
Adding Galois action to the quasi-isomorphism from Lemma \ref{paris3a}, we get  
$$
\R\Gamma(G_K, \R\Hom_{\FF}(i_{\infty,*}\B^{[0,i[}_{\dr}, i_{\infty,*}\B^{[a,b[}_{\pdr}) ) \simeq [\R\Gamma(G_K,\B^{[a,b[}_{\pdr})\to  \R\Gamma(G_K,\B^{[a-i,b-i[}_{\pdr})]
$$
Now, for $a<0< b-i$,  by  Remark \ref{special} and the generalized Tate formula from \cite[Th. 1.4]{CGN3},  which gives that  $$ \R\Gamma(G_K,\B^{[b-i,b[}_{\pdr})\simeq  \R\Gamma(G_K,\B^{[a-i,a[}_{\pdr})\simeq 0.$$
we get the first quasi-isomorphism of the lemma. 

For the  second quasi-isomorphism  of the  lemma, the above implies that 
$$ \R\Gamma(G_K,\R\Hom_{\FF}(i_{\infty,*}\B^{[0,i[}_{\dr},i_{\infty,*}F^a\B_{\pdr}))\simeq  \R\Gamma(G_K,\B^{[a-i,a[}_{\pdr}))[-1]\simeq 0,
$$
as wanted. 
The third quasi-isomorphism follows by a colimit argument from this since $i_{\infty,*}\B^{[0,i[}_{\dr}$ is a compact object and $G_K$ is profinite. 
\end{proof}
\begin{corollary} \label{paris4} 
Let $i\geq 0, a\leq 0$ or $a=-\infty$. Let $W$ be a space of compact type over $K$ (with a trivial $G_K$-action). Then we have in $\sd(\Q_{p,\Box})$
$$
\R\Gamma(G_K, \R\Hom_{\FF}(W\otimes^{\Box}_Ki_{\infty,*}\B^{[0,i[}_{\dr},i_{\infty,*}F^a\B_{\pdr}))\simeq 0.
$$
\end{corollary} 
\begin{proof} Write $W=\colim_nW_n$, where each $W_n$ is Smith. Limit-arguments show that we may assume $W$ itself is Smith. Now, we compute 
\begin{align*}
\R\Gamma(G_K, & \R\Hom_{\FF}(W\otimes^{\Box}_K i_{\infty,*}\B^{[0,i[}_{\dr},i_{\infty,*}F^a\B_{\pdr}))\\
& \simeq     \R\Gamma(G_K,  W^{*}\otimes^{\LL_{\Box}}_K\R\Hom_{\FF}(i_{\infty,*}\B^{[0,i[}_{\dr},i_{\infty,*}F^a\B_{\pdr}))\\
	& \simeq     W^{*}\otimes^{\LL_{\Box}}_K\R\Gamma(G_K, \R\Hom_{\FF}(i_{\infty,*}\B^{[0,i[}_{\dr},i_{\infty,*}F^a\B_{\pdr}))\simeq 0
\end{align*}
The first quasi-isomorphism uses the fact that $W$ is Smith; second --  follows from \cite[Prop. 2.3]{DN1}; and  the last one --  from Lemma \ref{paris3}. 
\end{proof}

\subsubsection{$\B_{\pdr}$-dual of Hyodo--Kato cohomology}     
Let $X$ be a smooth partially proper  rigid analytic variety over $K$.
\begin{lemma} \label{paris9}Let $r\geq d, a\in\Z\cup\{-\infty\}$. Assume that  de Rham cohomology of $X$ has  finite rank. 
\begin{enumerate}
\item We have in ${\rm QCoh}(X_{\FF})$
$$
\R\Hhom_{\FF}(\se_{\hk,c}(X_C,r)[2d],i_{\infty,*}F^{a}\B_{\pdr}) \simeq  \R\Gamma_{\dr}(X)\otimes^{\LL_{\Box}}_{K} i_{\infty,*}F^{a}\B_{\pdr}.
$$
The same formula holds with \(\B_{\mathrm{dR}}\) in place of \(\B_{\mathrm{pdR}}\).
\item We have in $\sd(\Q_{p,\Box})$
$$
\R\Gamma(G_K,    \R\Hom_{\FF}(\se_{\hk,c}(X_C,r)[2d],i_{\infty,*}F^{a}\B_{\pdr}))\simeq \begin{cases}  \R\Gamma_{\dr}(X) & \mbox{ if  } a\leq  0,  \\
0 & \mbox{ if  } a>  0. 
 \end{cases}
$$
\end{enumerate}
\end{lemma}
\begin{proof}   Since $\se_{\hk,c}(X_c,r)$ is perfect, we have (for $r^{\prime}\geq d, s=r+r^{\prime}-d$)
\begin{align*}
\R\Hhom_{\FF}(\se_{\hk,c}(X_C,r)[2d],i_{\infty,*}F^{a}\B_{\dr}) & \simeq   \R\Hhom_{\FF}(\se_{\hk,c}(X_C,r)[2d],\so)\otimes^{\LL}_{\so} i_{\infty,*}F^{a}\B_{\dr}\\
& \simeq \R\Hhom_{\FF}(\se_{\hk,c}(X_C,r)[2d],\so(s))(-s)\otimes^{\LL}_{\so} i_{\infty,*}F^{a}\B_{\dr}\\
& \simeq \se_{\hk}(X_C,r^{\prime})\otimes^{\LL}_{\so} i_{\infty,*}F^{a}\B_{\dr}\\
 & \simeq \R\Gamma_{\dr}(X)\otimes^{\LL_{\Box}}_{K} i_{\infty,*}F^{a}\B_{\dr}.
 \end{align*}
We used here that $\wh{\so(-s)}_{\infty}\simeq \B_{\dr}^+$ as Galois equivariant modules. The fourth  quasi-isomorphism follows from \eqref{paris7}.  Analogously for $F^{a}\B_{\pdr}$. This proves the first claim of our lemma.

Hence, using \cite[Prop. 2.3]{DN1}  we have 
\begin{align*}
\R\Gamma(G_K,    \R\Hom_{\FF}(\se_{\hk,c}(X_C,r)[2d],i_{\infty,*}F^{a}\B_{\pdr})) & \simeq \R\Gamma(G_K,\R\Gamma_{\dr}(X)\otimes^{\LL_{\Box}}_K F^{a}\B_{\pdr})\\
& \simeq \R\Gamma_{\dr}(X)\otimes^{\LL_{\Box}}_K \R\Gamma(G_K, F^{a}\B_{\pdr}),
\end{align*}
as wanted. The second claim of our lemma now follows from the generalized Tate theorem (see \cite[Th. 1.4]{CGN3}). 
\end{proof}

\subsubsection{$\B_{\pdr}$-dual of de Rham cohomology}    
Let $X$ be a smooth Stein   rigid analytic variety over $K$.
Let  $d\leq r$.  Recall that we have
\begin{align}\label{paris5} 
\R\Gamma_{\dr,c}&(X_C,r)\simeq\\
&\big(H^d_c(X,\so)\otimes^{\Box}_K\B^{[0,r[}_{\dr} 
 \to H^d_c(X,\Omega^1)\otimes^{\Box}_K\B^{[0,r-1[}_{\dr} \to\cdots\to 
 H^d_c(X,\Omega^d)\otimes^{\Box}_K\B^{[0,r-d[}_{\dr}\big)[-d].\notag
\end{align}
\begin{lemma} \label{paris8}
Let $r\geq d$, $a\in \Z$. Set ${\cal E}:=\se_{\dr,c}(X_C,r)[2d]$ and $\iota:=i_{\infty,*}$.  
\begin{enumerate}
\item We have in ${\rm QCoh}(X_{\FF})$
\begin{equation}\label{essen1a}
\R\Hhom_{\FF}( {\cal E},\iota F^a\B_{\dr}) \simeq  \iota (F^{a-r+d}(\R\Gamma_{\dr}(X)\otimes_K^{\LL_{\Box}}\B_{\dr})/\R\Gamma_{\dr}(X)\otimes_K^{\LL_{\Box}}F^{a}\B_{\dr})[-1].
\end{equation}
\item We have in $\sd(\Q_{p,\Box})$
$$
\R\Gamma(G_K, \R\Hom_{\rm FF}( {\cal E}, \iota F^a\B_{\pdr}))\simeq \begin{cases} 0 & \mbox{ if } a\leq 0,\\
F^{a-r+d}\R\Gamma_{\dr}(X)[-1]  &  \mbox{ if } a > 0.
\end{cases}
$$
\end{enumerate}
\end{lemma}
\begin{proof} Computing as in the proof of Corollary \ref{paris4}, we have (note that $F^a\B_{\dr}$ is Fr\'echet)
\begin{align*}
\R\Hhom_{\FF}(H^d_c(X,\Omega^i)\otimes^{\Box}_{K}\iota  \B^{[0,r-i[}_{\dr},\iota F^a\B_{\dr})& \simeq H^d_c(X,\Omega^i)^{*}\otimes^{\LL_{\Box}}_{K}  \R\Hhom_{\FF}(\iota \B^{[0,r-i[}_{\dr},\iota F^a\B_{\dr})\\
& \simeq \Omega^{d-i}(X)\otimes^{{\Box}}_K  \iota  \B^{[a-r+i,a[}_{\dr}[-1].
\end{align*}
The last quasi-isomorphism is by Serre duality and \eqref{paris6}. 
Hence
\begin{align}\label{essen1}
\R\Hhom_{\FF}( {\cal E},\iota F^a\B_{\dr}) & \simeq (\so(X)\otimes^{\Box}_K  \iota  \B^{[a-r+d,a[}_{\dr}\to\cdots\\
&\to  \Omega^d(X)\otimes^{\Box}_K  \iota  \B^{[a-r,a[}_{\dr})[-1]\notag \\
& \simeq \iota (F^{a-r+d}(\R\Gamma_{\dr}(X)\otimes_K^{\LL_{\Box}}\B_{\dr})/\R\Gamma_{\dr}(X)\otimes_K^{\LL_{\Box}}F^{a}\B_{\dr})[-1].\notag
\end{align}
This gives the first claim of the lemma.

  Now,  we claim that, for \(0\leq i\leq d\) and a space \(W\) of compact type over \(K\), endowed with the trivial \(G_K\)-action, we have functorially
$$
\R\Gamma(G_K, \R\Hom_{\FF}(W\otimes^{\Box}_{K}\iota  \B^{[0,r-i[}_{\dr},\iota F^a\B_{\pdr})) \simeq \begin{cases}
W^*[-1] & \mbox{ if } 1\leq a\leq r-i,  \\
0 & \mbox{ otherwise}.
\end{cases}
$$
Indeed, assuming first that $W$ is Smith we compute
 as in the proof of Corollary \ref{paris4}: 
\begin{align*}
\R\Hhom_{\FF}(W\otimes^{\Box}_{K}\iota  \B^{[0,r-i[}_{\dr},\iota F^a\B_{\pdr})& \simeq W^{*}\otimes^{\LL_{\Box}}_{K}  \R\Hhom_{\FF}(\iota \B^{[0,r-i[}_{\dr},\iota F^a\B_{\pdr})\\
& \simeq W^*\otimes^{{\Box}}_K  \iota  \B^{[a-r+i,a[}_{\pdr}[-1].
\end{align*}
The last quasi-isomorphism is by Lemma \ref{paris3a}. Taking global sections we get
$$
\R\Hom_{\FF}(W\otimes^{\Box}_{K}\iota  \B^{[0,r-i[}_{\dr},\iota F^a\B_{\pdr}) \simeq W^*\otimes^{{\Box}}_K    \B^{[a-r+i,a[}_{\pdr}[-1].
$$
Applying $\R\Gamma(G_K,-)$  and the generalized Tate formula we get our claim. The case of general $W$ follows by writing $W=\colim_nW_n$ as a colimit of Smith spaces and passing to limits in the proved  formula for Smith spaces. 

   We now apply the above to the terms of $\se$ to compute $\R\Gamma(G_K,\R\Hom_{\FF}( {\cal E},  \iota F^a\B_{\pdr})) $.   For $a\leq 0$, this implies  the acyclicity we want in the lemma. For $a> 0$,  we get
\begin{align*}
\R\Gamma(G_K,\R\Hom_{\FF}( {\cal E},  \iota F^a\B_{\pdr}))  
& \simeq (\Omega^{a-r+d}(X)\to\cdots\to \Omega^d(X))[-a+r-d-1]\\
& \simeq F^{a-r+d}\R\Gamma_{\dr}(X)[-1],
\end{align*}
as wanted. 
    \end{proof}
    \begin{remark} \label{banana1}  (1) The first claim of Lemma \ref{paris8} works with $\B^+_{\dr}$ in the target replaced by $\wt{\B}^+_{\pdr}$ defined in Appendix \ref{second1} (with the same proof). This is because the latter is Fr\'echet. An analogous remark applies to many situations in the rest of the paper. 
    
    (2) The computation in the proof of Lemma \ref{paris8}, shows that the canonical map
    $$
    \R\Gamma(G_K, \R\Hom_{\FF}( {\cal E},  \iota  F^a\B_{\dr})[\log t])\to \R\Gamma(G_K,\R\Hom_{\FF}( {\cal E},  \iota F^a\B_{\pdr}))
    $$
    is a quasi-isomorphism. 
    
    (3)  For $a>0$, 
the  quasi-isomorphism
$$
F^{a-r+d}\R\Gamma_{\dr}(X)[-1]\stackrel{\sim}{\to  } \R\Gamma(G_K,\R\Hom_{\FF}( \se_{\dr,c}(X_C,r)[2d],i_{\infty,*}F^a\B_{\pdr})) 
$$
from the proof of Lemma \ref{paris8} 
comes from the map (we use point (2))
$$
(\Omega^{a-r+d}(X)\to \cdots\to \Omega^d(X))[-a+r-d]\to (\so(X)\otimes^{\Box}_K  \B^{[a-r+d,a[}_{\dr}\to\cdots\to \\
\Omega^d(X)\otimes^{\Box}_K   \B^{[a-r,a[}_{\dr})
$$
induced by the canonical, Galois equivariant maps 
$$
\Omega^{i}(X)\to \Omega^{i}(X)\otimes^{\Box}_K   \B^{[a-r+d-i,a[}_{\dr},
$$
which induce the quasi-isomorphisms
$$
\Omega^{a-r+d+i}(X)\stackrel{\sim}{\to} \R\Gamma(G_K,\Omega^{a-r+d+i}(X)\otimes^{\Box}_K  \B^{[-i,a[}_{\pdr}).
$$
\end{remark}

\subsubsection{The main theorem} 
We are now ready to prove the main result of  Section \ref{paris10}. 
\begin{theorem} \label{essen3} {\rm (Conjecture $C_{\st}$: de Rham part)} Let $X$ be a smooth partially proper rigid analytic variety over $K$ of dimension $d$. Let $a\in\Z\cup\{-\infty\}$. 
\begin{enumerate}
\item If $X$ is Stein with finite rank de Rham cohomology then we have in ${\rm QCoh}(X_{\FF})$
\begin{equation}\label{gateway10}
\R\Hhom_{\FF}(\se_{\proeet,c}(X_C,\Q_p(d))[2d],i_{\infty,*}F^{a}\B_{\dr})\simeq i_{\infty,*}F^{a}(\R\Gamma_{\dr}(X)\otimes_K^{\LL_{\Box}}\B_{\dr}).
\end{equation}
\item We have in $\sd(K_{\Box})$
$$
\R\Gamma(G_K,\R\Hom_{\FF}(\se_{\proeet,c}(X_C,\Q_p(d))[2d],i_{\infty,*}F^{a}\B_{\pdr}))\simeq F^{a}\R\Gamma_{\dr}(X).
$$
\end{enumerate}
\end{theorem}
\begin{proof}  
Assume first that $X$ is Stein with finite rank de Rham cohomology. Take $r\geq 2d$ and $a\in\Z$. 
Set ${\cal E}_{\eet}:=\se_{\proeet,c}(X_C,\Q_p(r))[2d]$,
${\cal E}_{\dr}:=\se_{\dr,c}(X_C,r)[2d-1]$,
${\cal E}_{\hk}:=\se_{\hk,c}(X_C,r)[2d-1]$
  and $\iota:=i_{\infty,*}$.

It suffices to show that, in ${\rm QCoh}(X_{\FF})$ and $\sd(K_{\Box})$, respectively: 
\begin{align}\label{luminy1}
\R\Hhom_{\FF}({\cal E}_{\eet},\iota F^{a+r-d}\B_{\dr}) & \simeq \iota F^{a}(\R\Gamma_{\dr}(X)\otimes_K^{\LL_{\Box}}\B_{\dr}),\\
\R\Gamma(G_K,\R\Hom_{\FF}({\cal E}_{\eet},\iota F^{a+r-d}\B_{\pdr})) & \simeq F^{a}\R\Gamma_{\dr}(X).\notag
\end{align}

  From the basic comparison theorem (see \cite[Prop. 4.4]{CGN2}) we get the distinguished triangle 
\begin{equation}\label{paris20}
\xymatrix{
\se_{\proeet,c}(X_C,\Q_p(r))\ar[r] & \se_{\hk,c}(X_C, r)\ar[r] &  \se_{\dr,c}(X_C, r).
}
\end{equation}
Together with Lemma \ref{paris8} and Lemma \ref{paris9}, it yields the quasi-isomorphisms
(with $\otimes:=\otimes^{\LL_{\Box}}_K$, $\R\Gamma_{\dr}:=\R\Gamma_{\dr}(X)$,
$\R\Hhom:=\R\Hhom_{\FF}$):
\begin{align*}
\R\Hhom({\cal E}_{\eet}, & \iota F^{a+r-d}\B_{\dr})) \stackrel{\sim}{\to} 
\big[\R\Hhom({\cal E}_{\dr},\iota F^{a+r-d}\B_{\dr})
\stackrel{\iota_{\hk}^*}{\longrightarrow} \R\Hhom({\cal E}_{\hk},\iota F^{a+r-d}\B_{\dr})\big]\\
 &\simeq \big[\iota (F^a( \R\Gamma_{\dr}{\hskip.5mm{\otimes}\hskip.5mm} \B_{\dr})/\R\Gamma_{\dr}{\hskip.5mm{\otimes}\hskip.5mm} F^{a+r-d}\iota \B_{\dr})
\stackrel{\rm can}{\longrightarrow} \R\Gamma_{\dr}{\hskip.5mm{\otimes}\hskip.5mm} F^{a+r-d}\iota \B_{\dr})[1]\big]\\
& \stackrel{\sim}{\leftarrow} \iota F^a( \R\Gamma_{\dr}{\hskip.5mm{\otimes}\hskip.5mm} \B_{\dr}),
\end{align*}
as wanted in the first part of  \eqref{luminy1} for $a\in\Z$. 

Now assume that $a=-\infty$. We use again the distinguished triangle \eqref{paris20}. The de Rham term is zero:  its source is built from  \(t\)-torsion period pieces, while \(t\) is invertible in \(\B_{\mathrm{dR}}\).
The Hyodo--Kato term cohomology sheaf on the Fargues-Fontaine curve  is perfect under the finite-rank hypothesis, we can invoke Lemma \ref{paris9}. 

    We now prove the second claim  for every \(a\in\mathbf Z\cup\{-\infty\}\), still assuming that \(X\) is Stein with finite-rank de Rham cohomology. We claim that we have the canonical quasi-isomorphism
   $$
   \R\Gamma(G_K,\R\Hom_{\FF}({\cal E}_{\eet},  \iota F^{a+r-d}\B_{\dr})[\log t])\stackrel{\sim}{\to}\R\Gamma(G_K,\R\Hom_{\FF}({\cal E}_{\eet},  \iota F^{a+r-d}\B_{\pdr})).
   $$
   Indeed, using the exact sequence \eqref{paris20} it is enough to show  this for ${\cal E}_{\hk}$ and ${\cal E}_{\dr}$ in place of ${\cal E}_{\eet}$. But our map in the Hyodo-Kato case  is a quasi-isomorphism already before taking Galois invariants because the sheaf ${\cal E}_{\hk}$ is perfect;  the case of ${\cal E}_{\dr}$ follows from Remark \ref{banana1} for $a$ finite and "the $t$ invertible on the target" argument from above for $a=-\infty$. 
Now we can  compute using the first part 
\begin{align*}
\R\Gamma(G_K,\R\Hom_{\FF}({\cal E}_{\eet},\iota F^{a+r-d}\B_{\pdr})) & \simeq \R\Gamma(G_K,\R\Hom_{\FF}({\cal E}_{\eet},\iota F^{a+r-d}\B_{\dr})[\log t]) \\  & \simeq \R\Gamma(G_K,F^a( \R\Gamma_{\dr}(X)\otimes^{\LL_{\Box}}_K \B_{\pdr}))
  \simeq F^{a}\R\Gamma_{\dr}(X),
\end{align*}
as wanted. 

To treat the general case of the second claim of our theorem,  write $X\stackrel{\sim}{\leftarrow} U_{\jcdot}$, where $U_{\jcdot}$ is an analytic hypercovering with each $U_n$ a disjoint union of Stein varieties with finite rank de Rham cohomology. 
 We compute (using the special case treated above), setting $\se_{\eet}(-):=\se_{\proeet,c}(-,\Q_p(d))[2d]$,
 \begin{align*}
 \R\Gamma(G_K, \R\Hom_{\FF}(\se_{\eet}(X_C), & \iota F^a\B_{\pdr}) ) \stackrel{\sim}{\leftarrow} \R\Gamma(G_K,\R\Hom_{\FF}(\colim_n\se_{\eet}(U_{n,C}),\iota F^a\B_{\pdr}))\\
 & \simeq  \R\Gamma(G_K,\R\lim_n\R\Hom_{\FF}(\se_{\eet}(U_{n,C}),\iota F^a\B_{\pdr}))\\
& \simeq  \R\lim_n\R\Gamma(G_K,\R\Hom_{\FF}(\se_{\eet}(U_{n,C}),\iota F^a\B_{\pdr}))\\
 & \simeq \R\lim_nF^{a}\R\Gamma_{\dr}(U_n)
  \stackrel{\sim}{\leftarrow} F^{a}\R\Gamma_{\dr}(X),
\end{align*}
as wanted. 
\end{proof}

\subsection{Conjecture $C_{\st}$: Hyodo--Kato case} \label{paris33}
We will prove in this section the Hyodo--Kato version of Conjecture $C_{\st}$. 

\subsubsection{$\B_{\pFF}$-dual of Hyodo--Kato cohomology}
\begin{lemma} \label{paris12}  Let $X$ be a smooth partially proper rigid analytic variety over $K$ of dimension $d$ with finite rank de Rham  cohomology. Let $r\in \Z$. 
\begin{enumerate}
\item We have in ${\rm Mod}^{\prime}_{\phi,N,G_K}({\rm QCoh}(X_{\FF}))$ 
$$
\R\Hhom_{\FF}(\se_{\hk,c}(X_C,r)[2d],\mathcal{B}_{\pFF})  \simeq \colim_{n,m,L} \R\pi_{\FF,*}(\R\Gamma_{\hk}(X_C)\{d-r\}\otimes^{\LL_{\Box}}_{C^{\nr}}\so^+_{n,m,L}),
$$
where  the $C^{\nr}(\phi,N)$-structure  on the left hand side of \eqref{paris25}  comes from  $\mathcal{B}_{\pFF}$. 
\item We have in $\sd^{\prime}_{\phi,N,G_K}(C^{\nr}_{\Box})$ 
\begin{equation}\label{paris25}
\R\Hom_{\FF}(\se_{\hk,c}(X_C,r)[2d],\mathcal{B}_{\pFF}\otimes\Q_p(r-d))^{\R G_K-{\rm sm}} \simeq \R\Gamma_{\hk}(X_C),
\end{equation}
where  $(-)^{\R G_K-{\rm sm}}:=\colim_H\R\Gamma(H, -)$ with the colimit taken over  open subgroups $H$ of $G_K$.      \end{enumerate}
\end{lemma}

\begin{proof} 
Set ${\cal E}_c:=\se_{\hk,c}(X_C,r)[2d]$, $\R\Hhom_{Y}:=\R\Hhom_{Y_{\FF}}$,
$\pi:=\pi_{\FF}$.

For the first claim, we compute in ${\rm Mod}^{\prime}_{\phi,N,G_L}({\rm QCoh}(X_{\FF}))$:
\begin{align}\label{paris30}
\R\Hhom_{\FF}({\cal E}_c,\R\pi_{*}\so^+_{n,m,L}) & \simeq  \R\pi_{*}\R\Hhom_{Y}(\LL\pi^*{\cal E}_c,\so^+_{n,m,L})\\
& \simeq  \R\pi_{*}\R\Hhom_{Y}([\R\Gamma_{\hk,c}(X_C)\{r\}\otimes^{\LL_{\Box}}_{C^{\nr}}\so^+_{\FF}]^{N=0}[2d],\so^+_{n,m,L})\notag \\
& \simeq  \R\pi_{*}(\R\Hhom_{Y}([\R\Gamma_{\hk,c}(X_C)\{r\}\otimes^{\LL_{\Box}}_{C^{\nr}}\so^+_{\FF}]^{N=0}[2d],\so)\otimes^{\LL}_{\so}\so^+_{n,m,L}) \notag \\
& \simeq \R\pi_{*}([\R\Gamma_{\hk}(X_C)\{d-r\}\otimes^{\LL_{\Box}}_{C^{\nr}}\so^+_{\FF}]^{N=0}\otimes^{\LL}_{\so}\so^+_{n,m,L}).\notag \\
& \simeq \R\pi_{*}(\R\Gamma_{\hk}(X_C)\{d-r\}\otimes^{\LL_{\Box}}_{C^{\nr}}\so^+_{n,m,L}).\notag 
\end{align}
Here, the second quasi-isomorphism is induced by the canonical quasi-isomorphism in ${\rm Mod}_{\phi,N,G_K}({\rm QCoh}(Y_{\FF}))$ 
$$
\LL\pi^*\se_{\hk,c}(X_C,r)\stackrel{\sim}{\to} [\R\Gamma_{\hk,c}(X_C)\{r\}\otimes^{\LL_{\Box}}_{C^{\nr}}\so^+_{\FF}]^{N=0}.
$$
The third quasi-isomorphism follows from the fact that $\R\Gamma_{\hk,c}(X_C)\{r\}\otimes^{\LL_{\Box}}_{C^{\nr}}\so$ is perfect by our assumption on Hyodo--Kato cohomology. The fourth quasi-isomorphism follows from 
Hyodo--Kato duality (see \eqref{paris31}). 

For the second claim of the lemma, taking global sections of \eqref{paris30}, we get  in $\sd^{\prime}_{\phi,N,G_L}(C^{\nr}_{\Box})$: 
\begin{align*}
\R\Hom_{\FF}({\cal E}_c,\R\pi_{*}\so^+_{n,m,L})     & \simeq \R\Gamma(X_{\FF},\R\pi_{\FF,*}(\R\Gamma_{\hk}(X_C)\{d-r\}\otimes^{\LL_{\Box}}_{C^{\nr}}\so^+_{n,m,L}))\\
& \simeq \R\Gamma(Y_{\FF},\R\Gamma_{\hk}(X_C)\{d-r\}\otimes^{\LL_{\Box}}_{C^{\nr}}\so^+_{n,m,L})\\
& \simeq     \R\Gamma_{\hk}(X_C)\{d-r\}\otimes^{\LL_{\Box}}_{C^{\nr}}\R\Gamma(Y_{\FF},\so^+_{n,m,L})
\end{align*}
which yields 
$$
\R\Hom_{\FF}({\cal E}_c,\mathcal{B}_{\pFF})   \simeq  \R\Gamma_{\hk}(X_C)\{d-r\}\otimes^{\LL_{\Box}}_{C^{\nr}}\B_{\pFF}.
$$
This implies that in $\sd^{\prime}_{\phi,N,G_K}(C^{\nr}_{\Box})$ we have
\begin{align*}
\R\Hom_{\FF}({\cal E}_c,\mathcal{B}_{\pFF}\otimes \Q_p(r-d))^{\R G_K-{\rm sm}}   & \simeq  (\R\Gamma_{\hk}(X_C)\{d-r\}\otimes^{\LL_{\Box}}_{C^{\nr}}\B_{\pFF}\otimes\Q_p(r-d))^{\R G_K-{\rm sm}}\\
 & \simeq \R\Gamma_{\hk}(X_C),
\end{align*}
as wanted (the last quasi-isomorphism uses the fact that the Galois action on Hyodo--Kato cohomology groups is smooth).
\end{proof}

\subsubsection{$\B_{\pFF}$-dual of de Rham  cohomology} 
\begin{corollary}\label{paris13}   Let $X$ be a smooth partially proper rigid analytic variety over $K$ of dimension $d$. Let $r\geq d$. 
We have in ${\rm QCoh}(X_{\FF})$ 
$$
\R\Hhom_{\FF}(\se_{\dr,c}(X_C,r)[2d],\mathcal{B}_{\pFF}) \simeq 0
$$
\end{corollary}
\begin{proof}  By localizing on $X$, we may assume that $X$ is Stein. 
Let $i\geq 0$.  By scalar tensor–Hom adjunction and  Lemma \ref{paris13a}  we have
$$
\R\Hhom_{\FF}(H^d_c(X,\Omega^i)\otimes^{\Box}_{K}i_{\infty,*} \B^{[0,r-i[}_{\dr},\mathcal{B}_{\pFF})\simeq  0
$$
giving us  our  corollary. 
\end{proof}

\subsubsection{The main theorem}  
We are now ready to prove the main result of  Section \ref{paris33}. 
\begin{theorem} {\rm (Conjecture $C_{\st}$: Hyodo--Kato part)} \label{luminy20}Let $X$ be a smooth partially proper rigid analytic variety over $K$ of dimension $d$ with finite rank de Rham cohomology. 
\begin{enumerate}
\item  We have in ${\rm Mod}^{\prime}_{\phi,N,G_K}({\rm QCoh}(X_{\FF}))$ 
\begin{equation}\label{luminy10}
\R\Hhom_{\FF}(\se_{\proeet,c}(X_C,\Q_p(d))[2d],\mathcal{B}_{\pFF})  \simeq \colim_{n,m,L}\R\pi_{\FF,*}(\R\Gamma_{\hk}(X_C)\otimes^{\LL_{\Box}}_{C^{\nr}}\so^+_{n,m,L}).
\end{equation}
\item We have in $\sd^{\prime}_{\phi,N,G_K}(C^{\nr}_{\Box})$
$$
\R\Hom_{\FF}(\se_{\proeet,c}(X_C,\Q_p(d))[2d],\mathcal{B}_{\pFF})^{\R G_K-{\rm sm}}\simeq \R\Gamma_{\hk}(X_C).
$$
\end{enumerate}
\end{theorem}
\begin{proof}   Take $r\geq 2d$. 
For the first claim of our theorem,  it suffices to show that
in ${\rm Mod}^{\prime}_{\phi,N,G_K}({\rm QCoh}(X_{\FF}))$  we have
$$
\R\Hhom_{\FF}(\se_{\proeet,c}(X_C,\Q_p(r))[2d],\mathcal{B}_{\pFF}\otimes\Q_p(r-d))  \simeq \colim_{n,m,L}\R\pi_{\FF,*}(\R\Gamma_{\hk}(X_C)\otimes^{\LL_{\Box}}_{C^{\nr}}\so^+_{n,m,L}).
$$
But, by  the distinguished triangle \eqref{paris20}  and   Corollary  \ref{paris13},  we have   the quasi-isomorphism in  ${\rm Mod}^{\prime}_{\phi,N,G_L}({\rm QCoh}(X_{\FF}))$ 
$$
\R\Hhom_{\FF}(\se_{\proeet,c}(X_C,\Q_p(r))[2d],\mathcal{B}_{\pFF}\otimes\Q_p(r-d))   \stackrel{\sim}{\leftarrow} \R\Hhom_{\FF}(\se_{\hk,c}(X_C,r)[2d],\mathcal{B}_{\pFF}\otimes\Q_p(r-d)). 
$$
and,  by Lemma \ref{paris12}, a sequence of quasi-isomorphisms in  ${\rm Mod}^{\prime}_{\phi,N,G_L}({\rm QCoh}(X_{\FF}))$ 
\begin{align*}
 &  \R\Hhom_{\FF}(\se_{\hk,c}(X_C,r)[2d],\R\pi_{\FF,*}\so^+_{n,m,L}\otimes\Q_p(r-d)) \\
& \simeq \R\pi_{\FF,*}(\R\Gamma_{\hk}(X_C)\{d-r\}\otimes^{\LL_{\Box}}_{C^{\nr}}\so^+_{n,m,L}(r-d))\\
& \simeq \R\pi_{\FF,*}(\R\Gamma_{\hk}(X_C)\otimes^{\LL_{\Box}}_{C^{\nr}}\so^+_{n+r-d,m,L}),
\end{align*}
as wanted.

For the second claim of our theorem, we take global sections of  \eqref{luminy10} to get in $\sd^{\prime}_{\phi,N,G_K}(C^{\nr}_{\Box})$ 
$$
\R\Hom_{\FF}(\se_{\proeet,c}(X_C,\Q_p(d))[2d],\mathcal{B}_{\pFF})\simeq \R\Gamma_{\hk}(X_C)\otimes^{\LL_{\Box}}_{C^{\nr}}\B_{\pFF}.
$$
Applying $ (-)^{\R G_K-{\rm sm}}$, 
the Galois acyclicity of $\B_{\pFF}$ together with the smoothness of the $G_K$--action on $\R\Gamma_{\hk}(X_C)$, gives   in $\sd^{\prime}_{\phi,N,G_K}(C^{\nr}_{\Box})$ 
$$
\R\Hom_{\FF}(\se_{\proeet,c}(X_C,\Q_p(d))[2d],\mathcal{B}_{\pFF})^{\R G_K-{\rm sm}}\simeq \R\Gamma_{\hk}(X_C).
$$
\end{proof}

\section{Arithmetic Topological Vector Spaces} 
In this section
we will  define and study the category of Galois equivariant Topological Vector Spaces, which we will call Arithmetic Topological Vector Spaces (${\rm ATVS}$'s for short).

\subsection{Topological Vector Spaces}
We start with recalling briefly basic facts concerning  Topological Vector Spaces  (${\rm TVS}$'s for short) introduced in  \cite{TVS}. 

The category of Topological Vector Spaces   is the $\infty$-category of $\underline{\Q}_p$-modules in the $\infty$-derived category 
$\underline{\sd}({\rm Spa}(C),{\rm Solid})$ of topologically enriched presheaves 
on ${\rm sPerf}_C$ -- the category of strictly totally disconnected affinoids
over $C$ -- with values in solid abelian groups.  We will denote it by  $\underline{\sd}({\rm Spa}(C),\Q_{p,\Box})$.  We will denote by 
$\underline{\sd}({\rm Spa}(C),\Q_p)$ the analogous category with ${\rm Solid}$ replaced by ${\rm CondAb}$. 

We have analogous categories $\underline{\sd}(S,\Q_{p,\Box})$ and $\underline{\sd}(S,\Q_p)$, for any $S\in {\rm sPerf}_C$.  The former we will call TVS's over $S$.
\begin{remark} Let $\scc:= {\rm CondAb}, {\rm Solid}$. 
We like  to  think of a topologically enriched presheaf $\sff: {\rm sPerf}^{\rm op}_C\to \scc$ as a presheaf 
$$
\sff:{\rm sPerf}_C^{\rm op}\mapsto {\rm Ob}\scc
$$
plus a  collection of (structure) maps in $\scc$ 
\begin{equation}\label{alt1}
\underline{\sff}^{\prime}_{Y,T}:\quad  \sff(Y\times T)\to  [T, \sff(Y)],\quad Y\in {\rm sPerf}^{\rm op}_C, T\in {\rm ProFin},
\end{equation}
such that we have  commutative diagrams
\begin{equation}\label{alt1b}
\xymatrix{
\sff(Y\times T)\ar[r]^{\underline{\sff}^{\prime}_{Y,T}} \ar[d]^{(f\times g)^*}& [T, \sff(Y)]\ar[d]^{(g^*,f^*)} \\
  \sff(Y_1\times T_1)\ar[r]^{\underline{\sff}^{\prime}_{Y,T} }&   [T_1, \sff(Y_1)],
}
\end{equation}
for all maps $f: Y_1\to Y, g: T_1\to T$, together with a composition rule. 

A natural transformation $f:\sff\to\sg$ between two such presheaves    $\sff$, $\sg$ is given by  functorial maps in $\scc$
$$
f_Y: \sff(Y)\to \sg(Y),\quad Y\in {\rm sPerf}_C,
$$
which are  compatible with the structure maps \eqref{alt1}. 
\end{remark}

We list  the following properties\footnote{We refer the reader to \cite{TVS} for the definitions of the functors $\R\pi_*, \R\tau_*, \R\tau^{\prime}_*$ and the notation for $\Hom$-functors. Here we recall that in the solid setting $\Hom$'s acquire a superscript $\Box$. }
\begin{proposition}{\rm (\cite[Th. 1.1]{TVS})} \label{recall1}  Let $S\in {\rm sPerf}_C$.
\begin{enumerate} 
\item {\rm(Enriched fully-faithfulness)}  The canonical projection  functor from Vector Spaces to Topological Vector Spaces
$$
\R\pi_*: \sd(S_{\proeet},\Q_p)\to \underline{\sd}(S,\Q_p)
$$
is  fully faithful. More precisely,  for   $\sff,\sg \in \sd(S_{\proeet},\Q_p)$,  the canonical   morphism in $\underline{\sd}(S,\Q_p)$
$$
\R\pi_*\R\Hhom_{S,\Q_p}(\sff,\sg)\to \R\Hhom_{S^{\rm top},\Q_p}(\R\pi_*\sff,\R\pi_*\sg)
$$
is a quasi-isomorphism. 
\item  {\rm(Fargues-Fontaine fully-faithfulness)}  The global sections  functor 
\begin{equation}\label{gl1}
\R\tau_*: {\rm QCoh}(X_{\FF, S^{\flat}})\to \underline{\sd}(S,\Q_{p})
\end{equation}
is fully faithful when restricted to perfect complexes. That is, for $\sff,\sg\in {\rm Perf}(X_{\FF, S^{\flat}})$, the natural morphism in  $\underline{\sd}(S,\Q_p)$
$$\R\tau_*\R\Hhom_{{\rm QCoh}(X_{\FF,S^{\flat}})}(\sff,\sg)\to \R\Hhom_{S^{\rm top},\Q_p}(\R\tau_*\sff,\R\tau_*\sg)
$$
is a quasi-isomorphism. 
\item  {\rm (Compatibility of the algebraic and topological projections)} The global sections functor 
\begin{equation}\label{gl2}
\R\tau^{\prime}_*: {\rm QCoh}(X_{\FF, S^{\flat}})\to {\sd}(S_{\proeet},\Q_{p})
\end{equation}
is compatible with the functor $\R\tau_*$ when restricted to nuclear sheaves. That is, the
 following diagram commutes
$$
\xymatrix{
 {\rm Nuc}(X_{\FF,S^{\flat}})\ar[r]^-{\R\tau^{\prime}_*}\ar[rd]_{\R\tau_*} &  \sd(S_{\proeet},{\Q}_p)\ar[d]^{\R\pi_*}\\
  & \underline{\sd}(S,{\Q}_p).
  }
  $$
\end{enumerate}
\end{proposition}

\subsection{Categories ${\rm Perf}_S$}\label{hott2}   
Let $S\in {\rm Perf}_C$. It is easy to see that we have  canonical equivalences of categories 
\begin{equation}\label{goraco1}
{\rm Perf}_S\simeq ({\rm Perf}_{C})_{/S}\simeq ({\rm Perf}_{K})_{/S}, 
\end{equation}
where ${\rm Perf}_{K}$ is the category of perfectoid affinoids over ${\rm Spa}(K)$. 

   We enrich the category ${\rm Perf}_{K}$ in ${\rm Cond}$ by setting
$$
\uHom_{K}(Y_1,Y_2)(T):=\Hom_{K}(Y_1\times_{K}\wh{T}_{K},Y_2), \quad T\in {\rm ProFin},
$$
where we put $\wh{T}_{K}={\rm Spa}(\scc(T,K),\scc(T,\so_K))$. Note that $Y_1\times_{K}\wh{T}_{K}$ is perfectoid. Moreover,  we have  $\wh{T}_{K}\times_{K}{\rm Spa}(C)\simeq \wh{T}$ (recall that $\wh{T}={\rm Spa}(\scc(T,C),\scc(T,\so_C))$), which yields
that, for $Y\in {\rm Perf}_C$, we have $Y\times_{K}\wh{T}_{K}\simeq Y\times_C\wh{T}$. It follows that the equivalences in \eqref{goraco1} are ${\rm Cond}$-enriched (with the slice categories inheriting the canonical enrichment).

  Let $S_1,S_2\in {\rm Perf}_C$ and let $f: S_1\to S_2$ be a map over $K$.  Consider the following two functors. 
\begin{align*}
& c_f: \quad ({\rm Perf}_{K})_{/S_1}\to ({\rm Perf}_{K})_{/S_2},\quad (Y\to S_1)\mapsto (Y\to S_1\stackrel{f}{\to} S_2);\\
& b_f: \quad ({\rm Perf}_{K})_{/S_2}\to ({\rm Perf}_{K})_{/S_1},\quad (Y\to S_2)\mapsto (S_1\times_{S_2}Y\to S_1)
\end{align*}
It is easy to see that 
the functors $c_f$ and $b_f$ satisfy the enriched  adjunction $c_f \dashv b_f$.

   Let $\scc= {\rm CondAb},\, {\rm Solid}$. For $S \in {\rm Perf}_C$, let $\se(S)$ be the category of ${\rm Cond}$-enriched presheaves on ${\rm Perf_S}$ with values in $\scc$. For a map $f: S_1\to S_2$ as above we define two functors
\begin{align*}
& f_*:\quad \se(S_1)\to \se(S_2),\quad f_*\sff=\sff \circ b_f^{\rm op},\\
& f^*:\quad  \se(S_2)\to \se(S_1),\quad f^*\sff=\sff \circ c_f^{\rm op}.
\end{align*}
Explicitly, we set
\begin{align*}
(f_*\sff)(Y\to S_2) & =\sff(S_1\times_{S_2}Y\to Y_1),\\
(f^*\sff)(Y\to S_1) & =\sff(Y\to S_1\stackrel{f}{\to} S_2). 
\end{align*}
 Since $b_f^{\rm op} \dashv c_f^{\rm op}$,  it is easy to see that 
the functors $f^*$ and $f_*$ satisfy the usual  adjunction $f^* \dashv f_*$. The fact that this is also an enriched adjunction follows easily from the formula $$\uHom_C(\sff_1,\sff_2)(T)=\Hom_C(\sff_1\otimes_{\Z}\Z[T],\sff_2),\quad T\in {\rm ProFin}, \,\scc={\rm CondAb}.
$$

In the above we can replace $ {\rm CondAb},\, {\rm Solid}$ with the corresponding derived categories and $f^*, f_*$ with $ f^*,\R f_*$,  we can  pass to $\Q_p$-coefficients, and 
we can replace perfectoid spaces  with strictly totally disconnected spaces. 

\subsection{Arithmetic Topological Vector Spaces}  
We are now ready  to introduce  the category of  Arithmetic Topological Vector Spaces. 

\subsubsection{Definitions}
The Galois group $G_K$ acts continuously on $C$. Thus,  if we set 
$$\wh{G}_K:={\rm Spa}(\scc(G_K,K),\scc(G_K,\so_K)),$$ $\wh{G}_K$ acts on ${\rm Spa}(C)$ over $K$. We denote by $[{\rm Spa}(C)/\wh{G}_K]$ the simplicial perfectoid affinoid  space over  $K$, which is given in degree $n$ by ${\rm Spa}(C)\times_{{\rm Spa}(K)}\wh{G}^n_K$ with face and degeneracy maps given by the usual bar construction for the right $\wh{G}_K$--action on ${\rm Spa}(C)$.  

We define the category of ATVS's (the limit is  taken\footnote{We note here that for a $K$-linear map $f: S_1\to S_2$ between two perfectoid affinoids over $C$ the pullback map $ f^*$ on the respective categories in \eqref{austria1} is symmetric monoidal} in ${\rm CAlg}({\rm Cat}_{\infty})$)
\begin{align}\label{austria1}
\underline{\sd}({\rm Spa}(C),\Q_{p,\Box})^{G_K}  :=\underline{\sd}({\rm Spa}(C),\Q_{p,\Box})([{\rm Spa}(C)/\wh{G}_K]):=\lim_{\Delta}\underline{\sd}({\rm Spa}(C)\times_{{\rm Spa}(K)}\wh{G}^{\jcdot}_K,\Q_{p,\Box}),
\end{align}
where the  maps are defined via the pullback maps $f^*$ from Section \ref{hott2}. Since all transition functors are slice restrictions, they preserve colimits, tensor products, and internal Homs. It follows that the limit category is closed symmetric monoidal, with colimits, tensor products, and internal Homs computed levelwise.  (In the $\Hom$'s notation we add the subscript $G_K$ and drop ${\rm Spa}(K)$ and ${\rm Spa}(C)$ for the sake of brevity.) 
\begin{remark}[Change of the base field]
For a finite extension \(L/K\) contained in \(\overline K\), write
$$
\mathrm{ATVS}_L
:=
\underline{\sd}(\Spa(C),\mathbf Q_{p,\Box})^{G_L}
$$
for the category obtained by replacing \(G_K\) with
\(G_L=\Gal(\overline K/L)\). The inclusion \(G_L\subset G_K\)
induces a restriction functor
$$
\mathrm{Res}_{L/K}:\mathrm{ATVS}_K\longrightarrow\mathrm{ATVS}_L.
$$
For \(\sff,\sg\in\mathrm{ATVS}_K\), we set
$$
\R\underline{\Hom}^{\mathrm{sm}}_{\mathrm{ATVS}}(\sff,\sg)
:=
\operatorname*{colim}_{[L:K]<\infty}
\R\underline{\Hom}_{\mathrm{ATVS}_L}
\bigl(\mathrm{Res}_{L/K}\sff,\mathrm{Res}_{L/K}\sg\bigr),
$$
where, for \(L\subset L'\), the transition morphism is induced by
restriction from \(G_L\) to \(G_{L'}\). It is enough to take the
colimit over finite Galois extensions \(L/K\); the resulting object
carries a natural smooth \(G_K\)-action.
\end{remark}We have the canonical projections
$$
{\rm ev}_0:\quad  \underline{\sd}({\rm Spa}(C),-)^{G_K}\to \underline{\sd}({\rm Spa}(C),-),\quad \sff\mapsto \sff_0.
$$
We define similarly the equivariant categories   $$\underline{\sd}({\rm Spa}(C),\Q_{p})^{G_K},  \quad  \underline{\sd}({\rm Spa}(C),{\rm CondAb})^{G_K},  \quad   \underline{\sd}({\rm Spa}(C),{\rm Solid})^{G_K}.$$
\begin{example}[Examples of ATVS's] \label{expresheaves} Let $\scc= \Q_p, \Q_{p,\Box}, {\rm CondAb}, {\rm Solid}$. We define the $\infty$-category $ \underline{\sd}({\rm Spa}(K),\scc) $ analogously to $ \underline{\sd}({\rm Spa}(C),\scc) $. To generate examples of 
ATVS's, we will use the following: 
\begin{lemma} \label{berlin1} The natural functor
$$
\varepsilon^*: \quad  \underline{\sd}({\rm Spa}(K),\scc) \to \underline{\sd}({\rm Spa}(C),\scc)
$$
factors through $ \underline{\sd}({\rm Spa}(C),\scc)^{G_K}$. 
\end{lemma}
\begin{proof} The simplicial object $[\mathrm{Spa}(C)/\wh{G}_K]$  is augmented over $\mathrm{Spa}(K)$. For $\sff\in\underline{\sd}(\mathrm{Spa}(K),\scc)$, pull back $\sff$ along the structural maps
$Y_n=\mathrm{Spa}(C)\times_K\wh{G}_K^n\to \mathrm{Spa}(K).$
Functoriality of pullback provides the face and degeneracy compatibilities, hence an object in the limit category. Its level-zero term is $\varepsilon^* \sff$.
\end{proof}

We present now some 
examples of   enriched presheaves in $ \underline{\sd}(\mathrm{Spa}(C),\Q_{p,\Box})^{{G}_K}$.
\vskip2mm
(1)  {\em The constant presheaf $\underline{\Q}_p$.}  Recall that the constant sheaf $\underline{\Q}_p$ on ${\rm Perf}_{K,\proeet}$ is defined via 
$$\underline{\Q}_p: \{Y\mapsto {\scc(|Y|,\Q_p)}\}, \quad Y\in {\rm Perf}_K.$$ One can see that, a priori just a presheaf,  it is actually a sheaf (see  \cite[Lemma 3.14]{MW}). We define the presheaf in $\underline{\sd}(\mathrm{Spa}(K),{\rm Solid})$
$$\underline{\Q}_p: \{Y\mapsto \underline{\scc(|Y|,\Q_p)}\},$$
where  $\scc(|Y|,\Q_p)$ is equipped with the compact open topology.
 We have   $\R\pi_*\underline{\Q}_p\simeq\underline{\Q}_p$ (see \cite[Ex. 2.1]{TVS})), where $\R\pi_*: \sd(\mathrm{Spa}(K)_{\proeet},{\rm Ab})\to \underline{\sd}(\mathrm{Spa}(K),{\rm CondAb})$ is the canonical pushforward  from the pro-\'etale site ${\rm Perf}_{K,\proeet}$.  It follows that $\underline{\Q}_p$ is a sheaf for pro-\'etale topology (since so is $\pi_*\underline{\Q}_p$). We have  $\underline{\Q}_p\in  \underline{\sd}(\mathrm{Spa}(K),\Q_{p,\Box})$ and, by Lemma \ref{berlin1},  $\underline{\Q}_p\in  \underline{\sd}({\rm Spa}(C),\Q_{p,\Box})^{G_K}$ (we will call $\varepsilon^*(\underline{\Q}_p)$ simply $\underline{\Q}_p$).

If $W$ is  a complete Hausdorff  locally convex topological vector space  over $\Q_p$, we can analogously 
define constant presheaves $\underline{W}$ in  $\underline{\sd}({\rm Spa}(C),\Q_{p,\Box})^{G_K} $. 
In what follows, if this does not cause confusion, we simply write $W$ in place of $\underline{W}$. 
\vskip2mm
(2) {\em The presheaf ${\mathbb G}_a$.}  On ${\rm Perf}_{K,\proeet}$ and in $\underline{\sd}(\mathrm{Spa}(K),{\rm Solid})$ we have the presheaves 
$${\mathbb G}^{\rm alg}_a: \{Y={\rm Spa}(R_Y,R^+_Y)\mapsto {R_Y}\},\quad {\mathbb G}_a: \{Y={\rm Spa}(R_Y,R^+_Y)\mapsto \underline{R}_Y\},
$$ respectively,  where $R_Y$ is equipped with its natural  Banach topology.  The first presheaf is a sheaf by \cite[Th. 8.7]{SchD}. 
We have  $\R\pi_*({\mathbb G}^{\rm alg}_a)\simeq {\mathbb G}_a$ (see \cite[Ex. 2.1, Lemma 2.5]{TVS}).  It follows that ${\mathbb G}_a$  is a sheaf for pro-\'etale topology (as so is $\pi_*{\mathbb G}^{\mathrm{alg}}_a$).  We have ${\mathbb G}_a\in \underline{\sd}({\rm Spa}(C),\Q_{p,\Box})^{G_K}$ by Lemma \ref{berlin1}. 

\vskip2mm
(3) {\em De Rham period presheaves.}  For $i\geq 0$, in $\underline{\sd}(\mathrm{Spa}(K),\Q_{p,\Box})$  we have the presheaves 
\begin{equation}\label{periods1}
\mathbb{B}^+_{\dr}/t^i :\{Y\mapsto { \B}^+_{\dr}(Y)/t^i\},
\end{equation}
where the (objectwise)  period rings are equipped with their canonical topologies. 
Similarly, we have period sheaves $(\mathbb{B}^+_{\dr}/t^i)^{\mathrm{alg}}\in \underline{\sd}(\mathrm{Spa}(K),\Q_{p,\Box})$  (by forgetting the topology on the period rings in \eqref{periods1}). By \cite[Ex. 2.1]{TVS}, there are natural isomorphisms
\begin{equation}\label{periods2}
\R\pi_*(\mathbb{B}^+_{\dr}/t^i)^{\mathrm{alg}}\stackrel{\sim}{\leftarrow}\mathbb{B}^+_{\dr}/t^i,
\end{equation}
which implies that
the topological period presheaves $\mathbb{B}^+_{\dr}/t^i$  are sheaves for pro-\'etale topology.  They are  objects of  $\underline{\sd}({\rm Spa}(C),\Q_{p,\Box})^{G_K}$ by Lemma \ref{berlin1}. 

By passing to the limit over $i$, we get objects
$F^a\mathbb{B}_{\dr}\in \underline{\sd}({\rm Spa}(K),\Q_{p,\Box})$. By taking colimits we get  objects $\mathbb{B}_{\dr}, \mathbb{B}_{\pdr}\in \underline{\sd}({\rm Spa}(K),\Q_{p,\Box})$. They all have algebraic avatars, satisfy \eqref{periods2}, and yield objects 
in $\underline{\sd}({\rm Spa}(C),\Q_{p,\Box})^{G_K}$. 

\vskip2mm
(4) {\em $\FF$-period presheaves.}  In a similar way we define the objects $\mathbb{B}, \mathbb{B}[1/t], \mathbb{B}^+_{\FF}, \mathbb{B}_{\FF},\mathbb{B}_{\pFF}\in \underline{\sd}({\rm Spa}(K),\Q_{p,\Box})$ and their corresponding algebraic avatars. They all satisfy \eqref{periods2} (since we have $H^i_{\proeet}(Y,\mathbb{B}_I)=0$, for $i\geq 1$, any  closed interval $I\subset (0,\infty)$ with rational endpoints,  and $Y\in \mathrm{sPerf}_K$ (see \cite[Prop. C.6.2]{Zav}).) Using Lemma \ref{berlin1} we get analogous objects in $\underline{\sd}({\rm Spa}(C),\Q_{p,\Box})^{G_K}$.

\vskip2mm
(5) {\em Enriched Yoneda presheaves.}  Let $T\in {\rm ProFin}$. Consider the following enriched Yoneda presheaves (see \cite[Sec. 3.1.5]{TVS} for the nonequivariant version)
$$
h^{\rm top}_{T,K}:=h^{\rm top}_{\wh{T}_K},\quad h^{\rm top}_T:=h^{\rm top}_{\wh{T}}.
$$
Hence $h^{\rm top}_T=h^{\rm top}_{T,K}|_{{\rm Spa}(C)}$.  By Lemma \ref{berlin1}, the enriched presheaf $\underline{\Q}_p[h^{\rm top}_{T,K}]$ is in $\underline{\sd}({\rm Spa}(C),\Q_{p,\Box})^{G_K}$. 
\end{example}

\subsubsection{Equivariant $\Hom$'s}
For $\sff\in \underline{\sd}({\rm Spa}(C),\Q_{p})^{G_K}$, we define  the equivariant global section functor as 
$$
\R\Gamma_{G_K}({\rm Spa}(C),\sff):=\R\Gamma([{\rm Spa}(C)/\wh{G}_K],\sff).
$$ For $\sff_1,\sff_2\in \underline{\sd}({\rm Spa}(C),\Q_{p})^{G_K}$ and $\sh:=\R\Hhom_{G_K}(\sff_1,\sff_2)$, we define the global $\Hom$ as
$$
\R\uHom_{G_K}(\sff_1,\sff_2):=\R\Gamma_{G_K}({\rm Spa}(C), \sh).
$$
Similarly in the solid setting. 

We have the usual tensor-hom adjunctions: 
\begin{lemma}\label{adjoint1}
Let $\sff_1,\sff_2,\sff_3\in  \underline{\sd}({\rm Spa}(C),\Q_p)^{G_K}$. We have  functorial quasi-isomorphisms  in $\underline{\sd}({\rm Spa}(C),\Q_p)^{G_K}$ and in $\sd(\Q_{p})$,
respectively:
\begin{align*}
\R\Hhom_{G_K} (\sff_1\otimes_{\Q_p}\sff_2,\sff_3)\simeq \R\Hhom_{G_K}(\sff_1,\R\Hhom_{G_K}(\sff_2,\sff_3)),\\
\R\uHom_{G_K} (\sff_1\otimes_{\Q_p}\sff_2,\sff_3)\simeq \R\uHom_{G_K}(\sff_1,\R\Hhom_{G_K}(\sff_2,\sff_3)).
\end{align*}
Similarly in the solid setting. 
\end{lemma}
\begin{proof} The first quasi-isomorphism follows from  construction. For the second quasi-isomorphism, recall that, for $T\in {\rm ProFin}$, we have 
$$\R\uHom_{{\rm Spa}(C)^{\rm top},\Q_p} (\sff_1,\sff_2)(T)\simeq \R\Hom_{{\rm Spa}(C)^{\rm top},\Q_p} (\Q_p[T]\otimes_{\Q_p} \sff_1,\sff_2).
$$  It follows that, for all $T\in {\rm ProFin}$,  we have 
\begin{align*}
\R\uHom_{G_K} (\sff_1\otimes_{\Q_p}\sff_2,\sff_3)(T) & \simeq \R\Hom_{G_K} (\Q_p[T]\otimes_{\Q_p}\sff_1\otimes_{\Q_p}\sff_2,\sff_3)\\
 & \simeq \R\Gamma_{G_K}({\rm Spa}(C),\R\Hhom_{G_K} (\Q_p[T]\otimes_{\Q_p}\sff_1\otimes_{\Q_p}\sff_2,\sff_3))\\
 & \simeq \R\Gamma_{G_K}({\rm Spa}(C),\R\Hhom_{G_K}(\Q_p[T]\otimes_{\Q_p}\sff_1,\R\Hhom_{G_K}(\sff_2,\sff_3))\\
 & \simeq \R\Hom_{G_K}(\Q_p[T]\otimes_{\Q_p}\sff_1,\R\Hhom_{G_K}(\sff_2,\sff_3))\\
 & \simeq \R\uHom_{G_K}(\sff_1,\R\Hhom_{G_K}(\sff_2,\sff_3))(T),
\end{align*}
as wanted. 
Here  the third  quasi-isomorphism follows from what we have already proved. 
\end{proof}
\begin{remark}
  Let $\sff_1,\sff_2\in  \underline{\sd}({\rm Spa}(C),\Q_p)^{G_K}$. The equivariant structures on $\sff_1, \sff_2$ induce the canonical equivariant structure on $\R\Hhom_{{\rm Spa}(C)^{\rm top},\Q_p}(\sff_{1,0},\sff_{2,0})$ -- denote it by $\R\Hhom_{G_K}(\sff_{1,0},\sff_{2,0})$ --  and we have a canonical quasi-isomorphism 
  $$\R\Hhom_{G_K}(\sff_{1,0},\sff_{2,0})\stackrel{\sim}{\to} \R\Hhom_{G_K}(\sff_1,\sff_2).
  $$
   Similarly in the solid setting. 
  \end{remark}

\subsubsection{Equivariant cohomology}
  Let $X:={\rm Spa}(C)$.  Let $\sff\in \underline{\sd}(X,\Q_{p,\Box})^{G_K}$ be a $G_K$--equivariant  enriched presheaf. We claim that, if $\sff_0\in \underline{\sd}(X,\Q_{p,\Box})$  is derived rigid\footnote{See \cite[Sec. 3.2.8]{TVS} for the definition.}, then $\R\Gamma(X,\sff_0)$ can be canonically upgraded to $\R\Gamma(X,\sff_0)_{\sff}\in {\rm Mod}_{\Q_{p,\Box}[G_K]}(\sd({\rm Solid}))$. Indeed, the cosimplicial structure $$[n]\mapsto \R\Gamma(Y_n,\sff_n), \quad Y_n:=X\times_K\wh{G}^n_K,
$$ through the quasi-isomorphisms obtained from the rigid enriched structure on $\sff_0$ (recall that $\sff_n\simeq p_n^*\sff_0$, for $p_n: Y_n\to X$ the canonical projection): 
$$
\R\Gamma(Y_n,\sff_n)\simeq \R[G_K^n,\R\Gamma(X,\sff_0)]
$$
yields a structure of a derived solid $G_K$-module on $\R\Gamma(X,\sff_0)$, as wanted. Here, for a profinite set \(T\) and \(M\in \sd(\mathbf Q_{p,\square})\), we set
$\R[T,M]:=
\R\underline{\operatorname{Hom}}_{\mathbf Q_{p,\square}}
\bigl(\mathbf Q_{p,\square}[T],M\bigr).$

\begin{proposition}\label{lozko1}
Let $\sff\in \underline{\sd}(X,\Q_{p,\Box})^{G_K}$ be a $G_K$--equivariant  enriched presheaf. If $\sff_0\in \underline{\sd}(X,\Q_{p,\Box})$  is derived rigid, then there is a natural quasi-isomorphism in $\sd(\Q_{p,\Box})$
$$
\R\Gamma_{G_K}(X,\sff)\simeq\R\Gamma(G_K,\sff_0(X)_{\sff}).
$$
\end{proposition}
\begin{proof}
This basically follows from the paragraph above the proposition: We compute 
\begin{align*}
\R\Gamma_{G_K}(X,\sff) & \simeq \lim_{[n]\in\Delta}\R\Gamma(Y_n,\sff_n)\simeq \lim_{[n]\in\Delta}\R[G_K^n,\R\Gamma(X,\sff_0)]\\
 & \simeq\R\Gamma(G_K,\sff_0(X)_{\sff}),
\end{align*}
as wanted. 
\end{proof}

\subsubsection{Arithmetic Vector Spaces}
We define the category of Arithmetic Vector Spaces (AVS's  for short) by  the identical procedure to the one used for TVS's above. Briefly, the category of AVS's is defined as (the limit is  taken in ${\rm CAlg}({\rm Cat}_{\infty})$): 
\begin{align}\label{austria1a}
{\sd}({\rm Spa}(C)_{\proeet},\Q_{p})^{G_K}  :={\sd}({\rm Spa}(C)_{\proeet},\Q_{p})([{\rm Spa}(C)/\wh{G}_K]):=\lim_{\Delta}{\sd}({\rm Spa}(C)_{\proeet}\times_{{\rm Spa}(K)}\wh{G}^{\jcdot}_K,\Q_{p}),
\end{align}
where the  maps are defined via the pullback maps $f^*$ analogous to those in  Section \ref{hott2}. This category is closed symmetric monoidal, with tensor products and internal Homs computed levelwise.
 We have the canonical projections
$$
{\rm ev}_0:\quad  {\sd}({\rm Spa}(C)_{\proeet},\Q_p)^{G_K}\to {\sd}({\rm Spa}(C)_{\proeet},\Q_p),\quad \sff\mapsto \sff_0.
$$
We define similarly the equivariant category   $ {\sd}({\rm Spa}(C)_{\proeet},{\rm Ab})^{G_K}.$

\begin{lemma} \label{berlin1a}  Let $\scc= \Q_p, {\rm Ab}$.  Let $\varepsilon: {\rm Spa}(C)\to {\rm Spa}(K)$ be the canonical map. The pullback  functor
$$
\varepsilon^*: \quad  {\sd}({\rm Spa}(K)_{\proeet},\scc) \to {\sd}({\rm Spa}(C)_{\proeet},\scc)
$$
factors through $ {\sd}({\rm Spa}(C)_{\proeet},\scc)^{G_K}$ and induces an equivalence of closed symmetric monoidal $\infty$-categories 
$$
\varepsilon^*: \quad  {\sd}({\rm Spa}(K)_{\proeet},\scc) \stackrel{\sim}{\to} {\sd}({\rm Spa}(C)_{\proeet},\scc)^{G_K}
$$
\end{lemma}
\begin{proof} This follows from the fact that $[{\rm Spa}(C)/\wh{G}_K]$ is a pro-\'etale hypercovering of ${\rm Spa}(K)$ and we work with pro-\'etale hypersheaves. 
\end{proof}
It follows that the algebraic versions of the presheaves discussed in Example \ref{expresheaves} are AVS's. 

We have the canonical projection  functor from Arithmetic Vector Spaces to Arithmetic Topological Vector Spaces\footnote{Sometimes, abusively, we will  refer to objects of $\underline{\sd}(S,\Q_p)^{G_K}$ as ATVS's; we hope this does not cause confusion.} ($S={\rm Spa}(C)$):
$$
\R\pi_*: \quad \sd(S_{\proeet},\Q_p)^{G_K}\to \underline{\sd}(S,\Q_p)^{G_K}
$$
defined by applying the projections $
\R\pi_*:  \sd(Y_{n,\proeet},\Q_p)\to \underline{\sd}(Y_n,\Q_p)
$ levelwise in $[{\rm Spa}(C)/\wh{G}_K]$:  The functor $\R\pi_{Y_n,*}$ 
 is defined by condensed pro-\'etale cohomology: 
 $$\big(\R\pi_{Y_n,*}(\sff)(S\stackrel{f}{\to}Y_n)\big)(T)=\R\Gamma((S\times T)_{\proeet}, {\rm pr}^*_{Y_n}f^*\sff);$$
hence its definition is tautologically compatible with base change in $Y_n$ yielding  the equivariant functor $\R\pi_*$.

\subsubsection{Equivariant solid quasi-coherent sheaves on the Fargues-Fontaine curve} 
Let us review two definitions of $G_K$-equivariant solid quasi-coherent sheaves on the Fargues-Fontaine curve $X_{\FF}$. 
\vskip2mm
($\bullet$) {\em  Classical  definition.}   For a Huber pair $(A,A^+)$ with an action of the condensed group $G_K$, we write $\sd((A,A^+)_{\Box})^{G_K}$ for the derived $\infty$-category of solid $(A,A^+)_{\Box}$-modules with a semilinear action of $G_K$. Covering $X_{\FF}$ with such affinoids and gluing the respective categories of $G_K$-equivariant modules we obtain the category ${\rm QCoh}(X_{\FF})^{G_K}_{\rm cl}$ of $G_K$-equivariant solid quasi-coherent sheaves on $X_{\FF}$ (seen classically). 
\vskip2mm
($\bullet$) {\em  Simplicial definition.}  For $n\geq 0$, set $X_{\FF,n}:=X_{\FF, Y^{\flat}_n}$, the relative Fargues-Fontaine curve over the tilt of $Y_n={\rm Spa}(C)\times_K \wh{G}^n_K$. Denote by $[X_{\FF}/G_K]$ the simplicial adic   space, which is given in degree $n$ by $X_{\FF,n}$ with face and degeneracy maps given by the usual bar construction for the right $\wh{G}_K$--action on ${\rm Spa}(C)$.  

   We define the category of $G_K$-equivariant solid quasi-coherent sheaves on $X_{\FF}$ as the limit in ${\rm CAlg}({\rm Cat}_{\infty})$:
$$
{\rm QCoh}(X_{\FF})^{G_K}:={\rm QCoh}([X_{\FF}/G_K])=\lim_{[n]\in\Delta} {\rm QCoh}(X_{\FF,n}). 
$$
This category is closed symmetric monoidal, with tensor products computed levelwise.
Internal $\Hom$'s  with perfect first argument are also computed levelwise.
Similarly we can define the categories ${\rm Nuc}(X_{\FF})^{G_K}$, ${\rm Perf}(X_{\FF})^{G_K}$ of $G_K$-equivariant nuclear  solid sheaves and perfect complexes. 
\begin{remark}We note that, since $[X_{\FF}/G_K]$ is a simplicial space induced from a  pro-\'etale hypercovering of ${\rm Spa}(K)$, we have an equivalence
of closed symmetric monoidal $\infty$-categories
$$
{\rm Perf}(X_{\FF, K})\stackrel{\sim}{\to} {\rm Perf}(X_{\FF})^{G_K},
$$
where we wrote $X_{\FF, K}$ for the relative Fargues-Fontaine curve over $\mathrm{Spa}(K)^{\diamond}$. 
\end{remark}
We have the global sections functors 
\begin{align*}
 & \R\tau^{\prime}_*:\quad {\rm QCoh}(X_{\FF})^{G_K}\to {\sd}({\rm Spa}(C)_{\proeet},\Q_{p})^{G_K},\\
& \R\tau_*:\quad {\rm QCoh}(X_{\FF})^{G_K}\to \underline{\sd}({\rm Spa}(C),\Q_{p,\Box})^{G_K}
\end{align*}
defined levelwise by the functors $\R\tau^{\prime}_*$ and $\R\tau_*$ from \eqref{gl2} and \eqref{gl1}, respectively:
Since $\R\tau_{S,*}(\sff)$  is defined by $\R\tau_{S,*}(\sff)(Y\stackrel{f}{\to} S)=\R\Gamma(X_{\FF,Y^{\flat}},\LL f^*\sff)$, it is tautologically compatible with base change in $S$; hence
the levelwise functors $\R\tau_{Y_n,*}(-)$ define a morphism of simplicial diagrams and thus induce $\R\tau_*$ (similarly for $\R\tau^{\prime}_*$).

   These functors commute with colimits: this is true for their nonequivariant versions (because the solid tensor product commutes with colimits and so does the functor $\R\Gamma(X_{\FF,S^{\flat}},-)$) and the equivariant versions inherit 
 this property (since the simplicial transition functors  preserve colimits).
\vskip2mm
 ($\bullet$) {\em  Comparison.}   The two constructions above agree on nuclear objects:
 \begin{lemma}
 There is an equivalence of closed symmetric monoidal $\infty$-categories
 $${\rm Nuc}(X_{\FF})^{G_K}_{\rm cl}\stackrel{\sim}{\to} {\rm Nuc}(X_{\FF})^{G_K},$$
 where \(\operatorname{Nuc}(X_{\mathrm{FF}})_{\mathrm{cl}}^{G_K}\)  is the full subcategory of \(\operatorname{QCoh}(X_{\mathrm{FF}})_{\mathrm{cl}}^{G_K}\) consisting of objects whose underlying solid quasi-coherent sheaves are nuclear.
 \end{lemma}
 \begin{proof}  For a compact interval $I\subset(0,\infty)$ with rational
endpoints,  $A:=\B_I$, and  a nuclear \(A_{\mathrm{an}}\)-complex \(M\), we have  functorial quasi-isomorphisms
\[
\scc(G_K^n,A)\otimes^{\LL}_{A_{\mathrm{an}}}M
\simeq
\R\underline{\operatorname{Hom}}_{\mathbf Q_{p,\square}}
\bigl(\mathbf Q_{p,\square}[G_K^n],M\bigr).
\]
Together with \(\B_I(Y_n)_{\mathrm{an}}\simeq \scc(G_K^n,A)_{\mathrm{an}}\), these identify the derived bar construction for continuous semilinear \(G_K\)-actions with the simplicial descent construction. Compatibility with tensor products, restrictions between annuli, and Frobenius yields the claimed equivalence.
 \end{proof}

\subsubsection{Fully-faithfulness results} 
\begin{proposition}{\rm (\cite[Th. 1.1]{TVS})} \label{recall1a}  Let $S={\rm Spa}(C)$.
\begin{enumerate} 
\item {\rm(Enriched fully-faithfulness)}  The canonical projection  functor from Arithmetic Vector Spaces to Arithmetic Topological Vector Spaces
$$
\R\pi_*: \quad \sd(S_{\proeet},\Q_p)^{G_K}\to \underline{\sd}(S,\Q_p)^{G_K}
$$
is  fully faithful. More precisely,  for   $\sff,\sg \in \sd(S_{\proeet},\Q_p)^{G_K}$,  the canonical   morphism in $\underline{\sd}(S,\Q_p)^{G_K}$
$$
\R\pi_*\R\Hhom_{\rm AVS}(\sff,\sg)\to \R\Hhom_{\rm ATVS}(\R\pi_*\sff,\R\pi_*\sg)
$$
is a quasi-isomorphism. 
\item  {\rm(Fargues-Fontaine fully-faithfulness)}  The global sections  functor 
$$
\R\tau_*:\quad  {\rm QCoh}(X_{\FF})^{G_K}\to \underline{\sd}(S,\Q_{p})^{G_K}
$$
is fully faithful when restricted to perfect complexes. That is, for $\sff,\sg\in {\rm Perf}(X_{\FF})^{G_K}$, the natural morphism in  $\underline{\sd}(S,\Q_p)^{G_K}$
$$\R\tau_*\R\Hhom_{{\FF}}(\sff,\sg)\to \R\Hhom_{\rm ATVS}(\R\tau_*\sff,\R\tau_*\sg)
$$
is a quasi-isomorphism. 
\item  {\rm (Compatibility of the algebraic and topological projections)} The global sections functor 
$$
\R\tau^{\prime}_*:\quad  {\rm QCoh}(X_{\FF})^{G_K}\to {\sd}(S_{\proeet},\Q_{p})^{G_K}
$$
is compatible with the functor $\R\tau_*$ when restricted to nuclear sheaves. That is, the
following diagram commutes
$$
\xymatrix{
{\rm Nuc}(X_{\FF})^{G_K}\ar[r]^-{\R\tau^{\prime}_*}\ar[rd]_{\R\tau_*} &  \sd(S_{\proeet},{\Q}_p)^{G_K}\ar[d]^{\R\pi_*}\\
& \underline{\sd}(S,{\Q}_p)^{G_K}.
}
$$
\end{enumerate}
\end{proposition}
\begin{proof}  This follows directly from Proposition \ref{recall1},  the definitions of equivariant categories, and 
Proposition \ref{old2} cited below. 
\end{proof}
We quote the  following result\footnote{Clearly, this also follows from Proposition \ref{tired1}.} from \cite[Cor. 3.10]{AnLB}: 
\begin{proposition}[Equivariant Algebraic Fargues-Fontaine fully-faithfulness] \label{old2} The global sections  functor 
$$
\R\tau^{\prime}_*:\quad  {\rm Perf}(X_{\FF})^{G_K}\to {\sd}({\rm Spa}(C)_{\proeet},\Q_{p})^{G_K}
$$
is fully faithful. That is, for $\sff,\sg\in {\rm Perf}(X_{\FF})^{G_K}$, the natural morphism in  ${\sd}({\rm Spa}(C)_{\proeet},\Q_p)^{G_K}$
$$\R\tau^{\prime}_*\R\Hhom_{\FF}(\sff,\sg)\to \R\Hhom_{\rm AVS}(\R\tau^{\prime}_*\sff,\R\tau^{\prime}_*\sg)
$$
is a quasi-isomorphism. 
\end{proposition}
\begin{remark} \label{Fontaine} In \cite{FoC} Fontaine has introduced the category $\scc(G_K)$ of almost C-representations.   Its derived $\infty$-category is equivalent to that of perfect complexes on the Fargues-Fontaine curve (see \cite[Th. 4.23]{Li}). Hence by Proposition \ref{recall1a} and Proposition \ref{old2}
it is a full subcategory of the ATVS's and AVS's. 
\end{remark}

We will also  need the following fully-faithfulness results: 
\begin{lemma}\label{last1}
Let \(S=\Spa(C)\) and let
$\mathcal F,\sg \in \sd(S_{\proeet},\Q_p)^{G_K}. $ 
There is a natural quasi-isomorphism in \(\sd(\Q_p)\)
\[
\R\operatorname{Hom}_{\mathrm{AVS}}
   (\mathcal F,\mathcal G)
\xrightarrow{\ \sim\ }
\R\underline{\operatorname{Hom}}_{\mathrm{ATVS}}
   (\R\pi_*\mathcal F,\R\pi_*\mathcal G)(*),
\]
where \(*\) denotes the one-point profinite set. This
quasi-isomorphism is compatible with restrictions from \(G_K\) to
\(G_L\), for every finite extension \(L/K\).
\end{lemma}

\begin{proof}
Set $
\mathcal H:=
\R\mathcal{H}om_{\mathrm{AVS}}(\mathcal F,\mathcal G). $
By Proposition~\ref{recall1a}, there is a natural quasi-isomorphism
in \(\underline{\sd}(S,\Q_p)^{G_K}\)
\[
\R\pi_*\mathcal H
\xrightarrow{\ \sim\ }
\R\mathcal{H}om_{\mathrm{ATVS}}
   (\R\pi_*\mathcal F,\R\pi_*\mathcal G).
\]
It follows from the definition of the global $\Hom$ in ATVS's that
\[
\R\underline{\operatorname{Hom}}_{\mathrm{ATVS}}
   (\R\pi_*\mathcal F,\R\pi_*\mathcal G)(*)
\simeq
\R\Gamma_{G_K}(S,\R\pi_*\mathcal H)(*).
\]

For \(n\geq 0\), set $Y_n:=S\times_{\Spa(K)}\widehat G_K^{\,n},$
and let \(\mathcal H_n\) denote the component of \(\mathcal H\) over
\(Y_n\). Since limits of condensed complexes are computed pointwise,
evaluation at \(*\) commutes with the limit over \(\Delta\). Moreover,
by the definition of \(\R\pi_{Y_n,*}\), we have
\[
(\R\pi_{Y_n,*}\mathcal H_n)(Y_n)(*)
\simeq
\R\Gamma(Y_{n,\proeet},\mathcal H_n).
\]
Hence we have: 
\[
\begin{aligned}
\R\Gamma_{G_K}(S,R\pi_*\mathcal H)(*)
&\simeq
\left(
\lim_{[n]\in\Delta}
(\R\pi_{Y_n,*}\mathcal H_n)(Y_n)
\right)(*)                                                     \simeq
\lim_{[n]\in\Delta}
\R\Gamma(Y_{n,\proeet},\mathcal H_n)         \\ &                    \simeq
\R\operatorname{Hom}_{\mathrm{AVS}}
   (\mathcal F,\mathcal G),
\end{aligned}
\]
as wanted.

All the constructions above are levelwise and commute with the
restriction maps associated with \(G_L\subset G_K\). This proves the
compatibility with restriction. In particular, since evaluation at
\(*\) commutes with filtered colimits, the same comparison holds after
passing to smooth Hom's.
\end{proof}
  
\begin{corollary}\label{deszcz2}
Let $S\in \mathrm{sPerf}_C$.  Let $\sff\in {\rm Perf}(X_{\FF,S^{\flat}})$. The natural morphisms in $\underline{\sd}(S,\Q_{p,\Box})$
\begin{align*}
\R\tau_*\R\Hhom_{\FF,S^{\flat}}(\sff,i_{\infty,*}F^a\mathbb{B}_{\pdr}(S)) & \to \R\Hhom_{{\rm TVS(S)}}(\R\tau_*\sff,\R\tau_*i_{\infty,*}F^a\mathbb{B}_{\pdr}(S)),\quad a\in\Z\cup\{-\infty\}, \\
\R\tau_*\R\Hhom_{\FF,S^{\flat}}(\sff,\mathcal{B}_{\rm pFF}) & \to \R\Hhom_{{\rm TVS}(S)}(\R\tau_*\sff,\R\tau_*\mathcal{B}_{\rm pFF})
\end{align*}
are quasi-isomorphisms. 
\end{corollary}
\begin{proof} We will give an argument for  the first map (the argument for the second map is analogous). 
We may pass to the condensed setting, that is, we want to show that, for  $\sff \in {\rm Perf}(X_{\FF, S^{\flat}})$,   the canonical map: 
$$
\R\tau_*\R\Hhom_{\FF,S^{\flat}}(\sff,i_{\infty,*}F^a\mathbb{B}_{\pdr}(S))\to \R\Hhom_{{\rm TVS}(S)}(\R\tau_*\sff,\R\tau_*i_{\infty,*}F^a\mathbb{B}_{\pdr}(S))
$$
is a quasi-isomorphism in $\underline{\sd}(S,\Q_{p})$. For that, consider its factorization:
  \begin{align*}
\R\pi_*\R\tau^{\prime}_*\R\Hhom_{{\FF,S^{\flat}}}(\sff,i_{\infty,*}F^a\mathbb{B}_{\pdr}(S)) &\stackrel{\sim}{ \to }
\R\pi_*\R\Hhom_{{\rm VS}(S)}(\R\tau^{\prime}_*\sff,\R\tau^{\prime}_*i_{\infty,*}F^a\mathbb{B}_{\pdr}(S))\\
 & \to 
\R\Hhom_{{\rm TVS}(S)}(\R\tau_*\sff,\R\tau_*i_{\infty,*}F^a\mathbb{B}_{\pdr}(S)).
\end{align*}
The first morphism is a quasi-isomorphism by Lemma \ref{deszcz1}. The first morphism is a quasi-isomorphism by
Lemma~\ref{deszcz1}, and the second one by
Proposition~\ref{recall1}.
   \end{proof}

\begin{remark} \label{duck1r} Consider the "evaluation on $\mathrm{Spa}(C)$" functor  from ATVS's to Galois equivariant solid $\Q_p$-modules:
$$
\mathrm{ev}_{C,*}: \quad \underline{\sd}(\mathrm{Spa}(C),\Q_p)^{G_K}\to \Q_{p,\Box}[G_K],\quad \sff\mapsto \sff(\mathrm{Spa}(C)),
$$
and the induced morphism
$$
\R\uHom_{\mathrm{ATVS}}(\sff_1,\sff_2)\to \R\uHom_{G_K}(\mathrm{ev}_{,*}\sff_1,\mathrm{ev}_{C,*}\sff_2).
$$
We will show in the Appendix \ref{firstA} that this morphism,   in general, is not a quasi-isomorphism. However, in the case when $\sff_1=\underline{\Q}_p$ it clearly is.
\end{remark}

\section{Conjecture $C_{\st}$}
We now descend  the curve-level  internal $C_{\st}$-comparison isomorphisms from Theorem \ref{essen3} and Theorem \ref{luminy20} of Section \ref{better1} to Arithmetic Topological Vector Spaces.
This is accomplished via Proposition \ref{film1} below, which extends the perfect fully-faithfulness of Corollary \ref{deszcz2}  to the  non-perfect coefficient objects used in the main theorem.
Taking equivariant global sections and using the Galois acyclicity of $\B_{\pdr}$ and $\B_{\pFF}$ then yields the desired $C_{\st}$-conjecture.

\subsection{Descent to Arithmetic Topological Vector Spaces}
Our final fully-faithfulness result lifts Corollary \ref{deszcz2} from  perfect complexes  to  pro-\'etale sheaves. 
\begin{proposition}[Fully-faithfulness for pro-\'etale sheaves] \label{film1} Let $X$ be a smooth partially proper rigid analytic variety over $K$ of dimension $d$. Let $a\in\Z\cup\{-\infty\}$ and let $\spp\in \{i_{\infty,*}F^a{\B}_{\pdr}, \mathcal{B}_{\rm pFF},\so\}$. 
The natural morphism in $\underline{\sd}({\rm Spa}(C),\Q_{p,\Box})^{G_K}$ (we set $\se_{\proeet,c}:=\se_{\proeet,c}(X_C,\Q_p(d))[2d]$) 
$$
\R\tau_*\R\Hhom_{\rm FF}(\se_{\proeet,c},\spp)  \to   \R\Hhom_{\rm ATVS}(\R\tau_*\se_{\proeet,c},\R\tau_*\spp)
$$
is a  quasi-isomorphism. 
\end{proposition}
\begin{proof} We start with the de Rham claim. Take $r\geq 2d$. Let $S\in \mathrm{sPerf}_C$ and set $\se_{r}:=\se_{\proeet,c}(X_S,\Q_p(r))[2d]$. It suffices to show that the natural morphism in ${\rm TVS}(S)$'s
$$\R\tau_*\R\Hhom_{\FF,S^{\flat}}(\se_{r},i_{\infty,*}F^{a+r-d}\mathbb{B}_{\pdr}(S))\to \R\Hhom_{{\rm TVS}(S)}(\R\tau_*\se_{r},\R\tau_*i_{\infty,*}F^{a+r-d}\mathbb{B}_{\pdr}(S))$$
is a quasi-isomorphism. 
For that, using  the distinguished triangle (see \eqref{paris20})
\begin{equation}\label{kolobrzeg1}
\se_{\proeet,c}(X_S,\Q_p(r))\to\se_{\hk,c}(X_S,r)\to\se_{\dr,c}(X_S,r),
\end{equation}
we see  that  it suffices to show that, for any $i\in\Z$,   the natural morphism in ${\rm TVS}$'s   
\begin{align}\label{tea1}
\R\tau_*\R\Hhom_{\FF,S^{\flat}}(\se,i_{\infty,*}F^i\mathbb{B}_{\pdr}(S)) & \to \R\Hhom_{\rm TVS(S)}(\R\tau_*\se, \R\tau_*i_{\infty,*}F^i\mathbb{B}_{\pdr}(S))
\end{align}
is a  quasi-isomorphism, where $\se=\se_{\hk,c}:=\se_{\hk,c}(X_S,r)$ or $\se=\se_{\dr,c}:=\se_{\dr,c}(X_S,r)$.

For $\se=\se_{\hk,c}$,  assume first that the Hyodo--Kato cohomology of $X_C$ is of finite rank (this is always true locally).  Then $\se_{\hk,c}$ is a perfect complex 
and   we can invoke   Corollary \ref{deszcz2}. The general case follows by writing $\se_{\hk,c}$ as a countable colimit of perfect complexes and using the fact that $\R\tau_*$ commutes with colimits and countable limits of nuclear sheaves (see \cite[Prop. 5.4.12]{ALBM}). 
For $\se=\se_{\dr,c}$ in \eqref{tea1}, arguing analogously, we may assume  that $X$ is Stein. Then, using \eqref{paris5}, it suffices to show that the natural morphism in ${\rm TVS}$'s
\begin{equation}\label{doing1}
\R\tau_*\R\Hhom_{\FF, S^{\flat}}(i_{\infty,*}(W\otimes^{\LL_{\Box}}_KR_S),\spp)\to  \R\Hhom_{{\rm TVS}(S)}(\R\tau_*i_{\infty,*}(W\otimes^{\LL_{\Box}}_KR_S),\R\tau_*\spp).
\end{equation}
for $W$ -- a space of compact type over $K$ and $\spp=i_{\infty,*}F^i\mathbb{B}_{\pdr}(S)$, is a quasi-isomorphism.  Below we will also allow $\spp=\mathcal{B}_{\rm pFF}$ or $\spp=\so$. 

    Assume first that $W$ is Smith. We compute in ${\rm TVS}(S)$'s
\begin{align*}
\R\tau_*\R\Hhom_{\FF,S^{\flat}}(i_{\infty,*}(W\otimes^{\LL_{\Box}}_KR_S),\spp) & \simeq \R\tau_*(W^*\otimes^{\LL_{\Box}}_K\R\Hhom_{\FF,S^{\flat}}(i_{\infty,*}R_S,\spp))\\
& \simeq W^*\otimes^{\LL_{\Box}}_K \R\tau_*\R\Hhom_{\rm FF}(i_{\infty,*}R_S,\spp));\\
\R\Hhom_{{\rm TVS}(S)}(\R\tau_*i_{\infty,*}(W\otimes^{\LL_{\Box}}_KR_S),\R\tau_*\spp) & \simeq  \R\Hhom_{{\rm TVS}(S)}(W\otimes^{\LL_{\Box}}_K \R\tau_*i_{\infty,*}R_S,\R\tau_*\spp)\\
&\stackrel{\simeq}{\leftarrow} W^*\otimes^{\LL_{\Box}}_K\R\Hhom_{{\rm TVS}(S)}(\R\tau_*i_{\infty,*}R_S,\R\tau_*\spp).
\end{align*}
The above quasi-isomorphisms are compatible with the morphism \eqref{doing1}, 
hence we are reduced to showing that the natural morphism in ${\rm TVS}(S)$'s
$$
\R\tau_*\R\Hhom_{\FF,S^{\flat}}(i_{\infty,*}R_S,\spp)\to \R\Hhom_{{\rm TVS}(S)}(\R\tau_*i_{\infty,*}R_S,\R\tau_*\spp)
$$
is a quasi-isomorphism. But we have shown this in Corollary \ref{deszcz2}.  For a general compact type $W$, write it as a colimit of Smith spaces: $W=\colim_nW_n$ and reduce to the Smith case using the fact that $\R\tau_*$ commutes with colimits and countable limits of nuclear sheaves.

We pass now to the Hyodo--Kato  claim. Take $r\geq 2d$. Let $S\in \mathrm{sPerf}_C$. We want  to show that the natural morphism in ${\rm TVS}(S)$'s
$$\R\tau_*\R\Hhom_{\FF,S^{\flat}}(\se_{\proeet,c},\mathcal{B}_{\rm pFF})\to \R\Hhom_{{\rm TVS}(S)}(\R\tau_*\se_{\proeet,c},\R\tau_*\mathcal{B}_{\rm pFF})$$
is a quasi-isomorphism. 
For that, using  the distinguished triangle  \eqref{kolobrzeg1}, 
we see  that  it suffices to show that  the natural morphisms in ${\rm TVS}(S)$'s
\begin{align*}
\R\tau_*\R\Hhom_{\FF,S^{\flat}}(\se,\mathcal{B}_{\rm pFF}) & \to \R\Hhom_{{\rm TVS}(S)}(\R\tau_*\se,\R\tau_*\mathcal{B}_{\rm pFF})
\end{align*}
is a quasi-isomorphism for $\se$ as above. 

For $\se=\se_{\hk,c}$, arguing as above, we may assume that   the Hyodo--Kato cohomology of $X_C$ is of finite rank.  Then   $\se_{\hk,c}$ is a perfect complex and we can use Corollary \ref{deszcz2}. 
For $\se=\se_{\dr,c}$, we may assume  that $X$ is Stein.  By Lemma  \ref{paris13a}, the left hand side is trivial. Hence it suffices to show that so is the right hand side. 
For that, using \eqref{paris5}, it suffices to show that  in ${\rm TVS}(S)$'s
\begin{equation}\label{doing1a}
\R\Hhom_{{\rm TVS}(S)}(\R\tau_*i_{\infty,*}(W\otimes^{\LL_{\Box}}_KR_S),\R\tau_*\mathcal{B}_{\rm pFF})\simeq 0,
\end{equation}
for $W$ -- a space of compact type over $K$. But we may assume that $W$ is Smith and, by the above, 
we have  in ${\rm TVS}(S)$'s
\begin{align*}
\R\Hhom_{{\rm TVS}(S)}(\R\tau_*i_{\infty,*}(W\otimes^{\LL_{\Box}}_KR_S),\R\tau_*\mathcal{B}_{\rm pFF})  
&\simeq W^*\otimes^{\Box}_K\R\Hhom_{{\rm TVS}(S)}(\R\tau_*i_{\infty,*}R_S,\R\tau_*\mathcal{B}_{\rm pFF})\\
&\simeq W^*\otimes^{\Box}_K\R\Hhom_{{\rm TVS}(S)}(\mathbb{G}_a,\mathbb{B}_{\rm pFF}) \simeq 0.
\end{align*}
The second quasi-isomorphism follows from    Proposition \ref{app1}.

Finally, for $\spp=\so$,  we argue as above using 
for the Hyodo--Kato part fully-faithfulness for perfect complexes if the de Rham cohomology has finite rank.
\end{proof}

\subsection{Conjecture $C_{\st}$}
We are now ready to prove the main theorem of this paper. 

We start with recalling briefly the definition of pro-\'etale cohomology $\rm TVS$-presheaves from \cite[Sec. 6.1]{CGN2}.
Let $X$ be a smooth partially proper  rigid analytic variety over $K$.  
For $?=-,c$, we define the enriched presheaves 
$$
 {\mathbb R}_{\proeet,?}(X_C,\Q_p)\in \underline{\sd}({\rm Spa}(C),\Q_{p,\Box}):\quad S\to \R\Gamma_{\proeet,?}(X_{S}, \Q_p).
 $$
A priori, the topology on pro-\'etale cohomology is induced via \v{C}ech procedure from $p$-adic topologies but we have a natural quasi-isomorphism in $\underline{\sd}({\rm Spa}(C),\Q_{p,\Box})$
$$
{\mathbb R}_{\proeet,?}(X_C,\Q_p)\simeq \R\pi_*{\mathbb R}^{\rm alg}_{\proeet,?}(X_C,\Q_p),
$$
where  the sheaf $ {\mathbb R}^{\rm alg}_{\proeet,?}(X_C,\Q_p)$ is the
algebraic version of ${\mathbb R}_{\proeet,?}(X_C,\Q_p)$ (defined as ${\mathbb R}_{\proeet,?}(X_C,\Q_p)$ but with discrete topology objectwise). 
It was shown  in \cite[Prop. 6.2]{CGN2}, via comparison theorems,  that if $X$ is Stein,  there exists a natural quasi-isomorphism in $\underline{\sd}({\rm Spa}(C),\Q_{p,\Box})$
$$
\R\tau_*\se_{\proeet,?}(X_C,\Q_p) \simeq {\mathbb R}_{\proeet,?}(X_C,\Q_p),\quad ?=-,c.
$$
By taking colimits this also holds for the compactly supported cohomology and general smooth partially proper $X$. 

 For $r\in\Z$, we can define in a similar way ${\mathbb R}_{\proeet,?}(X,\Q_p(r))\in \underline{\sd}({\rm Spa}(K),\Q_{p,\Box})$ such that  ${\mathbb R}_{\proeet,?}(X_C,\Q_p(r))\simeq \varepsilon^*{\mathbb R}_{\proeet,?}(X,\Q_p(r))$.  Hence
 ${\mathbb R}_{\proeet,?}(X_C,\Q_p(r))\in \underline{\sd}({\rm Spa}(C),\Q_{p,\Box})^{G_K}$ by Lemma \ref{berlin1}. We have analogous constructions in the algebraic setting and all the above properties hold. 

  Our main theorem is: 
\begin{theorem}[Conjecture $C_{\st}$] \label{duck1}Let $X$ be a smooth partially proper rigid analytic variety over $K$ of dimension $d$.
\begin{enumerate}
\item {\rm (Filtered de Rham cohomology.)} There is a natural quasi-isomorphism\footnote{We defined the augmented filtered derived category:  
$
\sff^{\mathrm{aug}}\sd(K_\square)
:=
\mathrm{Fun}\bigl(
(\mathbf Z\cup\{-\infty\})^{\rm op},
\sd(K_\square)
\bigr). $}  in
\(\sff^{\mathrm{aug}}\sd(K_\square)\)
\[
F^\bullet \R\Gamma_{\mathrm{dR}}(X)
\simeq
\R{\uHom}_{\mathrm{ATVS}}
\bigl(
\mathbb{R}_{\mathrm{pro\acute et},c}
  (X_C,\mathbf Q_p(d))[2d],
F^\bullet \mathbb{B}_{\mathrm{pdR}}
\bigr).
\]
Here the right-hand side denotes the filtered object whose
\(a\)-th term is
\[
\R\uHom_{\mathrm{ATVS}}
\bigl(
\mathbb{R}_{\mathrm{pro\acute et},c}
  (X_C,\mathbf Q_p(d))[2d],
F^a\mathbb{B}_{\mathrm{pdR}}
\bigr),
\qquad
a\in\mathbf Z\cup\{-\infty\}.
\]
\item  {\rm (Hyodo--Kato  cohomology.)} Assume that de Rham  cohomology of $X$ has finite rank. There is a natural quasi-isomorphism in $\sd^{\prime}_{\phi,N, G_K}(C^{\nr}_{\Box})$ 
$${\R}\Gamma_{\hk}(X_C)  \simeq\R\uHom^{\mathrm{sm}}_{\rm ATVS}({\mathbb R}_{\proeet,c}(X_C,\Q_p(d))[2d],{\mathbb B}_{\rm pFF}).
$$
\end{enumerate}
\end{theorem}

  As an immediate consequence of Theorem \ref{duck1} we obtain the following: 
 \begin{corollary}[Conjecture $C_{\st}$ for Hodge cohomology] \label{gateway1} Let $X$ be a smooth partially proper rigid analytic variety over $K$ of dimension $d$. Let $a\in\Z$. 
 We have in $\sd(K_{\Box})$
\begin{equation}\label{gateway1a}
\mathrm{gr}^a_F{\R}\Gamma_{\dr}(X)  \simeq\R\uHom_{\rm ATVS}({\mathbb R}_{\proeet,c}(X_C,\Q_p(d))[2d], \mathbb{G}_a(a)[\log t]).
\end{equation}
In particular, if $X$ is Stein, we have in $\sd(K_{\Box})$
$$\Omega^a(X)[-a]  \simeq\R\uHom_{\rm ATVS}({\mathbb R}_{\proeet,c}(X_C,\Q_p(d))[2d], \mathbb{G}_a(a)[\log t]).
$$
Moreover, the de Rham differential $
d:\Omega^a\to \Omega^{a+1} $ is then induced by the boundary map from the distinguished triangle
\begin{equation}\label{nie1}
\mathbb{G}_a(a+1)[\log t]\to (F^a/F^{a+2}){\mathbb B}_{\pdr}\to \mathbb{G}_a(a)[\log t]\stackrel{\partial}{\to} \mathbb{G}_a(a+1)[\log t][1]
\end{equation}
 \end{corollary}
 \begin{proof} Apply \(\mathrm{gr}_F^a\) to the filtered quasi-isomorphism of
Theorem~\ref{duck1}(1) to get the first claim of the corollary. For the differential,  apply Theorem \ref{duck1} to  the distinguished triangle
 \eqref{nie1} 
  to yield the distinguished triangle
  $$
 \Omega^{a+1}[-a-1]\to  (F^a/F^{a+2})\R\Gamma_{\dr}(X)\to \Omega^a[-a]\stackrel{d}{\to}\Omega^{a+1}[-a].
  $$
 \end{proof}
 \begin{remark} In Section \ref{soso12} we will show (conditionally) that the quasi-isomorphism \eqref{gateway1a} is compatible with multiplicative structures. 
 \end{remark}

\begin{proof}[Proof of Theorem \ref{duck1}] Let $a\in\Z\cup \{-\infty\}$. Write for simplicity ${\mathbb R}_{\proeet,c}:={\mathbb R}_{\proeet,c}(X_C,\Q_p(d))[2d]$ and $\se_{\proeet,c}:=\se_{\proeet,c}(X_C,\Q_p(d))[2d]$. We start with the de Rham case. Assume first that $X$ is Stein with finite rank de Rham cohomology. Since $\R\tau_*i_{\infty,*}F^a{\B}_{\pdr}\simeq F^a{\mathbb B}_{\pdr}$, we compute in ${\rm ATVS}$'s
\begin{align}\label{dzisiaj1}
\R\Hhom_{\rm ATVS}({\mathbb R}_{\proeet,c}, F^a{\mathbb B}_{\pdr}) & \simeq 
\R\Hhom_{\rm ATVS}(\R\tau_*\se_{\proeet,c},\R\tau_*i_{\infty,*}F^a{\B}_{\pdr}) \\
 & \stackrel{\simeq}{\leftarrow}  \R\tau_*\R\Hhom_{\rm FF}(\se_{\proeet,c},i_{\infty,*}F^a{\B}_{\pdr})\notag 
\end{align}
 The second   quasi-isomorphism follows   from Proposition \ref{film1}.
 Taking  global sections of   this computation yields in $\sd(K_{\Box})$
\begin{align}\label{lion1}
\R\uHom_{\rm ATVS}({\mathbb R}_{\proeet,c}, F^a{\mathbb B}_{\pdr}) & \simeq \R\Gamma_{G_K}({\rm Spa}(C), \R\Hhom_{\rm ATVS}({\mathbb R}_{\proeet,c}, F^a{\mathbb B}_{\pdr}) )\\
 & \stackrel{\simeq}{\leftarrow} \R\Gamma_{G_K}({\rm Spa}(C),\R\tau_*\R\Hhom_{\rm FF}(\se_{\proeet,c},i_{\infty,*}F^a{\B}_{\pdr})) \notag\\
  & \simeq  \R\Gamma(G_K, \R\Hom_{\rm FF}(\se_{\proeet,c},i_{\infty,*}F^a{\B}_{\pdr}))\notag\\
  & \simeq   F^a{\R\Gamma}_{\dr}(X),\notag
\end{align}
as wanted.
 Here the fourth quasi-isomorphism follows from   Theorem \ref{essen3}; the third one  from Proposition \ref{lozko1}:
 the d\'evissage in the proof of Proposition \ref{film1}  shows that the sheaf 
$$ \mathcal H_a:=
\R\mathcal Hom_{\mathrm{FF}}
(\mathcal E_{\mathrm{pro\acute et},c},i_{\infty,*}F^a\B_{\mathrm{pdR}}) $$
is nuclear: the Smith-stage Homs are nuclear, and nuclear objects are stable under countable limits. Hence \(R\tau_*\mathcal H_a\) is derived rigid by \cite[Rem.  4.17]{TVS}.
 
 For the general case, we use analytic hypercovers with Stein spaces having finite rank de Rham cohomologies, use quasi-isomorphisms \eqref{lion1} on finite levels, and pass to limit.  Moreover, all  the morphisms in the above computations  are natural in the
coefficient object. Hence, applied to the filtered object
\(F^\bullet \mathbb{B}_{\mathrm{pdR}}\), they define a morphism of filtered
objects. Its component at every
\(a\in\mathbf Z\cup\{-\infty\}\) is a quasi-isomorphism  hence it is a quasi-isomorphism in the filtered
derived category.
  
     We pass now to the Hyodo--Kato case.   Since $\R\tau_*\mathcal{B}_{\rm pFF}\simeq {\mathbb B}_{\rm pFF}$, we compute in ${\rm ATVS}$'s 
\begin{align}\label{dzisiaj2}
\R\Hhom_{\rm ATVS}({\mathbb R}_{\proeet,c}, {\mathbb B}_{\rm pFF}) & \simeq 
\R\Hhom_{\rm ATVS}(\R\tau_*\se_{\proeet,c},\R\tau_*\mathcal{B}_{\rm pFF}) \\
 & \stackrel{\simeq}{\leftarrow} \R\tau_*\R\Hhom_{\rm FF}(\se_{\proeet,c},\mathcal{B}_{\rm pFF}),\notag \\
  & \simeq \colim_{n,m,L}\R\tau_*\R\pi_*(\R\Gamma_{\hk}(X_C)\otimes^{\LL_{\Box}}_{C^{\nr}}\so^+_{n,m,L}).\notag
\end{align}
The second quasi-isomorphism follows from Proposition \ref{film1};  the third one --   from Theorem \ref{luminy20}. 
Taking  global sections of   this computation yields in $\sd(F_{\Box})$
\begin{align*}
\R\uHom_{\rm ATVS}({\mathbb R}_{\proeet,c}, {\mathbb B}_{\rm pFF}) & \simeq \R\Gamma_{G_K}({\rm Spa}(C), \R\Hhom_{\rm ATVS}({\mathbb R}_{\proeet,c}, {\mathbb B}_{\rm pFF}) )\\
& \simeq  \R\Gamma_{G_K}({\rm Spa}(C), \colim_{n,m,L }\R\tau_*\R\pi_*(\R\Gamma_{\hk}(X_C)\otimes^{\LL_{\Box}}_{C^{\nr}}\so^+_{n,m,L})) \\
& \simeq  \R\Gamma(G_K, \colim_{n,m,L}\R\Gamma(X_{\FF},\R\pi_*(\R\Gamma_{\hk}(X_C)\otimes^{\LL_{\Box}}_{C^{\nr}}\so^+_{n,m,L}))\\
& \simeq  \R\Gamma(G_K, \R\Gamma_{\hk}(X_C)\otimes^{\LL_{\Box}}_{C^{\nr}}\B_{\rm pFF}).
\end{align*}
Here the third quasi-isomorphism follows from Proposition \ref{lozko1} (since, clearly, $\R\tau_*\R\pi_*(\R\Gamma_{\hk}(X_C)\otimes^{\LL_{\Box}}_{C^{\nr}}\so^+_{n,m,L})$ is derived rigid). 

Varying $K$ over all finite extensions of $K$, this implies in $\sd^{\prime}_{\phi,N,G_K}(C^{\nr}_{\Box})$
\begin{align*}
\R\uHom^{\mathrm{sm}}_{\rm ATVS}({\mathbb R}_{\proeet,c}, {\mathbb B}_{\rm pFF}) & \simeq ( \R\Gamma_{\hk}(X_C)\otimes^{\LL_{\Box}}_{C^{\nr}}\B_{\rm pFF})^{\R G_K-{\rm sm}} \\
& \simeq  {\R\Gamma}_{\hk}(X_C),
\end{align*}
as wanted. 
\end{proof}
\begin{remark}[Uniqueness of Hyodo--Kato cohomology]
The proof of Theorem \ref{duck1}  provides a uniqueness criterion for Hyodo--Kato cohomology.
 Suppose that a second Hyodo--Kato theory, together with its compactly supported version, satisfies the same basic comparison theorem and compatible Poincar\'e duality, and that the relevant complexes are perfect locally for the analytic topology.
 The reconstruction argument then identifies the two theories locally.
 Provided these identifications are compatible with restriction and both theories satisfy analytic descent, they yield a natural equivalence globally.

This gives a possible approach to comparing the theory used here with the theory associated with the Hyodo--Kato stack of \cite{ABL}; applying the criterion requires establishing the corresponding comparison theorem and compatibilities.
\end{remark}

\subsection{Examples}
To illustrate  Theorem \ref{duck1}, let us look  in more detail at two simple  examples. 

\subsubsection{Stein curves}
We start with   the case of  smooth Stein curves and de Rham cohomology\footnote{The calculation below is intended as a consistency check rather than an independent proof of Theorem \ref{duck1}; we indicate the relevant d\'evissage and suppress some routine compatibility checks.}
. This calculation was one of the original motivations for the compact-support $C_{\st}$-conjecture  (see \cite{CGN3}).

Let $X$ be a smooth Stein rigid analytic curve defined over a finite extension $K$. 
Assume  that it is geometrically connected and that the de Rham cohomology groups 
$H^i_{\rm dR}(X)$ are of finite rank over $K$. We note that  they are always trivial  for $i\geq 2$. 

Comparison theorems for pro-\'etale cohomology with compact support yield the following (see \cite[Sec. 7.2]{AGN}),  where we  use the shorthand 
$\mathbb{H}^i_c:=\mathbb{H}^i_{{\proeet},c}(X_{C},\Q_p(1))$: 
\begin{enumerate}[leftmargin=*]
\item We have isomorphisms in ATVS's:
\begin{align*}
\mathbb{H}^i_c=\begin{cases}    \mathbb{V} & \text{ if  } i=1,\\
0 & \text{ if  $i\neq 1,2$},
\end{cases}
\end{align*}
where  $\mathbb{V}:=(H^1_{{\rm HK},c}(X_{C})\otimes^{\Box}_{C^{\nr}} t\mathbb{B}^+_{\FF})^{N=0,\varphi=p}$  
is finite rank over  $\Q_p$  and  de Rham
as a representation of  $G_K$. 
\item  We have exact sequences (see \cite[(7.9), (7.12), (7.13)]{AGN})
\begin{align*}
&    0\to \mathbb{V}_1 \to   \mathbb{H}_c^2\to \Q_p\to 0,\\
&   0\to \mathbb{W}/\mathbb{V}\to H^1_{c}(X,\so)\otimes_K^{\Box}\mathbb{G}_a\to \mathbb{V}_1\to 0,
\end{align*}
where $\mathbb{W}:=(H^1_{{\rm HK},c}(X_{C})\otimes^{\Box}_{C^{\nr}} \mathbb{B}^+_{\FF})^{N=0,\varphi=p}$.
\end{enumerate}

Set $\mathbb{R}_{{\proeet},c}:=\mathbb{R}_{{\proeet},c}(X_{C},\Q_p(1))[2]$. Consider the spectral sequence
$$
E^{a,b}_2=\underline{\Ext}^a_{\rm ATVS}(\mathbb{R}^{-b}_{{\proeet},c},\mathbb{B}_{\pdr})\Rightarrow H^{a+b}(\underline{{\R\Hom}}_{\rm ATVS}(\mathbb{R}_{{\proeet},c},\mathbb{B}_{\pdr}))
$$
Set  ${\Ext}^i(-):=\underline{{\Ext}}^i_{\mathrm{ATVS}}(-,\mathbb{B}_{\pdr})$ and $\Ext^i:=H^{i}(\underline{{\R\Hom}}_{\rm ATVS}(\mathbb{R}_{{\proeet},c},\mathbb{B}_{\pdr}))$. By Theorem \ref{duck1} we have in $\sd(K_{\Box})$ the isomorphism
$
\Ext^i\simeq H^i_{\dr}(X),$ for $ i\in \Z. 
$
We would like to see how it fits with the above computations.  

Since $\Ext^i=0,$ for $i\geq 2$. For $i=0,1$, 
the  above  spectral sequence yields exact sequences (note that the relevant $\Ext^2(-)=0$)
\begin{equation}\label{rhom}
0\to {\rm Ext}^1(\mathbb{H}^{3-i}_c)\to {\Ext}^i\to {\rm Ext}^0(\mathbb{H}^{2-i}_c)\to 0
\end{equation}
Since  ${\rm Ext}^0(\mathbb{G}_a)=0$,  ${\rm Ext}^1(\mathbb{G}_a)=0$, and $H^1_c(X,\so)$ is of compact type, 
we have 
\begin{align*}
& {\Ext}^0(\mathbb{V}_1)=0,\quad {\rm Ext}^1(\mathbb{V}_1)={\rm Ext}^0(\mathbb{W}/\mathbb{V}),\\ 
& {\Ext}^0(H^1_c(X,\so)\otimes_K^{\LL_{\Box}}\mathbb{G}_a)=0,\quad {\Ext}^1(H^1_c(X,\so)\otimes_K^{\LL_{\Box}}\mathbb{G}_a)=0
\end{align*}
Since ${\Ext}^0(\mathbb{H}_c^2)  \simeq K\simeq H^0_{\rm dR}(X) $, this gives
$
H^0_{\rm dR}(X) \stackrel{\sim}{\to}\Ext^0,
$
as wanted. 

For $\Ext^1$, we will use the exact sequence
\begin{equation}\label{goraco2}
0\to {\rm Ext}^1(\mathbb{H}^{2}_c)\to {\Ext}^1\to {\rm Ext}^0(\mathbb{H}^{1}_c)\to 0
\end{equation}
To compute the side terms, first, a small digression.  Recall that we have Hyodo--Kato isomorphism $H^i_{\rm dR,?}(X)\simeq ( H^i_{\rm HK,?}(X_{C})\otimes_{C^{\nr}}^{\Box}C)^{G_K}$. The slopes of  Frobenius on 
$H^1_{\rm HK}(X_{C})$ and  $H^1_{{\rm HK},c}(X_{C})$ are in the interval 
$[0,1]$;  denote by  $H^1_{{\rm HK},?}(X_{C})^{(r)}$
the subspace of slope $r$ and set 
$H^1_{{\rm dR},?}(X)^{(r)}:=(H^1_{{\rm HK},?}(X_{C})^{(r)}\otimes_{C^{\nr}}^{\Box}C)^{G_K}$.
The duality between  $H^1_{\hk}(X)$ and  $H^1_{{\hk},c}(X)$ induces duality
between $H^1_{\rm dR}(X)^{(r)}$ and  $H^1_{{\rm dR},c}(X)^{(1-r)}$.  We also have 
$\mathbb{V}\simeq t(H^1_{{\rm HK},c}(X_{C})^{(0)}\otimes^{\Box}_{C^{\nr}}\mathbb{B})^{\varphi=1}$.

To show that $\Ext^1\simeq H^1_{\dr}(X)$, it suffices to show that
\begin{align*}
{\Ext}^0(\mathbb{H}_c^1) & \simeq 
(H^1_{{\rm dR},c}(X)^{(0)})^\dual \simeq H^1_{\rm dR}(X)^{(1)}
\subset H^1_{\rm dR}(X),\\
{\Ext}^1(\mathbb{H}_c^2) & \simeq {\rm Ker}\big(H^1_{\rm dR}(X)\to H^1_{\rm dR}(X)^{(1)}\big)
\end{align*}
The first isomorphism follows from the fact that $\mathbb{H}^1_c\simeq \mathbb{V}$ is de Rham associated to the pair $(H^1_{\hk,c}(X_C)^{(0)},H^1_{\dr,c}(X)^{(0)})$, where de Rham filtration is concentrated in degree $0$. 
For the second one, we use that
$$
\Ext^1(\mathbb{H}_c^2)\simeq \Ext^1(\mathbb{V}_1)\simeq \Ext^0(\mathbb{W}/\mathbb{V})\simeq \mathrm{Ker}\big(\Ext^0(\mathbb{W})\to\Ext^0(\mathbb{V})\big)
$$
and invoke\footnote{Forget first the $G_K$-action in \cite[Prop. 5.18]{CN5}. Since the resulting isomorphism is induced by the canonical evaluation pairing, it is functorial in $M$ and thus  $G_K$--equivariant.}  \cite[Prop. 5.18]{CN5}.   Hence the edge exact sequence \eqref{goraco2} identifies with
$$
0\to  {\rm Ker}\big(H^1_{\rm dR}(X)\to H^1_{\rm dR}(X)^{(1)}\big) \to H^1_{\dr}(X)\to H^1_{\dr}(X)^{(1)}\to 0
$$

\begin{remark} If we had used  $\mathbb{B}_{\dr}$ instead of  $\mathbb{B}_{\pdr}$, there would have been  parasitic terms coming from  
$\underline{{\Ext}}^1_{\rm ATVS}(\mathbb{V},\mathbb{B}_{\dr})$ and  $\underline{{\Ext}}^1_{\rm ATVS}(\Q_p,\mathbb{B}_{\dr})$ contributing respectively to 
$\underline{{\Ext}}^1_{\rm ATVS}(\mathbb{H}^1_c,\mathbb{B}_{\dr})$ and  $\underline{{\Ext}}^1_{\rm ATVS}(\mathbb{H}^2_c,\mathbb{B}_{\dr})$.
\end{remark}

\subsubsection{Proper rigid analytic varieties}
  Let $X$ be a smooth proper rigid analytic variety over $K$ of dimension $d$.
 Let $a\in \Z\cup\{-\infty\}$. We compute  in $\sd(K_{\Box})$:
\begin{align*}
F^a{\R}\Gamma_{\dr}(X)  & \simeq\R\uHom_{\rm ATVS}({\mathbb R}_{\proeet,c}(X_C,\Q_p(d))[2d], F^a{\mathbb B}_{\pdr})\\
 & \simeq \R\uHom_{G_K}({ \R}\Gamma_{\proeet,c}(X_C,\Q_p(d))[2d], F^a{\B}_{\pdr})\\
 & \simeq \R\Gamma(G_K, { \R}\Gamma_{\proeet}(X_C,\Q_p)\otimes^{\LL_{\Box}}_{\Q_p}F^a{\B}_{\pdr})
\end{align*}
 The first quasi-isomorphism follows from Theorem \ref{duck1},  the second one - from Remark \ref{duck1r}, and the third one -- from Poincar\'e duality for pro-\'etale cohomology and the fact that this cohomology is dualizable.  We recovered the known $(C_{\dr})$ statement.

 Similarly, for the Hyodo--Kato cohomology.  Assume that de Rham  cohomology of $X$ has finite rank. We compute  in $\sd^{\prime}_{\phi,N, G_K}(C^{\nr}_{\Box})$: 
\begin{align*}
{\R}\Gamma_{\hk}(X_C)  &  \simeq\R\uHom^{\mathrm{sm}}_{\rm ATVS}({\mathbb R}_{\proeet,c}(X_C,\Q_p(d))[2d],{\mathbb B}_{\rm pFF})\\
& \simeq \R\uHom^{\mathrm{sm}}_{G_K}({\R}\Gamma_{\proeet,c}(X_C,\Q_p(d))[2d],{ \B}_{\rm pFF})\\
 & \simeq ({\R}\Gamma_{\proeet}(X_C,\Q_p)\otimes^{\LL_{\Box}}_{\Q_p}{ \B}_{\rm pFF})^{\R G_K-{\rm sm}}.
\end{align*}
Again: The first quasi-isomorphism follows from Theorem \ref{duck1}, the second one -- from Remark \ref{duck1r}, and the third one -- from Poincar\'e duality for pro-\'etale cohomology and the fact that this cohomology is dualizable.
  We recovered the known $(C_{\st})$ statement.

Hence, for proper varieties, Theorem \ref{duck1}  is equivalent to the usual $C_{\dr}$-- and  $ C_{\st}$--reconstruction formulas after applying geometric Poincar\'e duality.
 In this case the derived Hom introduces no additional extension terms.
 This is in sharp contrast with the Stein case discussed above, where the higher $\Ext$ groups are necessary to recover the full de Rham complex.

\subsection{Algebraic Conjecture $C_{\st}$} 
Theorem \ref{duck1} has an algebraic analog though it   looses the topological structure on the de Rham cohomology. 
 \begin{corollary} [Algebraic Conjecture $C_{\st}$] \label{duck1a}Let $X$ be a smooth partially proper rigid analytic variety over $K$ of dimension $d$.
 \begin{enumerate}
\item {\rm (Filtered de Rham cohomology.)} There is a natural quasi-isomorphism in
\(\sff^{\mathrm{aug}}\sd(K)\)
\[
F^\bullet \R\Gamma_{\mathrm{dR}}(X)
\simeq
\R{\uHom}_{\mathrm{AVS}}
\bigl(
\mathbb{R}^{\rm alg}_{\mathrm{pro\acute et},c}
  (X_C,\mathbf Q_p(d))[2d],
F^\bullet \mathbb{B}_{\mathrm{pdR}}
\bigr).
\]
\item  {\rm (Hyodo--Kato  cohomology.)} Assume that $X$ has  finite rank de Rham cohomology.
 There is a natural quasi-isomorphism in $\sd^{\prime}_{\phi,N, G_K}(C^{\nr})$ 
$${\R}\Gamma_{\hk}(X_C)  \simeq\R\Hom^{\rm sm}_{\rm AVS}({\mathbb R}^{\rm alg}_{\proeet,c}(X_C,\Q_p(d))[2d],{\mathbb B}_{\rm pFF}).
$$
\end{enumerate}
 \end{corollary}
 \begin{proof} Let $a\in\Z\cup \{-\infty\}$. For the de Rham claim, by Theorem \ref{duck1}, it suffices to show that in $\sd(K)$ 
 $$
 \R\uHom_{\rm ATVS}({\mathbb R}_{\proeet,c}(X_C,\Q_p(d))[2d], F^a{\mathbb B}_{\pdr})(*)\simeq \R\Hom_{\rm AVS}({\mathbb R}^{\rm alg}_{\proeet,c}(X_C,\Q_p(d))[2d], F^a{\mathbb B}_{\pdr}).
 $$
But this follows immediately from Lemma \ref{last1}.
  The comparisons of Lemma \ref{last1} are natural in the second variable and thus compatible with the transition maps \(F^{a+1}\mathbb B_{\mathrm{pdR}}\to F^a\mathbb B_{\mathrm{pdR}}\); hence the pointwise quasi-isomorphisms assemble into the claimed filtered quasi-isomorphism.

The proof of the Hyodo--Kato case is analogous. 
\end{proof}

\subsection{Conditional multiplicative refinement}\label{soso12} 
We construct an intrinsic convolution algebra structure on the derived Hom appearing in the curve-level reconstruction theorem.
 Conditional on compatibility of the counit in the six-functor formalism with the geometric trace, we show that the reconstruction quasi-isomorphisms are multiplicative.

\subsubsection{Input from the six-functor formalism}
\label{soso10}
We recall here some of  the properties of the six-functor formalism from 
\cite{ALBM} that will be used below\footnote{
As in \cite[Sec.~6.2]{ALBM}, we do not distinguish notationally
between an analytic adic space, or a perfectoid space, and its
associated diamond. All fiber products and all six operations in
Sections \ref{soso10}, \ref{soso1}, and \ref{soso2}  are understood in the category of diamonds. }.

Let \(S\in\mathrm{Perf}_C\)  and assume that \(S^\flat\) admits a
morphism of finite \(\mathrm{dim.trg}\) to a totally disconnected
perfectoid space\footnote{For our applications $S$ strictly totally disconnected is enough.}; we will  write \(S\) also for the associated diamond.
The category
$
\mathcal D_{\mathrm{FF}}(S,\Q_p)
$
from  \cite[Def.~4.2.5]{ALBM} identifies  with
the category $\mathcal D_{\widehat{\square}}
(X_{\mathrm{FF},S^\flat})$ of modified solid quasi-coherent complexes on the
relative Fargues--Fontaine curve \(X_{\mathrm{FF},S^\flat}\);
see \cite[Rem.~4.2.6(i)]{ALBM}. The comparison functors
\[
\alpha^*\colon\quad
\mathrm{QCoh}(X_{\mathrm{FF},S^\flat})
\rightleftarrows
\mathcal D_{\widehat{\square}}
   (X_{\mathrm{FF},S^\flat})
\quad\colon \alpha_*
\]
induce a symmetric monoidal equivalence on nuclear objects by
\cite[Prop.~2.1.2(ii),  Rem.~2.2.2]{ALBM}. We denote the resulting
equivalence by
\begin{equation*}
\Phi_S\colon\quad
\mathcal D_{\mathrm{FF}}(S,\Q_p)^{\mathrm{nuc}}
\xrightarrow{\ \sim\ }
\mathrm{Nuc}(X_{\mathrm{FF},S^\flat}).
\end{equation*}
This equivalence is functorial with respect to pullback in \(S\).

We will use  the following properties of
\(\mathcal D_{\mathrm{FF}}(-,\Q_p)\):

\begin{enumerate}
\item[(1)] ({\em Projection formula}) 
For every morphism \(g\), there is an adjunction
$
\LL g^*\dashv \R g_*.
$
If \(g\) is \(!\)-able, there is moreover an adjunction
$
\R g_!\dashv g^!,
$
and the usual compatibilities with composition, base change, and the
projection formula hold. In particular, there is a natural
quasi-isomorphism
\[
\R g_!\bigl(M\otimes^{\LL}\LL g^*N\bigr)
\simeq
\R g_!M\otimes^{\LL} N.
\]
For a proper lpbc\footnote{Here  lpbc means locally
\(p\)-bounded and compactifiable in the sense of
\cite[Def.~3.2.3]{ALBM}.}  morphism \(g\), one has \(\R g_!\simeq \R g_*\). See \cite[Th.~4.3.3, Rem.~4.3.12]{ALBM}.
 \item[(2)] (\emph{\(!\)-base change})
Consider a Cartesian square of small \(v\)-stacks
\[
\begin{tikzcd}
Y' \ar[r,"g'"] \ar[d,"f'"'] &
Y \ar[d,"f"] \\
X' \ar[r,"g"'] &
X .
\end{tikzcd}
\]
If \(f\) is \(!\)-able, then \(f'\) is \(!\)-able and there is a 
natural base-change quasi-isomorphism
\[
\LL g^*\R f_!
\simeq
\R f'_!\LL g'^*.
\]
See \cite[Th.~4.3.3, Rem.~4.3.12]{ALBM}.
\item[(3)] ({\em Purity}) 
If \(g\colon Y\to Z\) is a  smooth equidimensional morphism of
analytic adic spaces over \(\Q_p\), of relative dimension \(e\), then
\(g\) is cohomologically smooth and there is a natural purity
quasi-isomorphism
\begin{equation*}
g^!M\simeq \LL g^*M(e)[2e].
\label{ALBM-purity}
\end{equation*}
In particular, we have $g^!\mathbf 1_Z\simeq\mathbf 1_Y(e)[2e].$
See \cite[Th.~4.4.4, Th.~4.4.6]{ALBM}.

\item[(4)]({\em Poincar\'e duality and its monoidal refinement}) 
Let \(p\colon Y\to T\) be a cohomologically smooth \(!\)-able morphism of small \(v\)-stacks.  Poincar\'e duality gives a natural
quasi-isomorphism
\begin{equation}\label{ogrod1}
\R p_*\mathbf 1_Y
\simeq
\R\mathcal{H}om_T
   \bigl(\R p_!p^!\mathbf 1_T,\mathbf 1_T\bigr).
\end{equation}
The six-functor  formalism canonically
endows \(\R p_*\mathbf 1_Y\) with a commutative algebra structure,
\(\R p_!p^!\mathbf 1_T\) with a cocommutative counital coalgebra structure, and
$
\R \Hhom_T(\mathcal \R p_!p^!\mathbf 1_T,\mathbf 1_T)
$
with the  induced convolution commutative algebra structure. The 
Poincar\'e duality morphism \eqref{ogrod1}
is a morphism of commutative
algebra objects. See
\cite[Def.~A.0.1, Th.~4.3.3,  Th.~6.2.1]{ALBM}.
\item[(5)]({\em Compactly supported  pro-\'etale cohomology sheaves})
Let \(f\colon X\to\Spa(C)\) be partially proper rigid analytic space. Choose an exhaustive strict covering $\{Y_n\}$, $n\in\N$,  of $X$ by quasi-compact open subspaces and let
$
f_n\colon Y_n\to\Spa(C),
h_n\colon X\setminus Y_n\to\Spa(C) $
be the corresponding structural morphisms. Excision gives
\begin{equation}
\R f_!\mathbf 1_X
\simeq
\underset{n}{\operatorname{colim}}\,\R f_{n,!}\mathbf 1_{Y_n}
\simeq
\operatorname{fib}\left(
\R f_*\mathbf 1_X
\longrightarrow
\underset{n}{\operatorname{colim}}\,
\R h_{n,*}\mathbf 1_{X\setminus Y_n}
\right).
\label{ALBM-excision}
\end{equation}
See \cite[Prop.~6.2.5]{ALBM}. Moreover, if \(X\) is smooth, then
\(\R f_!\mathbf 1_X\) is nuclear; see
\cite[Prop.~6.2.10]{ALBM}.

\item[(6)] ({\em Pro-\'etale cohomology sheaves})  Let $Z$ be an  analytic adic space with a structure
map 
\(g\colon Z\to S\) that is  \(p\)-bounded. Let
\(I\subset(0,\infty)\) be a compact interval with rational
endpoints.  Denote by \(M_I\) the \(\B_I(S)\)-module obtained from
\(M\in\mathcal D_{\mathrm{FF}}(S,\Q_p)\) by pulling back to the
\(Y_{\FF,S^{\flat}}\)-curve and restricting to the annulus
\(Y_{\FF,S^\flat,I}\). Specializing
\cite[Rem.~6.2.8]{ALBM} to \(N=\B_I(S)\), and  using
the \(\mathcal D_{\mathrm{nuc}}(\Q_p)\)-linearity of \(\R g_*\) from
\cite[Prop.~4.2.7(i)]{ALBM},  we get\footnote{We write "pro-\'etale topology"  for "quasi-pro-\'etale topology".} a natural quasi-isomorphism  of solid \(\B_I(S)\)-complexes
\begin{equation}
\bigl(\R g_*\mathbf 1_Z\bigr)_I
\simeq
\R\Gamma_{\proeet}\bigl(Z,\mathbb {B}_I\bigr),
\label{ALBM-annulus}
\end{equation}
where we wrote $
(\R g_*\mathbf 1_Z\bigr)_I:=
\R\Gamma\!\left(
Y_{\mathrm{FF},S^\flat,I},
\LL\pi_S^*\R g_*\mathbf 1_Z
\right).$
Moreover, \(\R g_*\mathbf1_Z\) is nuclear
and its formation commutes with base change in \(S\).

More generally,  combined with \cite[Lemma 6.2.6 ]{ALBM}, the above yields the same description, together with nuclearity and base change, for quasi-separated spaces admitting a countable analytic exhaustion by \(p\)-bounded pieces.
\end{enumerate}

\subsubsection{Example} \label{soso1}
Let $X$ be a smooth partially proper rigid analytic space over $K$ of pure dimension $d$,  let $S$ be a strictly totally disconnected affinoid perfectoid over $C$, and denote by $p: X_S\to S$ the structural morphism. The morphism \(p\) is cohomologically
smooth and \(!\)-able, so we may apply point~{\rm (4)} of
Section~\ref{soso10} to $
T=S, Y=X_S $. 
Let
\[
q\colon\quad  X_S\times_S X_S\longrightarrow S,
\qquad
\Delta\colon\quad  X_S\longrightarrow X_S\times_S X_S
\]
be the structural morphism and the diagonal. Since
\(p=q\Delta\), the counit
$
\R\Delta_{!}\Delta^!\longrightarrow\mathrm{Id}
$
induces the binary coproduct 
\[
\begin{aligned}
\delta_{X_S} \colon\quad  \mathcal C_S(X_S)
&\simeq
\R q_!\R\Delta_!\Delta^!q^!\mathbf 1_S
\to 
\R q_!q^!\mathbf 1_S=:C_S(X_S\times_SX_S)\\ 
 & \xrightarrow[\sim]{\ \kappa_!^{-1}\ }
\mathcal C_S(X_S)\otimes^{\LL} \mathcal C_S(X_S),
\end{aligned}
\]
where $\kappa_!: C_S(X_S)\otimes^{\LL}C_S(X_S)\to C_S(X_S\times_SX_S)$ is the K\"unneth map and $C_S(X_S):=\R p_{!}p^!\mathbf 1_S$. 
The coalgebra counit is the counit of the adjunction
\(\R p_{!}\dashv p^!\):
\[
\epsilon_{X_S}\colon\quad 
C_S(X_S)=\R p_{!}p^!\mathbf 1_S
\longrightarrow
\mathbf 1_S.
\]

 These morphisms
are the binary coproduct and counit of the canonical cocommutative
counital coalgebra structure on \(C_S(X_S)\). All higher
coherences are supplied by the correspondence formalism. This
construction is compatible with pullback in \(S\).
\begin{remark}
Here and below, all base change compatibilities with respect to $S$  are understood with respect
to \(K\)-linear morphisms \(u\colon S'\to S\). 
\end{remark}

\subsubsection{Comparison of cohomology sheaves} \label{soso2}
We now compare the objects discussed in Section \ref{soso10}  with the solid quasi-coherent
cohomology sheaves used in this paper. Let $X$  and $S$ be as in Section \ref{soso1}. 
The annular description \eqref{ALBM-annulus} gives the quasi-isomorphism
\begin{equation}
\Phi_S\bigl(\R p_{*}\mathbf 1_{X_S}\bigr)
\simeq
\mathcal E_{\mathrm{pro\acute et}}(X_S,\Q_p).
\label{tak1}
\end{equation}
Under this  comparison  the  canonical commutative algebra structure on 
\(\R p_{*}\mathbf1_{X_S}\) is carried to the 
usual pro-\'etale cup-product structure. 

Combining \eqref{ALBM-annulus} with the relative version of the excision formula
\eqref{ALBM-excision} (use \(!\)-base change) gives similarly a quasi-isomorphism
\begin{equation}
\Phi_S\bigl(\R p_{!}\mathbf 1_{X_S}\bigr)
\simeq
\mathcal E_{\mathrm{pro\acute et},c}(X_S,\Q_p).
\label{tak2}
\end{equation}
Indeed, after restriction to \(Y_{\FF,I,S^\flat}\), both sides of
\eqref{tak2} are represented by
$
\R\Gamma_{\mathrm{pro\acute et},c}(X_S,\mathbb B_I).
$
The identifications \eqref{tak1} and
\eqref{tak2} are compatible with base change in \(S\).

 Since  the functor $\Phi_S$ is  symmetric monoidal,  the comparison \eqref{tak2} yields a
coalgebra structure on $\se_{\proeet,c}(X_S,\Q_p(d))[2d]$ and, varying $S$, 
on $
\mathcal E_{\mathrm{pro\acute et},c}
(X_C,\Q_p(d))[2d]
$
in \(\mathrm{Nuc}(X_{\mathrm{FF}})^{G_K}\).

\subsubsection{Multiplicative form of Conjecture  $C_{\st}$}  
We endow the augmented filtered derived category
$
\sff^{\mathrm{aug}}\sd(K_{\square})
=
\operatorname{Fun}\bigl(
(\mathbf Z\cup\{-\infty\})^{\mathrm{op}},
\sd(K_{\square})
\bigr) $
with the symmetric monoidal structure: 
\begin{align*}
\bigl(F^\bullet A\otimes_K^{\LL_{\Box}}F^\bullet B\bigr)^n
& :=
\operatorname*{colim}_{r+s\geq n}
\bigl(F^rA\otimes_K^{\LL_{\Box}}F^sB\bigr),\quad n\in\Z, \\
\bigl(F^\bullet A\otimes_K^{\LL_{\Box}}F^\bullet B\bigr)^{-\infty}
&  :=
A\otimes_K^{\LL_{\Box}}B.
\end{align*}

Let $X$ be as in Section \ref{soso1}. 
In Section \ref{soso2},  we have 
endowed $\mathcal C_X:=
\mathcal E_{\mathrm{pro\acute et},c}
   (X_C,\Q_p(d))[2d] $
with a natural  structure of a  coalgebra object in
\(\mathrm{Nuc}(X_{\mathrm{FF}})^{G_K}\).
For every $G_K$-equivariant augmented filtered commutative
\(\mathcal O\)-algebra \(P^\bullet\), the object $
\R\Hhom_{\mathrm{FF}}(\mathcal C_X,P^\bullet) $
   carries  canonical convolution structure as an augmented
filtered commutative algebra;  applying equivariant global sections to the internal $\Hom$  gives the algebra structure on
$\R\uHom_{\FF,G_K}(\mathcal C_X,P^{\bullet})$.

For the rest of Section \ref{soso12} we will assume the following:
\begin{hypothesis}[Trace compatibility]
\label{trace1}

Under the equivalence
$\Phi_S $ the image of the counit $
\R p_{!}p^!\mathbf 1_S\hskip.5mm{\to}\hskip.5mm\mathbf 1_S
$
agrees,  functorially in $X$, with the geometric trace morphism
$$
\operatorname{tr}_{X,S}\colon\quad 
\mathcal E_{\mathrm{pro\acute et},c}
(X_S,\Q_p(d))[2d]
\to
\mathcal O.
$$
We assume that these identifications are compatible with the
\(K\)-linear base changes in \(S\).
\end{hypothesis}

\begin{theorem}[Multiplicative form of Conjecture  \(C_{\mathrm{st}}\)]
\label{mult-cst}
Let \(X\) be a smooth partially proper rigid analytic variety over
\(K\) of pure dimension \(d\). Assume Hypothesis \ref{trace1}.

\begin{enumerate}
\item
The quasi-isomorphism of Theorem \ref{essen3}   canonically refines
to a quasi-isomorphism  in
$
\operatorname{CAlg}
\bigl(\sff^{\mathrm{aug}}\sd(K_{\square})\bigr):
$
\[
F^\bullet \R\Gamma_{\mathrm{dR}}(X)
\xrightarrow{\ \sim\ }
\R{\uHom}_{\FF,G_K}
\left(
\se_{\mathrm{pro\acute et},c}
\bigl(X_C,\mathbf Q_p(d)\bigr)[2d],
i_{\infty,*}F^\bullet\B_{\mathrm{pdR}}
\right).
\]
\item
Assume that the de Rham cohomology of \(X\) has finite rank. The
quasi-isomorphism of Theorem \ref{luminy20}  canonically refines to a quasi-isomorphism
of commutative algebra objects
\[
\R\Gamma_{\mathrm{HK}}(X_C)
\xrightarrow{\ \sim\ }
\R{\uHom}_{\FF,G_K-\mathrm{sm}}
\left(
\se_{\mathrm{pro\acute et},c}
\bigl(X_C,\mathbf Q_p(d)\bigr)[2d],
\mathcal B_{\mathrm{pFF}}
\right)
\]
in \(\sd^{\prime}_{\varphi,N,G_K}(C^{\mathrm{nr}}_{\square})\). In particular,
\(\varphi\) and \(G_K\) act by algebra automorphisms and \(N\) acts as
a derivation.
\end{enumerate}
\end{theorem}

\begin{proof}  ($\bullet$) {\em De Rham cohomology.} For the de Rham claim, we first show that if \(X\) is Stein and has finite-rank de Rham
cohomology  then 
the morphism  (see \eqref{gateway10})
\begin{equation}\label{gateway11}
  \alpha^{\bullet}_{\dr,X}:\quad  i_{\infty,*}F^{\bullet}(\R\Gamma_{\dr}(X)\otimes_K^{\LL_{\Box}}\B_{\dr})\simeq \R\Hhom_{\FF}(\mathcal C_X,i_{\infty,*}F^{\bullet}\B_{\dr})
\end{equation}
is a quasi-isomorphism of augmented filtered commutative algebra objects.

 For that,  set
$ \mathcal A_X
 :=
 \mathcal E_{\mathrm{pro\acute et}}(X_C,\mathbf Q_p), $ $\mathcal P^\bullet:=i_{\infty,*}F^{\bullet}\B_{\dr}$.
The canonical pairing
$
 \mathcal A_X
 \otimes_{\mathcal O}^{\mathrm L}
 \mathcal C_X
 \longrightarrow
 \mathcal O,
$
after tensoring with \(\mathcal P^\bullet\) and applying the
tensor--Hom adjunction,  induces a canonical morphism of augmented filtered commutative algebra objects (use monoidal Poincar\'e duality  and Hypothesis \ref{trace1})
\begin{equation*} 
 \Theta_{\dr,X}^\bullet:\quad
 \mathcal A_X
 \otimes_{\mathcal O}^{\mathrm L}
 \mathcal P^\bullet
 \to
 \R\mathcal Hom_{{\mathrm{FF}}}
 \bigl(\mathcal C_X,\mathcal P^\bullet\bigr).
\end{equation*}
The multiplicative period comparison quasi-isomorphism (see Theorem \ref{hot1K}) yields a quasi-isomorphism
\begin{equation*} \label{mult-period}
\beta^{\bullet}_{\dr,X}:  \quad \mathcal A_X
 \otimes_{\mathcal O}^{\mathrm L}
 \mathcal P^\bullet
 \xrightarrow{\ \sim\ }
 i_{\infty,*}
 F^\bullet\bigl(
   \R\Gamma_{\mathrm{dR}}(X)
   \otimes_K^{\LL_{\Box}}\B_{\mathrm{dR}}
 \bigr)
\end{equation*}
of augmented filtered commutative algebra objects.

 By construction of the morphism $\alpha_{\mathrm{dR},X}^\bullet$ and compatibility of
the period comparison with products and traces, we have 
$$ \Theta^\bullet_{\dr, X} = \alpha^\bullet_{\mathrm{dR},X}\circ\beta^\bullet_{\mathrm{dR},X}. $$
Hence, the above  observations  show that \(\alpha^\bullet_{\mathrm{dR},X}\) is a quasi-isomorphism of commutative algebra objects.

    Attaching  $\log t$ to both domain and target of \(\alpha^\bullet_{\mathrm{dR},X}\) , then applying the lax symmetric monoidal functor of equivariant global
sections,  passing (multiplicatively) to $\B_{\pdr}$, and  using the
multiplicative form of the filtered quasi-isomorphism
\[
\R\Gamma\left(
G_K,
F^\bullet\bigl(
\R\Gamma_{\mathrm{dR}}(X)
 \otimes_K^{\LL_{\Box}}\B_{\mathrm{pdR}}
\bigr)
\right)
\xrightarrow{\ \sim\ }
F^\bullet \R\Gamma_{\mathrm{dR}}(X)
\]
appearing in the proof of Theorem~\ref{essen3}(2) (as $a$ varies), gives the first
claim for Stein \(X\) with finite-rank de Rham cohomology.
 The general case follows by the same analytic hypercover
argument as in the proof of Theorem~\ref{essen3}  since local algebra equivalences are compatible with all cosimplicial transition maps and  limits of
commutative algebra objects are computed on the underlying objects.

  ($\bullet$) {\em Hyodo--Kato  cohomology.} We pass now  to the Hyodo--Kato claim, where we argue in a very similar manner. Set
$$
\mathcal H_X :=
\operatorname*{colim}_{n,m,L}
\R\pi_*
\left(
\R\Gamma_{\mathrm{HK}}(X_C)
\otimes_{C^{\mathrm{nr}}}^{\LL_{\Box}}
\mathcal O_{n,m,L}^+
\right).
$$
The Hyodo--Kato cup product, together with the multiplication on $\so_{\pFF}$
endows \(\mathcal H_X\) with the structure of a commutative algebra
object. This structure is compatible with \(\varphi\), \(N\), and
the \(G_K\)-action.

Tensoring  the multiplicative ordinary basic comparison triangle on the Fargues-Fontaine curve with \(\mathcal B_{\mathrm{pFF}}\)  we see that the de Rham term disappears (because it is \(t\)-torsion and \(t\) is invertible on the coefficients)  yielding a natural quasi-isomorphism
$$
\beta_{\mathrm{HK},X}\colon\quad 
\mathcal A_X\otimes_{\mathcal O}^L\mathcal B_{\mathrm{pFF}}
\xrightarrow{\ \sim\ }
\mathcal H_X
$$
of commutative algebra objects, compatible with
\(\varphi\), \(N\), and \(G_K\).
On the other hand, the  pairing
$
\mathcal A_X\otimes_{\mathcal O}^L \mathcal C_X\longrightarrow\mathcal O
$
induces, after tensoring with
\(\mathcal B_{\mathrm{pFF}}\) and applying the tensor--Hom
adjunction, a canonical morphism 
$$
\Theta_{\mathrm{HK},X}\colon\quad 
\mathcal A_X\otimes_{\mathcal O}^L\mathcal B_{\mathrm{pFF}}
\longrightarrow
\R\Hhom_{\mathrm{FF}}
   \left(\mathcal C_X,\mathcal B_{\mathrm{pFF}}\right).
$$
This is a morphism of commutative
algebra objects, where the right-hand side is endowed with the
convolution multiplication induced by the coalgebra structure on
\(\mathcal C_X\).

Let 
\begin{equation}
\alpha_{\mathrm{HK},X}\colon\quad 
\mathcal H_X
\xrightarrow{\ \sim\ }
\R\Hhom_{\mathrm{FF}}
   \left(\mathcal C_X,\mathcal B_{\mathrm{pFF}}\right)
\label{oldish2}
\end{equation}
denote the quasi-isomorphism of
Theorem~\ref{luminy20}(1), written in this direction. By construction of \(\alpha_{\mathrm{HK},X}\) and compatibility of the period comparison with products and traces, we have 
$$ \Theta_{\hk, X} = \alpha_{\mathrm{HK},X}\circ\beta_{\mathrm{HK},X}. $$
Hence, the above  observations  show that \(\alpha_{\mathrm{HK},X}\) is a quasi-isomorphism of commutative algebra objects.

 Finally, the lax symmetric monoidal smooth equivariant global-sections functor applied to $\alpha_{\mathrm{HK},X}$ and the multiplicative acyclicity of $\B_{\pFF}$  yield our claim  (see the proof of Theorem \ref{luminy20}). 
\end{proof}

\subsubsection{Multiplicative Hodge reconstruction}
We finish this section with the following result: 
\begin{corollary}[Multiplicative Hodge reconstruction]
\label{mult-hodge} Let \(X\) be a smooth partially proper rigid analytic variety over
\(K\) of pure dimension \(d\). Set $ \mathcal C_X:=\se_{\mathrm{\proeet},c}(X_C,\mathbf Q_p(d))[2d]$. 
Taking associated graded objects in
Theorem~\ref{mult-cst}(1) gives a quasi-isomorphism of graded commutative
algebras
\begin{equation}\label{oldish1}
 \operatorname{gr}_F^\bullet \R\Gamma_{\mathrm{dR}}(X)
 \xrightarrow{\ \sim\ }
 \bigoplus_{a\in\mathbf Z}
 \R{\uHom}_{\FF,G_K}
 \left(\mathcal C_X,i_{\infty,*} C(a)[\log t]
 \right).
\end{equation}

 If \(X\) is Stein, this becomes
\[
 \bigl(\Omega^\bullet(X)[-\bullet],d,\wedge\bigr)
 \xrightarrow{\ \sim\ }
 \bigoplus_{a=0}^d
 \R{\uHom}_{\FF,G_K}
 \left(\mathcal C_X, i_{\infty,*} C(a)[\log t]
 \right)
\]
as differential graded commutative algebra objects (the differentials on the right hand side are induced by applying  $\R{\uHom}_{\FF,G_K}(\mathcal C_X,i_{\infty,*}(-))$ to \eqref{nie1}).
\end{corollary}
\begin{proof}We will apply 
the symmetric monoidal associated graded functor
$
 \operatorname{gr}_F^\bullet:
 \sff^{\mathrm{aug}}\sd(K_{\square})
 \longrightarrow
 \operatorname{Gr}\sd(K_{\square})
$
 to Theorem~\ref{mult-cst}(1) to give us
a quasi-isomorphism  of graded commutative algebra objects.
If \(X\) is Stein, coherent acyclicity gives
\[
 \operatorname{gr}_F^a\R\Gamma_{\mathrm{dR}}(X)
 \simeq
 \R\Gamma(X,\Omega_X^a)[-a]
 \simeq
 \Omega^a(X)[-a].
\]
Under this identification
and under the  quasi-isomorphisms \eqref{oldish1}, the convolution multiplication on  the right hand side of the latter  corresponds to exterior multiplication. 
The differentials were identified in Corollary \ref{gateway1}.
\end{proof}
 \section{Duality theorems revisited}
Having the results of this paper we revisit  in this section the proofs of the geometric and arithmetic dualities for $p$-adic pro-\'etale cohomology of  rigid analytic varieties from \cite{CGN2} and \cite{Li}.

\subsection{Geometric duality}
The geometric level computations in this paper can be thought of as a refinement of those  appearing in the proof of geometric duality for $p$-adic pro-\'etale cohomology by Colmez-Gilles-Nizio{\l} \cite{CGN2} (see \cite{ALBM} for a different proof). Hence it is no surprise that the latter can be derived from the former. We present here a simple  argument. 

Recall the statement of geometric duality: If $X$ is a smooth partially proper rigid analytic variety over $K$ of dimension $d$ the natural pairing in TVS's
$$ \mathbb{R}_{\proeet}(X_C,\Q_p) \otimes_{\Q_p}^{\LL_{\Box}}\mathbb{R}_{\proeet,c}(X_C,\Q_p(d)) \to \mathbb{R}_{\proeet,c}(X_C,\Q_p(d))\xrightarrow[]{\mathrm{tr}_{X_C}}\Q_p[-2d]
$$
induces a duality quasi-isomorphism
\begin{equation}\label{kiln1}
\gamma_{X_C}:\quad  \mathbb{R}_{\proeet}(X_C,\Q_p)\stackrel{\sim}{\to}  \R\Hhom_{\mathrm{TVS}}(\mathbb{R}_{\proeet,c}(X_C,\Q_p(d))[2d],\Q_p)
\end{equation}

Here is how we can prove it using  the results from this paper.  For $r\geq 0$, we have the fundamental exact sequence of TVS's 
$$
0\to \Q_p(r)\to \mathbb{B}_{\FF}^{N=0,\phi=p^r}\to \mathbb{B}_{\dr}/F^r\to 0
$$
Take $r\geq 2d$ assume that the de Rham cohomology of $X$, Stein,  is of finite rank. Consider the  following diagram in which the top and bottom rows are   distinguished triangles: (We set $\mathbb{R}_{\proeet,c}:=\mathbb{R}_{\proeet,c}(X_C,\Q_p(d))[2d]$.) 
$$
\xymatrix@R=5mm@C=4mm{
\R\Hhom_{\mathrm{TVS}}(\mathbb{R}_{\proeet,c},\Q_p(r))\ar[r] &   [\R\Hhom_{\mathrm{TVS}}(\mathbb{R}_{\proeet,c}, \mathbb{B}_{\FF})]^{N=0,\phi=p^r}\ar[r]^-{\iota} &   \R\Hhom_{\mathrm{TVS}}(\mathbb{R}_{\proeet,c},\mathbb{B}_{\dr}/F^r)\\
\mathbb{R}_{\proeet}(X_C,\Q_p(r)) \ar[r] \ar[u]_{\gamma_{X_C}}& [\R\Gamma_{\hk}(X_C)\otimes^{\LL_{\Box}}_{C^{\nr}}\mathbb{B}_{\FF}]^{N=0,\phi=p^r}\ar[r]^-{\iota_{\hk}\otimes\iota} \ar[r] \ar[u]^{\wr}_{\alpha_{\hk}}& (\R\Gamma_{\dr}(X)\otimes_K^{\LL_{\Box}}\mathbb{B}_{\dr})/F^r\ar[u]^{\wr}_{\alpha_{\dr}}
}
$$
The bottom triangle is the $t$-inverted form of  the basic comparison triangle. The middle and the  right vertical arrows are 
defined in the following way. Define the map 
$$\alpha_{\hk}:\quad\R\Gamma_{\hk}(X_C)\otimes^{\LL_{\Box}}_{C^{\nr}}\mathbb{B}_{\FF}\to\R\Hhom_{\mathrm{TVS}}(\mathbb{R}_{\proeet,c},\mathbb{B}_{\FF})$$
as the adjoint of the Hyodo--Kato  pairing, using the map from  $\mathbb{R}_{\proeet,c}$ to  compactly supported Hyodo--Kato cohomology\footnote{There are some twists involved in doing that so one needs to be careful.}. Similarly, define the map 
$$\alpha_{\dr}:\quad \R\Gamma_{\dr}(X)\otimes^{\LL_{\Box}}_{K}\mathbb{B}_{\dr}\to \R\Hhom_{\mathrm{TVS}}(\mathbb{R}_{\proeet,c},\mathbb{B}_{\dr})$$
from the de Rham Poincar\'e pairing.  These maps are  quasi-isomorphisms via computations  \eqref{dzisiaj2} and \eqref{dzisiaj1}, respectively\footnote{More precisely, their TVS-version with $\mathbb B_{\FF}$ and filtered $\mathbb B_{\dr}$ coefficients.}.

  We claim that the above  diagram is a map of distinguished triangles. Indeed, the right square  commutes:  After adjunction, its commutativity follows from the fact  that the two maps
$$
(\R\Gamma_{\hk}(X_C)\otimes_{C^{\nr}}^{\LL_{\Box}}\mathbb{B}_{\FF})\otimes^{\LL_{\Box}}_{\Q_p} \mathbb{R}_{\proeet,c}\to \mathbb{B}_{\dr}
$$
are canonically homotopic. One route first applies the Hyodo--Kato pairing,  then the period morphism, and, finally, the projection to Hyodo--Kato cohomology; the other first applies the Hyodo--Kato morphism  and then the de Rham pairing.  These are
 homotopic by the familiar compatibilities (see \cite{CGN2} for a discussion):   the pro-\'etale pairing  is compatible with the syntomic pairing, 
the Hyodo--Kato and de Rham cup pairings  are compatible with each other  and  the syntomic pairing. Applying the functorial fiber construction yields a morphism of fiber sequences in the diagram.  Via  the same compatibilities, we see that, under the tensor-Hom adjunction, the induced morphism on the fibers 
is adjoint to the pro-\'etale  pairing, hence it is  $\gamma_{X_C}$, as wanted. 

   Now, it follows from the above diagram that the map  $\gamma_{X_C}$ is a quasi-isomorphism yielding another proof of the duality quasi-isomorphism \eqref{kiln1} in this case. The general case is obtained by passing to colimit in $\mathbb{R}_{\proeet,c}$ and limit in $\mathbb{R}_{\proeet}(X_C,\Q_p(r))$ over suitable analytic coverings with finite rank de Rham cohomologies. 

\subsection{Arithmetic duality}
Recall the statement of arithmetic  duality from \cite{Li}: If $X$ is a smooth partially proper rigid analytic variety over $K$ of dimension $d$  and $j\in\Z$ the natural pairing in $\sd(\Q_{p,\Box})$
$$ {\R}\Gamma_{\proeet}(X,\Q_p(j)) \otimes_{\Q_p}^{\LL_{\Box}}{\R}\Gamma_{\proeet,c}(X,\Q_p(d+1-j)) \to {\R}\Gamma_{\proeet,c}(X,\Q_p(d+1))\xrightarrow[]{\mathrm{tr}_{X}}\Q_p[-2d-2]
$$
induces a duality quasi-isomorphism 
\begin{equation}\label{kiln11}
\gamma_{X}:\quad  {\R}\Gamma_{\proeet}(X,\Q_p(j))\stackrel{\sim}{\to}  \R\uHom_{\Q_{p}}({\R}\Gamma_{\proeet,c}(X,\Q_p(d+1-j))[2d+2],\Q_p)
\end{equation}

Our formalism provides an intermediate  form of arithmetic duality. Namely, for $j\in\Z$, the above computations yield the duality quasi-isomorphism in ATVS's
$$
\gamma_{X_C}:\quad  \mathbb{R}_{\proeet}(X_C,\Q_p(j))\stackrel{\sim}{\to}  \R\Hhom_{\mathrm{ATVS}}(\mathbb{R}_{\proeet,c}(X_C,\Q_p(d-j))[2d],\Q_p). 
$$
Applying equivariant global sections yields a quasi-isomorphism in $\sd(\Q_{p,\Box})$
$$
\wt{\gamma}_{X}:\quad  \R\Gamma_{\proeet}(X,\Q_p(j))\stackrel{\sim}{\to}  \R\uHom_{\mathrm{ATVS}}(\mathbb{R}_{\proeet,c}(X_C,\Q_p(d-j))[2d],\Q_p). 
$$
Combining it with the  local Tate-duality comparison  from Lemma \ref{kicius1} below recovers the solid-dual formulation of arithmetic Poincar\'e duality from \eqref{kiln11}.
More specifically, let $\lambda_i$ 
denote the quasi-isomorphism of Lemma \ref{kicius1}, obtained by transporting 
the local Tate-duality morphism of~\cite{Li}
through the  curve-to-ATVS comparison. For $ i=d-j$, consider the composite
$\lambda_{d-j}[-2d]\wt{\gamma}_X$.
Under the tensor-Hom adjunction, this composite is induced by the pro-\'etale cup product, followed by the geometric trace,  and the local Galois trace. Hence it   agrees, up to the normalization of the trace, with the arithmetic duality morphism $\gamma_X$
(see \cite[proof of Th. 4.27 and  Lemma 4.28]{Li}).
\begin{lemma} [Local Tate duality] \label{kicius1} Let $i\in\Z$. We have a quasi-isomorphism in $\sd(\Q_{p,\Box})$
$$
\R\uHom_{\mathrm{ATVS}}(\mathbb{R}_{\proeet,c}(X_C,\Q_p(i)),\Q_p)\simeq D_{\Q_p}(\R\Gamma_{\proeet,c}(X,\Q_p(i+1))[2]),
$$
where we set $D_{\Q_p}(-):=\R\uHom_{\Q_{p}}(-,\Q_p)$. 
 \end{lemma}
 \begin{proof}  This is an ATVS's version of a Fargues-Fontaine curve level result proven 
in \cite[Th. 4.18]{Li}  and will be derived from it by  fully-faithfulness of the global section functor $\R\tau_*$. More specifically, we have the quasi-isomorphisms
 in $\sd(\Q_{p,\Box})$
 \begin{align*}
    \R\uHom_{\mathrm{ATVS}}(\mathbb{R}_{\proeet,c}(X_C,\Q_p(i),\Q_p)) &\simeq  \R\uHom_{\mathrm{ATVS}}(\R\tau_*\se_{\proeet,c}(X_C,\Q_p(i)),\R\tau_*\so)\\
     & \simeq \R\uHom_{\FF,G_K}(\se_{\proeet,c}(X_C,\Q_p(i)),\so)\\
     & \simeq D_{\Q_p}(\R\Gamma_{\proeet,c}(X,\Q_p(i+1))[2]),
    \end{align*}
    where the third quasi-isomorphism is the local Tate duality on the curve and the second one follows from Proposition \ref{film1}.
 \end{proof}
 \begin{remark} The proof of arithmetic duality in~\cite{Li}
uses the geometric duality on the Fargues-Fontaine curve from \cite{CGN2} plus local Tate duality on the curve. Here, we use the geometric duality on the level of TVS's from \cite{CGN2} as well and then a local Tate duality on the level of ATVS's. As both are descended from the curve via the global section functor and its fully-faithfulness one can say that the argument presented here is a descent version of the argument of~\cite{Li}. 
 \end{remark}

 \appendix
\section{Extensions of $C$ by $C$}\label{firstA}
Theorem~\ref{extCC} below
explains why, to formulate the $C_{\rm st}$-conjecture, we need to work in the category 
of Arithmetic Topological Vector Spaces and not, as one would naturally do, in the category
$\Q_{p,\square}[G_K]$ of $G_K$-equivariant solid $\Q_p$-modules.

Let $K$ be a finite extension of $\Q_p$.  We consider $C=\widehat{\Qbar_p}$ as a continuous $\Q_p$-Banach representation 
of $G_K$. 

\subsection{Almost $C$-representations versus Banach representations} 
Recall that Fontaine~\cite{FoC} defined the subcategory ${\cal C}(G_K)$
of $\Q_{p,\square}[G_K]$ 
of almost $C$-representations which, a posteriori, can be thought of  as the subcategory
of ATVS's consisting of finite Dimensional objects.  In this category, extensions of $C$ by $C$ are
almost split, and Fontaine proved that $\Ext^1_{{\cal C}(G_K)}(C,C)\simeq K$,
the non-trivial class being given by tensoring the representation $K\oplus K\log t$ with $C$.
Hence Theorem~\ref{extCC} below  implies\footnote{See Remark \ref{gateway7}.}, in particular,
that the natural functor between bounded derived categories
$$\sd^b({\cal C}(G_K))\to \sd^b(\Q_{p,\square}[G_K])$$
is not fully faithful.

Let $\Ext^1_{G_K}(C,C)$ be the Yoneda extension group in the category of Banach $(\Q_p,G_K)$-representations.
So an element in $\Ext^1_{G_K}(C,C)$ corresponds to a sequence of Banach $\Q_p$-spaces 
$$0\longrightarrow C\longrightarrow E\longrightarrow C\longrightarrow 0$$
which is strictly exact and  $G_K$-equivariant.

\begin{theorem}\label{extCC}
We have\footnote{The $K$-action comes from the second term $C$.}  $\dim_K \Ext^1_{G_K}(C,C)=+\infty$.
\end{theorem}

\begin{remark}\label{gateway7}
The extension groups $\Ext^1_{G_K}(C,C)$,
computed in the categories of locally convex $(\Q_p,G_K)$-representations 
and Banach $(\Q_p,G_K)$-representations, coincide.
Indeed, any extension of two Banach spaces in the category of locally convex topological vector spaces is itself a Banach space.    Moreover, this group clearly injects into its solid analog. 
\end{remark}

\subsection{Topological considerations}
In what follows, {\it group cohomology is continuous group cohomology}.

Let $G$ be a profinite group, and
let $W,V$ be two $\Q_p$-Banach spaces equipped with a continuous $G$-action.
Since $G$ is compact, one can choose $G$-invariant valuations on $W,V$ defining
the Banach topology; then the closed unit balls $W^+$ and $V^+$ are stable by $G$.
Let ${\rm Hom}(W,V)$ denote the space of continuous $\Q_p$-linear morphisms.
We have ${\rm Hom}(W,V)=\Q_p\otimes_{\Z_p}{\rm Hom}_{\Z_p}(W^+,V^+)$
since a $\Q_p$-linear map $\lambda:W\to V$ is continuous if and only if
it is Lipschitz.

Now, ${\rm Hom}_{\Z_p}(W^+,V^+)$ is the Tate module of ${\rm Hom}_{\Z_p}(W^+,V/V^+)$, which
is the projective limit over $n$ of the ${\rm Hom}_{\Z_p}(W^+,p^{-n}V^+/V^+)$. 
And ${\rm Hom}_{\Z_p}(W^+,p^{-n}V^+/V^+)$ is the projective limit of the
${\rm Hom}_{\Z_p}(M,p^{-n}V^+/V^+)$ where $M$ runs through the finite length sub-$\Z_p$-modules
of $W^+/p^n W^+$. Since the action of $G$ on $W^+/p^nW^+$ is smooth for all $n$, and $G$
is profinite, we can
take the limit over subobjects stable by $G$ in which case 
each ${\rm Hom}_{\Z_p}(M,p^{-n}V^+/V^+)$
is discrete and endowed with a smooth action of $G$.  If we endow
${\rm Hom}_{\Z_p}(W^+,p^{-n}V^+/V^+)$ with the projective limit topology, and
${\rm Hom}_{\Z_p}(W^+,V^+)$ also with the projective limit topology,
the action of $G$ remains continuous.
We then endow ${\rm Hom}(W,V)=\colim_n p^{-n}{\rm Hom}_{\Z_p}(W^+,V^+)$ with
the inductive limit topology, which turns it into a locally convex
$\Q_p$-topological vector space with a continuous $G$-action.

\begin{remark}\label{top5}
(i)
The action of $G$ on these Hom spaces is given by 
{\it the diagonal action}
$$(g\cdot\phi)(v)=g\cdot(\phi(g^{-1}\cdot v))$$
and the space ${\rm Hom}_G(W,V)$ of $G$-equivariant continuous morphisms $W\to V$ is
just the fixed part ${\rm Hom}(W,V)^G$ under this action.

(ii)
Another natural topology that can be put on ${\rm Hom}(W,V)$ is {\it the strong topology},
which turns ${\rm Hom}(W,V)$ into a $\Q_p$-Banach, but the action of $G$
on ${\rm Hom}(W,V)$, endowed with this topology, is not continuous in general.

(iii) If $(W_i)_{i\in I}$ is a directed system of finite
dimensional sub-$\Q_p$-vector spaces of $W$, such that $\colim_i W_i$ is dense
in $W$, we also have ${\rm Hom}_{\Z_p}(W^+,V^+)=\lim_i{\rm Hom}_{\Z_p}(W_i^+,V^+)$
as topological $\Z_p$-modules.

(iv) 
Let $\phi:G\to {\rm Hom}(W,V)$ be a continuous $1$-cocycle. Then the map
$$G\times W\to V;\ (g,w)\mapsto \phi(g)(w)$$
is continuous.
 For a fixed $w\in W$, the evaluation map ${\rm Hom}(W,V) \to V;\ f\mapsto f(w)$ is continuous.
Since $G$ is compact,
 the image of $G$ in $V$ under the map $\phi(-)(w)$ is compact, hence bounded.
 This implies the image $\phi(G)$ is bounded in ${\rm Hom}_s(W,V)$ under the pointwise topology.
 By Banach-Steinhaus theorem (\cite[Prop.\,6.15]{schneider}) 
$\phi(G)$ is equicontinuous, which implies the claim: Since $\phi(G)$  is equicontinuous, the maps
$ w\mapsto \phi(g)(w),$ $g\in G$,
are uniformly continuous at $0$;  together with the continuity of $ g\mapsto \phi(g)(w)$,  for fixed $w$, this proves the joint continuity of
$(g,w)\mapsto \phi(g)(w)$.
\end{remark}

\begin{lemma}\label{top2}
{\rm (i)} We have a natural injection:
$$H^1(G,\Hom(W,V))\hookrightarrow \Ext^1_G(W,V)$$

{\rm (ii)} We have a canonical  isomorphism:
$$\Q_p\otimes_{\Z_p} H^1(G,\Hom_{\Z_p}(W^+,V^+))\overset{\sim}{\to}
H^1(G,\Hom(W,V))$$.
\end{lemma}
\begin{proof}
This is standard. The map in (i) is given as follows: 
if $g\mapsto \phi_g$ 
is a continuous 1-cocycle with values in ${\rm Hom}(W,V)$,
we define an action $\rho_\phi$ on $W\oplus V$ by 
$$\rho_\phi(g)=
\begin{pmatrix}
\rho_V(g) & \phi_g\circ\rho_W(g)\\
0 & \rho_W(g)
\end{pmatrix}$$

Now, if $g\mapsto\phi_g$ is a continuous $1$-cocycle, with values in ${\rm Hom}_{\Z_p}(W^+,V^+)$
whose class is $0$ in $H^1(G,{\rm Hom}(W,V))$, there exists $\phi\in {\rm Hom}(W,V)$
with $\phi_g=(g-1)\cdot\phi$ for all $g\in G$.  
Hence the class of $g\mapsto\phi_g$ is killed by $p^N$,  
if $\phi\in p^{-N}{\rm Hom}_{\Z_p}(W^+,V^+)$; 
this implies injectivity in (ii).  For the surjectivity:  By Remark A.3(iv), the compact family $\{\phi_g\}, g\in G$, 
is equicontinuous; hence there exists $N$ such that $\phi_g(W^+)\subset p^{-N}V^+$, 
for all $g\in G$; thus the cocycle takes values in $p^{-N}{\rm Hom}_{\Z_p}(W^+,V^+)$.
\end{proof}

\subsection{Proof of Theorem~\ref{extCC}}
Let $\Delta$ be the torsion subgroup of $\Gal(K(\mu_{p^\infty})/K)$.
Set 
$$K_n=K(\mu_{p^n})^{\Delta},\hskip4mm K_{\infty}=\cup_n K_n,\hskip4mm
\Gamma=\Gal(K_\infty/K)\simeq \Z_p,\hskip4mm H=\Gal(\Qbar_p/{K_\infty})$$
By the Ax-Sen-Tate theorem, we also have $H={\rm Aut}_{\rm cont}(C/\widehat{K_\infty})$ and
$C^H=\widehat{K_\infty}$.

By inflation, we get an injection
$$H^1(\Gamma,\Hom_H(C,C))
\hookrightarrow H^1(G_K,\Hom(C,C))$$
Thus, by Lemma \ref{top2}, {\it it is enough to show that} 
$$\dim_{K} H^1(\Gamma,\Hom_H(C,C))=+\infty$$

Now, by Remark~\ref{top5} (iii),
\begin{equation}\label{gateway5}
\Hom_H(C,C)=\Hom_H(\O_C,{\goth m}_C)[\tfrac{1}{p}]\quad{\rm and}\quad
\Hom_H(\O_C,{\goth m}_C)=\lim
\Hom_H(\O_L,{\goth m}_C),
\end{equation}
where the  limit is over finite Galois extensions of $K$.
Now, if $L_1\subset L_2$ are finite Galois extensions of $K$,
we have an exact sequence $0\to \O_{L_1}\to \O_{L_2}\to X\to 0$ of $\O_K[G_K]$-modules.
This induces an exact sequence
$0\to {\rm Hom}(X,{\goth m}_C)\to {\rm Hom}(\O_{L_2},{\goth m}_C)
\to {\rm Hom}(\O_{L_1},{\goth m}_C)\to 0$ and passing to $H$-invariants
still gives an exact sequence, since $H^1(H,{\rm Hom}(X,{\goth m}_C))=0$
by almost \'etale descent.
(This is the reason for considering ${\goth m}_C$ instead of $\O_C$: for $\O_C$,
the above $H^1$ would only be almost~$0$.) 
Hence all the maps in the  limit from \eqref{gateway5} are surjective,
$\Hom_H(\O_C,{\goth m}_C)\to \Hom_H(\O_{K_n},{\goth m}_C)=
\Hom(\O_{K_n},{\goth m}_{\widehat{K_\infty}})$ is surjective for all $n$, and
we get a surjection
\begin{equation}\label{gateway6}
\Hom_H(C,C)\sur \Hom({K_n},\widehat{K_\infty}),\quad{\text{for all $n$.}}
\end{equation}

Since $\Gamma$ has cohomological dimension $1$, 
this induces a surjection (note that the target of the map in \eqref{gateway6} is finite free over $\widehat{K_\infty}$)
$$H^1(\Gamma,\Hom_H(C,C))\sur H^1(\Gamma,\Hom({K_n},\widehat{K_\infty})),
\quad{\text{for all $n$.}}$$ 
But Tate's normalized trace gives a $\Gamma$-equivariant 
decomposition $\widehat{K_{\infty}}\simeq K_n\oplus X_n$, 
hence we get a surjection
$$H^1(\Gamma,{\rm Hom}({K_{n}},\widehat{K_{\infty}}))\sur
H^1(\Gamma,{\rm Hom}(K_n,K_n)),\quad{\text{for all $n$.}}$$ 
Since $K$ is a finite extension of $\Q_p$, we have that ${\rm Hom}_K(K_n,K_n)\subset {\rm Hom}(K_n,K_n)$ is a $\Gamma$-equivariant  direct summand. Hence 
we get a surjection
$$H^1(\Gamma,{\rm Hom}({K_{n}},\widehat{K_{\infty}}))\sur
H^1(\Gamma,{\rm Hom}_K(K_n,K_n)),\quad{\text{for all $n$.}}$$

Now, $K_n$ is the regular representation of 
$\Gamma_n:={\rm Gal}(K_n/K)$. Hence, after extending scalars
to $\Qbar_p$, it becomes the direct sum of the characters of $\Gamma_n$.
It follows that, if $\gamma$ is a generator of $\Gamma$, the multiplicity
of the eigenvalue $1$ of $\gamma$ acting on ${\rm Hom}(K_n,K_n)$ is $|\Gamma_n|=[K_n:K]$
(each character of
$\Gamma_n$ contributing an eigenvector). Since $H^1(\Gamma,{\rm Hom}(K_n,K_n))=
{\rm Hom}(K_n,K_n)/(\gamma-1)$, we deduce that
$$\dim_{K}H^1(\Gamma,{\rm Hom}(K_n,K_n))=[K_n:K]$$
But it follows from the above discussion that
$H^1(\Gamma,{\rm Hom}(K_n,K_n))$ is a quotient of $H^1(\Gamma,{\rm Hom}_H(C,C))$.
Hence $H^1(\Gamma,{\rm Hom}_H(C,C))$ has quotients of unbounded dimensions, which implies
that it is infinite dimensional, and proves Theorem~\ref{extCC}.

\section{A Fr\'echet version of $\B_{\pdr}$}\label{second1}
In this appendix, we propose a Fr\'echet version of $\B_{\pdr}$ and compute its Galois cohomology.

Retain the notation \(K,\overline K,C,G_K\) fixed above.
Let $\chi$ be the cyclotomic character of $G_K:=\Gal(\overline{K}/K)$.
Let $K\{T\}:=\lim_n K\langle p^nT\rangle\cong \cO(\bA^{1,\an}_K)$; 
this is a $K$-Fr\'echet algebra.
Let $(\B_{\dr}^+/t^n)\{T\}:=(\B_{\dr}^+/t^n)\otimes^{\square}_K K\{T\}$ endowed with $G_K$-action given by natural action on $(\B_{\dr}^+/t^n)$ and $g(T)=T+\log(\chi(g))$ for $g\in G_K$
(i.e., $T=\log t$).
Then $(\B_{\dr}^+/t^n)\{T\}$ is a Fr\'echet $K$-algebra and the action of $G_K$ is continuous.
Let $\widetilde{\B}^+_{\pdr}:=\lim_n ((\B_{\dr}^+/t^n)\{T\})$ which is a Fr\'echet 
$K$-algebra, and
let $\widetilde{\B}_{\pdr}:=\widetilde{\B}^+_{\pdr}[1/t]$.

\begin{theorem}\label{li0}
We have 
$$\R\Gamma(G_K,\widetilde{\B}^+_{\pdr})\simeq K[0],\quad \R\Gamma(G_K,\widetilde{\B}_{\pdr})\simeq K[0].$$
\end{theorem}

\begin{proof}
We have $t^i\widetilde{\B}^+_{\pdr}/t^{j+1}\widetilde{\B}^+_{\pdr}=C(j)\{T\}$;
hence it is enough to show that
\begin{equation}
	\R\Gamma(G_K, C\{T\})=K[0]
\end{equation}
and, for $j\neq 0$,
\begin{equation}
	\R\Gamma(G_K,C(j)\{T\})\simeq 0.
\end{equation}
For $f=\sum_i a_iT^i\in C\{T\}$, let $||\sum_i a_iT^i||_R:=\max\{ |a_i|R^i\}$ be the norm on the disk of radius $R$.

Let $K_n=K(\mu_{p^n})$, $K_\infty=\cup K_n$, $\Gamma_K=\Gal(K_{\infty}/K)$, $\Gamma_n:=\Gal(K_\infty/K_n)$ and $H=\Gal(\overline{K}/K_\infty)$.
 We first choose $m$ such that $\log\chi:\Gamma_m\simeq p^m\Z_p$
 and, if $n\geq m$, we let $\gamma_n$ be the topological generator of $\Gamma_n$ such that
$\log\chi(\gamma_n)=p^n$.
 By finite Galois descent, it suffices to show the above two formulas for $K_m$ and we will thus assume $K=K_m$. Let $\Gamma:=\Gamma_m$. 
 Since $H$ acts trivially on $T$, by almost \'etale descent we have
$$\R\Gamma(G_K,C(j)\{T\})\simeq \R\Gamma(\Gamma,\R\Gamma(H,C(j)\{T\}))
\simeq \R\Gamma(\Gamma,\widehat{K}_\infty(j)\{T\}).$$
For $n\geq m$, Tate's normalized trace $\tr_n:\widehat{K}_\infty\to K_n$ gives a decomposition 
$$\widehat{K}_\infty\cong K_n\oplus X_n,\quad X_n=\ker(\tr_n)$$
 Moreover, $(1-\gamma_n)$ is an isomorphism on $X_n$ and there exists $d$ such that 
$||(\gamma_n-1)^{-1}||_{X_n}\leq p^d$ for all $n$ (\cite[Proposition 1.104]{fontaine2008theory}).
 We choose $n > d$ so that $|\chi(\gamma_n)-1|< p^{-d}$.
 Again, by Galois descent, it suffices to prove the theorem for $K=K_n$,
and the theorem follows from the computation of
$\rg(\Gamma_n,K_n(j)\{T\})$ (Lemma~\ref{li2}) and $\rg(\Gamma_n,X_n(j)\{T\})$ (Lemma~\ref{li3}).
\end{proof}

\begin{lemma}\label{li1}
Let $\alpha: X\{T\}\to X\{T\}$ be a continuous isomorphism of Fr\'echet $K$-spaces 
such that there exist $c$ such that for all $f\in X\{T\}$ 
and $R>R_0$ we have $ ||\alpha^{-1}(f)||_R\leq c||f||_R$.
 Let $\beta:X\{T\}\to X\{T\}$ be $K$-linear, 
 such that, for all $R>R_0$, we have 
$$||\alpha-\beta||_R\leq b<c^{-1}.$$
Then $\beta$ is an isomorphism.
\end{lemma} 
\begin{proof}
We have $\beta=\alpha-(\alpha-\beta)=\alpha(1-(1-\alpha^{-1}\beta))$; hence it is enough
to check that $1-(1-\alpha^{-1}\beta)$ is an isomorphism.
By our assumption, for $R>R_0$, we have $||(1-\alpha^{-1}\beta)f||_R\leq bc||f||_R$ where $bc<1$.
This implies that the series $\sum_{i\geq 0} (1-\alpha^{-1}\beta)^i(f)$
converges in $X\{T\}$ to some $u(f)$, and $u$ is an inverse to $1-(1-\alpha^{-1}\beta)$.
\end{proof}

We write $\gamma$ for $\gamma_n$ and $X$ for $X_n$; we have $\log\chi(\gamma)=p^n$.
\begin{lemma}\label{li3}
For all $j\in\Z$, we have
$$\R\Gamma(\Gamma_n,X(j)\{T\})\simeq 0.$$
\end{lemma}
\begin{proof}
We want to show that 
$$1-\gamma: X\{T\}(j)\to X\{T\}(j)$$
is an isomorphism.
 Let $\gamma_c:X\{T\}(j)\to X\{T\}(j)$; $\sum_i a_iT^i\mapsto \gamma(a_i)T^i$ (here $c$ means $\gamma$ only acts on coefficients).
 Since $(1-\gamma)$ is an isomorphism on $X$, 
$(1-\gamma_c)$ is an isomorphism on $X\{T\}(j)$ with inverse given by $\sum_i a_iT^i\mapsto ((1-\gamma)^{-1}a_i)T^i$.
 For any $R>0$, we have
$$||(1-\gamma_c)^{-1}(a_iT^i)||_R\leq p^d|a_i|R^i\leq p^d||a_iT^i||_R.$$
This implies for all $f\in X\{T\}$ and $R>0$, we have $||(1-\gamma_c)^{-1}f||_R\leq p^d||f||_R$.

Now, if $f=\sum_i a_iT^i$, 
then $(\gamma-\gamma_c)f=\sum_i \gamma(a_i)(\chi(\gamma)^j(T+p^n)^i-T^i)$.
 Since $\gamma$ is an isometry on $X$, we have, if $R>1$,
$$||\gamma(a_i)(\chi(\gamma)^j(T+p^n)^i-T^i)||_R
\leq |a_i|\cdot \max_{0\leq k\leq i-1}\{ |\chi(\gamma)^j-1|R^i, |p^n|^kR^{i-k} \}
\leq p^{-n}|a_i|R^i\leq p^{-n}||f||_R $$
Since $p^{-n}<p^{-d}$, we can use the above lemma (with $c=p^d$, $b=p^{-n}$, $R_0=1$)
to conclude that $\gamma-1$ is an isomorphism.
\end{proof}

\begin{lemma}\label{li2}
We have
\begin{equation}\label{j0}
\R\Gamma(\Gamma_n, K_n\{T\})=K_n[0]	
\end{equation}
and, for $j\neq 0$,
\begin{equation}\label{jn0}
	\R\Gamma(\Gamma_n,K_n(j)\{T\})\simeq 0.
\end{equation}
\end{lemma}
\begin{proof}
Since $\Gamma_n$ acts trivially on $K_n$, we have 
$$(\gamma-1)f=\chi^j(\gamma)f(T+p^n)-f(T)$$
Since $f$ is the sum of its Taylor expansion on all of $C$, 
we have $f(T+p^n)=\exp(p^n\frac{d}{dT})f$.
 Since $\chi^j(\gamma)=\exp( jp^n)$, we get
$$(\gamma-1)f=(\exp(D_j)-1)f,\quad
D_j:= p^n(\tfrac{d}{dT}+j): K_n\{T\}\to K_n\{T\}$$
Set $V_j=\sum_{i\geq 1} D_j^i/(i+1)!$ so that $\exp(D_j)-1=D_j(1+V_j)$.
We have $||D_j(f)||_R\leq p^{-n} ||f||_R$ if $R>1$, from which it follows
that $||V_j(f)||_R\leq p^{-n}||f||_R$ if $R>1$; hence
$1+V_j$ is invertible on $K_n\{T\}$.
We infer that:
$$[K_n(j)\{T\}\stackrel{\gamma-1}\longrightarrow K_n(j)\{T\}]
\simeq [K_n\{T\}\stackrel{D_j}\longrightarrow K_n\{T\}].$$
If $j=0$, \eqref{j0} follows from the computation of de Rham cohomology of $\bA^{1,\an}_{K_n}$.
 If $j\neq 0$, we note that $||j^{-1}\frac{d}{dT}f||_R\leq \frac{1}{R}|j|^{-1}||f||_R
\leq p^{-1}||f||_R$ if $R\geq p|j|^{-1}$.
 This implies that $D_j=p^nj(1+j^{-1}\frac{d}{dT})$ 
is invertible on $K_n\{T\}$, which proves \eqref{jn0}.
\end{proof}

\begin{remark}
The original $\B_{\pdr}^+$ uses $K[T]$ (corresponding to $\bA^{1,\alg}$) instead
of $K\{T\}$ (corresponding to $\bA^{1,\an}$).  One could be tempted to use
open balls instead of $\bA^{1,\an}$ as in the construction of the rings $\B_{{\rm Sen},n}$
(see~\cite{Csen}), but the resulting ring would contain $t^{p^n}\exp (-p^n T)$ for $n$ big enough,
which would produce elements fixed by $G_{K_n}$, and not belonging to $K_n$.
\end{remark}

  \end{document}